%% file: Article/Arxiv/xxx_arxiv.tex
\documentclass[11pt]{article}
\usepackage{amsmath,amssymb,amsthm,bm}
\usepackage{mathtools}
\usepackage{bbm}
\usepackage{dsfont}
\usepackage{verbatim}
\usepackage{graphicx,rotating}
\usepackage{wrapfig,sidecap}
\usepackage{xcolor}
\usepackage{natbib}
\usepackage{ulem} 
\usepackage{enumitem} 
\usepackage{tikz-cd}
\usepackage{anyfontsize}

\usepackage{floatflt}
\usepackage{xfrac}
\usepackage{caption}
\usepackage[refpage]{nomencl}

\def\@@@nomenclature[#1]#2#3#4{%
 \def\@tempa{#2}\def\@tempb{#3}\def\@tempc{#4}%
 \protected@write\@nomenclaturefile{}%
  {\string\nomenclatureentry{#1\nom@verb\@tempa @[{\nom@verb\@tempa}]%
      \begingroup\nom@verb\@tempb\protect\nomeqref{\theequation}%
        |nomlabelref}{\@tempc}}%
 \endgroup
 \@esphack}
\makenomenclature

\newtheorem{theorem}{Theorem}
\newtheorem{lemma}{Lemma}
\newtheorem{corollary}{Corollary}
\newtheorem{example}{Example}

\theoremstyle{definition}
\newtheorem{definition}{Definition}

\theoremstyle{plain}

\theoremstyle{remark}

\newtheorem*{remark}{Remark}

\usepackage{tikz}
\usepackage[framemethod=tikz]{mdframed}
\newmdenv[tikzsetting={draw=black, line width=0.5pt, dash pattern=on 1pt off 1pt, dashed}, linecolor=white, outerlinewidth=1pt]{examplebox}

\usepackage{xr-hyper} 
\begin{document}

\title{Spatial Confidence Regions for Piecewise Continuous Processes} 
\author{Thomas Maullin-Sapey$^1$,  Fabian J.E. Telschow$^2$ \\[5mm]
  $^1$School of Mathematics, Fry Building, University of Bristol\\
	$^2$Department of Mathematics, Humboldt University of Berlin 
}
\date{}
\maketitle

\begin{abstract}
	\input{Article/abstract.tex}

\end{abstract}

\section{Introduction}\label{sec:intro}

\input{Article/introduction.tex}

\section{Background}\label{sec:background}
\input{Article/background.tex}

\section{Preliminaries}\label{sec:preliminaries}
\input{Article/preliminaries.tex}

\section{Motivation}\label{sec:motivation}

\input{Article/motivation.tex}

\section{Convergence Notions}\label{sec:convergence}

\input{Article/convergence.tex}

\section{Restrained Convergence for Random Variables}\label{sec:random}
\input{Article/suprema.tex}

\section{Confidence Regions}\label{sec:piecewise_copes}

\input{Article/crs.tex}

\section{Illustrative Examples}\label{sec:applications}

\input{Article/applications.tex}


\section{Discussion}\label{sec:discussion}

\input{Article/discussion.tex}

\section*{Acknowledgments}
\input{Article/acknowledgement.tex}



\bibliographystyle{plain}
\bibliography{Article/paper-ref}


\end{document}

%% file: Article/abstract.tex
In scientific disciplines such as neuroimaging, climatology, and cosmology it is useful to study the uncertainty of excursion sets of imaging data. While the case of imaging data obtained from a single study condition has already been intensively studied, confidence statements about the intersection, or union, of the excursion sets derived from different subject conditions have only been introduced recently. Such methods aim to model the images from different study conditions as asymptotically Gaussian random processes with differentiable sample paths.

In this work, we remove the restricting condition of differentiability and only require continuity of the sample paths. This allows for a wider range of applications including many settings which cannot be treated with the existing theory. To achieve this, we introduce a novel notion of convergence on piecewise continuous functions over finite partitions. This notion is of interest in its own right, as it implies convergence results for maxima of sequences of piecewise continuous functions over sequences of sets. Generalizing well-known results such as the extended continuous mapping theorem, this novel convergence notion also allows us to construct for the first time confidence regions for mathematically challenging examples such as symmetric differences of excursion sets.

%% file: Article/introduction.tex
Throughout a broad range of statistical settings, researchers wish to assess the same fundamental question; at which spatial locations does an unknown signal exceed a predefined threshold? For instance, in geostatistical disease mapping, malaria prevalence rates which exceed expectation are used as indicators of an outbreak. In cosmological settings, the locations of celestial bodies are identified by assessing where heat signatures exceed a predefined temperature. And in fMRI neuroscience, researchers identify active regions of the brain by assessing where proxy measures of blood flow exceed predefined thresholds.

Such analyses all fall under the same theoretical framework; the study of random excursion sets. Formally, given a spatial domain, $\mathcal{S}$, in each scenario there is an unknown, spatially-varying, target function $\mu:\mathcal{S}\rightarrow \mathbb{R}$ (e.g. malaria prevalence over a geographic region, temperature over space, blood flow across the brain). The aim of the researcher is to estimate \textit{excursion sets} of $\mu$, i.e.~sets of the form $\mathcal{U}_c:=\{s \in \mathcal{S}: \mu(s) > c\}$ or $\mathcal{L}_c:=\{s \in \mathcal{S}: \mu(s) < c\}$, where $c$ is some predefined threshold. The conventional approach to studying such sets is to collect a sample of $n$, potentially noisy, observations of $\mu$, use these to derive an estimator $\hat{\mu}_n:\mathcal{S}\rightarrow \mathbb{R}$ and then report the excursion sets of $\hat{\mu}_n$, defined analogously to the above. This procedure provides estimates of the excursion sets, but no notion of the spatial uncertainty and reliability of such estimates. 

To address these issues, recent literature has introduced the notion of Confidence Regions (CRs). Confidence regions are random sets, denoted $\hat{\mathcal{U}}_c$ and $\hat{\mathcal{L}}_c$, designed to satisfy $\mathbb{P}[\hat{\mathcal{L}}_c \subseteq \mathcal{L}_c~ \wedge ~\hat{\mathcal{U}}_c \subseteq \mathcal{U}_c]=1-\alpha$ asymptotically for some predefined confidence level $\alpha$, e.g. $\alpha=5\%$. An illustration of CRs for a spatially varying signal $\mu:\mathbb{R}^2\rightarrow \mathbb{R}$ and thresholds $c\in\{0,1\}$ is provided by Fig. \ref{fig:crs}.

In applications where $\mu^{-1}(c)$ may be assumed to have zero measure, such regions facilitate the interpretation of results as they can be used to quantify the reliability of the estimated level set $\hat{\mu}_n^{-1}(c)$. Specifically, if the boundaries of the CRs closely resemble one another, this suggests that that the level set has been reliably estimated. Conversely, little resemblance between the CRs indicates either that there is high uncertainty in the estimated level set or the level set has non-zero measure. If the researcher is willing to assume the latter does not hold (that is, $\mu^{-1}(c)$ is a zero-measure set), it can be concluded that the estimated level set is unreliable, and that more data is required. For example, in Fig. \ref{fig:crs} Panel (c), the boundaries of $\hat{\mathcal{U}}_0$ and $\hat{\mathcal{L}}_0$ very closely match one another, accurately capturing the shape of the circular level set $\mu^{-1}(0)$ on the right side of the image. This close resemblance suggests that $\hat{\mu}_n^{-1}(0)$, which will lie in the region enclosed between the CR boundaries, is a reliable estimator of $\mu^{-1}(0)$.
\vspace{0.2cm}

\begin{centering}
\includegraphics[width=0.9\textwidth]{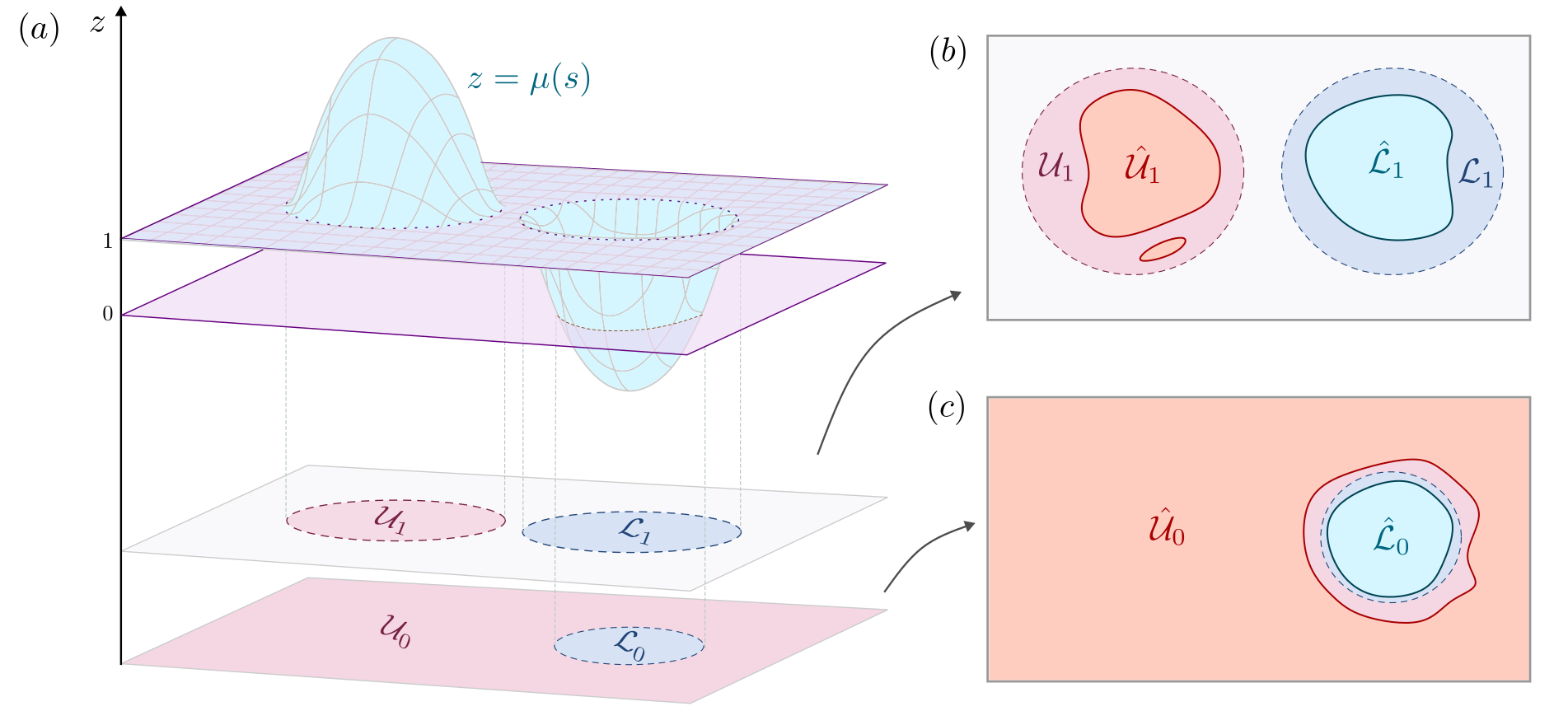}
    \captionof{figure}{(a) A spatially varying signal $\mu:\mathbb{R}^N\rightarrow\mathbb{R}$ (blue) thresholded at the levels $c=0$ and $c=1$ (purple). The sets $\mathcal{U}_c=\{s: \mu(s)>c\}$ and $\mathcal{L}_c=\{s: \mu(s)<c\}$, $c \in \{0,1\}$, are displayed below in dark red and dark blue, respectively. In practical settings, the signal $\mu$ and sets $\mathcal{U}_c$ and $\mathcal{L}_c$ are unknown and of empirical interest. (b) The sets $\mathcal{U}_1$ and $\mathcal{L}_1$, alongside potential CRs $\hat{\mathcal{U}}_1$ and $\hat{\mathcal{L}}_1$. In this case, since the level set $\{s: \mu(s)=1\}$ is a plateau, the boundaries of $\hat{\mathcal{U}}_1$ and $\hat{\mathcal{L}}_1$ do not resemble one another. (c) Potential CRs $\hat{\mathcal{U}}_0$ and $\hat{\mathcal{L}}_0$ overlayed on $\mathcal{U}_0$ and $\mathcal{L}_0$. As the boundaries of $\hat{\mathcal{U}}_0$ and $\hat{\mathcal{L}}_0$ are close to one another (both strongly resemble the circular level set on the right hand side of the image), we can infer that the point estimate for the level set, $\hat{\mu}_n^{-1}(0)$ (not depicted), which will lie within the `circular band' $\mathcal{S}\setminus(\hat{\mathcal{U}}_0\cup\hat{\mathcal{L}}_0)$ by construction, is a reliable estimate.}
    \label{fig:crs}
\end{centering}
\vspace{0.2cm}

When the level set $\mu^{-1}(c)$ has non-zero measure, such as when there is a plateau over which no signal occurs (e.g. $\mu=c=0$ over a measurable region), we can interpret CRs as providing a probabilistic bound on the region $\mu^{-1}(c)$. Specifically, CRs provide regions over which, with $(1-\alpha)\%$ confidence, we can state that $\mu$ does not equal $c$. Formally, we can assert that $\mu^{-1}(c)\subseteq \mathcal{S}\setminus (\hat{\mathcal{U}}_c\cup\hat{\mathcal{L}}_c)$ with confidence level $1-\alpha$. This is illustrated in Fig. \ref{fig:crs} Panel (b), where we conclude with confidence $1-\alpha$ confidence that $\mu$ is non-zero within the bright red and blue regions $\hat{\mathcal{U}}_1$ and $\hat{\mathcal{L}}_1$.

A conventional approach to constructing CRs is to define $\hat{\mathcal{L}}_c$ and $\hat{\mathcal{U}}_c$ as follows:
\begin{equation}\label{eq:cr_def1}
    \hat{\mathcal{L}}_c:=\bigg\{s \in \mathcal{S}: \hat{\mu}_n(s)< c+q\tau_n \sigma(s)\bigg\} \quad \text{ and }\quad \hat{\mathcal{U}}_c:=\bigg\{s \in \mathcal{S}: \hat{\mu}_n(s)> c+q\tau_n \sigma(s)\bigg\},
\end{equation}
where $\tau_n$ is a positive sequence such that $\tau_n\rightarrow 0$, $\sigma:\mathcal{S}\rightarrow \mathbb{R}$ represents standard deviation across samples, and $q$ is a quantile chosen to ensure the inclusion probability equals $1-\alpha$. Under this construction, the main focus of the CR literature is how best to choose the value of $q$ in practice.

Broadly speaking, most approaches to choosing the quantile $q$ share a common structure. First, a functional Central Limit Theorem (fCLT) or similar form of convergence is assumed. Following this, assumptions are placed on either (i) the continuity of the limiting process, or some variation thereof, or (ii) the convergence of the quantiles of some random variable, typically a suprema of the limiting process taken over carefully constructed spatial regions. Under these assumptions, the problem of estimating the quantile $q$ can be reduced to estimating the PDF of a suprema of a random process. This estimation procedure is performed empirically using computational methods (e.g. bootstrapping schemes, Bayesian updating, MCMC methods, etc). For instance, \cite{SSS} assumes an fCLT of the form $\sqrt{n}(\hat{\mu}_n-\mu)\xrightarrow{d}G$, continuity of the limiting process, $G$, and places several topological assumptions on the signal $\mu$. Under these assumptions, the required quantile $q$ can be proven to be equal to the $(1-\alpha)\%$ quantile of the PDF of $\sup_{s\in\mu^{-1}(0)}|G(s)|$, which is estimated using a bootstrap procedure. A detailed literature review of such constructions is given in Section \ref{sec:background_crs}.

Despite the wide-ranging applicability of the current CR methodology, many settings exist in which such CRs cannot be generated due to violation of the aforementioned assumptions. In this work, we are interested in the class of models where $\sqrt{n}(\hat{\mu}_n-\mu)$ satisfy a pointwise CLT, but the limiting process $G$ has piecewise continuous, rather than continuous, realizations. In such settings, the conventional assumptions of an fCLT, continuity and convergence of quantiles are all violated (c.f. Section \ref{sec:background_crs}). Such processes are easily constructed, with a particularly notable example being provided, alongside further detail, in Section \ref{sec:background_pieceproc} (see Fig. \ref{fig:absolute2}).

In the following text, we shall relax the conventional assumptions used to generate CRs in order to develop a CR theory for ``piecewise continuous processes''. We shall achieve this goal by first developing a framework to describe such processes formally, alongside a new notion of convergence, resembling that of the $M_1$ Skorokhod convergence employed in one dimensional time-series analysis (c.f. \cite{whitt2002stochastic,kern2022skorokhod}). Following this, we shall employ this convergence notion to define ``suprema-preserving" random processes; random processes which converge to a piecewise limit in a natural, application-driven, manner. This concept shall be utilised to define a `restrained'-CLT; a notion similar to, but weaker than, the convergence assumed by an fCLT. By generalising the extended continuous mapping theorem of \cite{van1996weak}, we then develop a theory of CRs for such suprema-preserving processes and describe how our theory may be applied to practical settings.

The outline of this paper is as follows. Section \ref{sec:background} provides background information and embeds our work within the current literature and Section \ref{sec:preliminaries} introduces notation used throughout the text. Following this, we begin our exploration of the topic by considering sequences of (non-random) continuous functions converging to piecewise continuous limits. Specifically, in Section \ref{sec:motivation}, we provide several motivating examples of such sequences, illustrating the convergence notion we wish to consider. Then, in Section \ref{sec:convergence}, we formalise this notion by introducing the central concept of \textit{`convergence with restraint'}. In Section \ref{sec:random}, we then shift attention to random processes, utilising the notion of \textit{`convergence with restraint'} to define \textit{`suprema preserving'} transformations; functions which transform uniformly convergent random processes in a manner that guarantees their suprema converge in a predictable manner. Building on this idea, we provide an extension of the continuous mapping theorem and introduce the notion of a `restrained' CLT. Using this machinery, in Section \ref{sec:piecewise_copes}, we then introduce CRs for processes satisfying restrained CLTs. Finally, in Section \ref{sec:applications}, we illustrate how our theory may be applied in a range of practical settings, drawing links to the existing literature.

%% file: Article/background.tex
\subsection{Piecewise Processes}\label{sec:background_pieceproc}

In Section \ref{sec:preliminaries}, we formally define the term \textit{`piecewise process'}. To aid intuition, however, we first provide an illustrative example. This example not only demonstrates how such a process can arise but also serves as a concrete point of reference for the discussions that follow.

\begin{examplebox}
    \begin{example}\label{example:bg_abs}
        Suppose $\mathcal{S}=[0,1]$, $\gamma:\mathcal{S}\rightarrow \mathbb{R}$ is the continuous function depicted in Fig. \ref{fig:absolute1} and $\hat{\gamma}_n:\mathcal{S}\rightarrow\mathbb{R}$ is a continuous function satisfying a functional CLT of the form $\hat{G}_n:=\tau_n^{-1}(\hat{\gamma}_n-\gamma)\xrightarrow{d} G$ where $\hat{G}_n$ and $G$ are both continuous random processes. 
        \vspace*{0.4cm}
        
{\centering
    \includegraphics[width=0.9\textwidth]{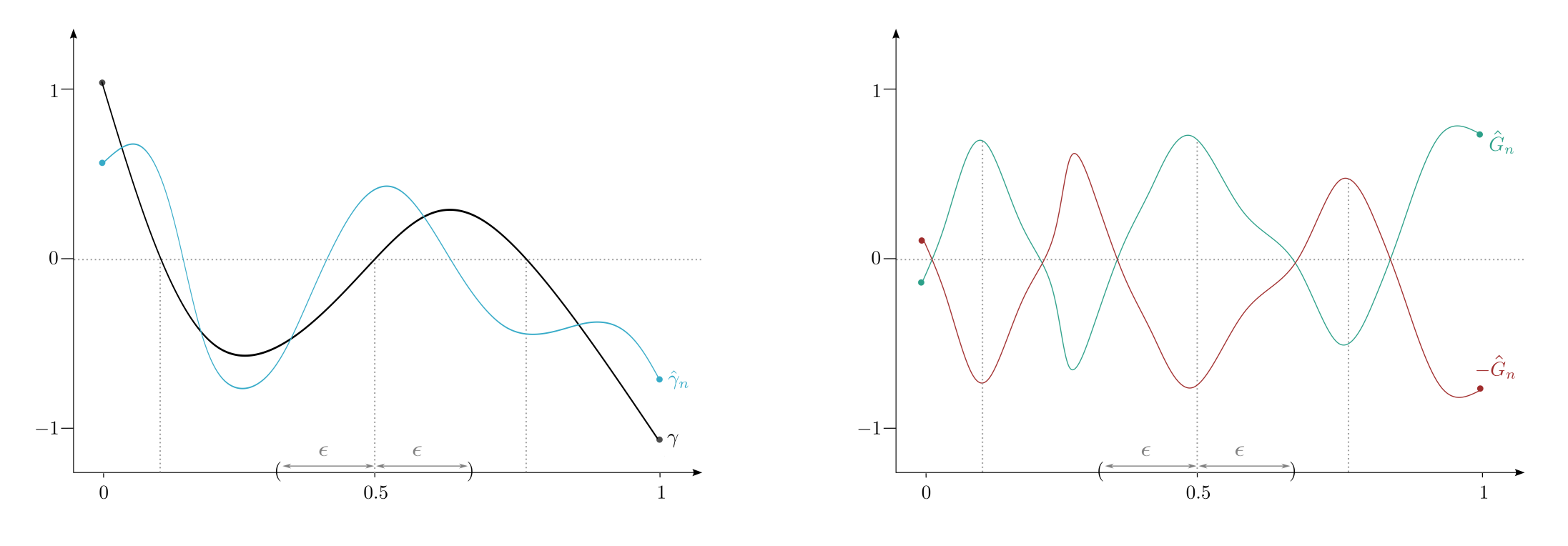}
    \captionof{figure}{Left: the function $\gamma$ and an instance of $\hat{\gamma}_n$. Right: an instance of the noise process $\hat{G}_n$, alongside $-\hat{G}_n$ (dashed). In both plots vertical lines highlight the points at which $\gamma=0$.}
    \label{fig:absolute1}}
    \vspace*{0.4cm}
        \noindent
        Assume further that $G$ satisfies $\mathbb{P}[G(s)=0]=0$ for any $s\in \mathcal{S}$. It is easily shown that, for each $s \in \mathcal{S}$, the below convergence holds pointwise:
        \begin{equation}\nonumber
            \hat{H}_n(s):=\tau_n^{-1}\big(|\hat{\gamma}_n(s)|-|\gamma(s)|\big) \xrightarrow{d} 
            H(s):=\begin{cases}
               G(s) & \quad\text{ if } \quad\gamma(s) > 0, \\
               -G(s) & \quad\text{ if } \quad\gamma(s) < 0, \\
               |G(s)| & \quad\text{ if } \quad\gamma(s) = 0. \\
            \end{cases}
        \end{equation}
        As realisations of $G$ are continuous, it follows that realisations of $H$ must be piecewise continuous (Fig. \ref{fig:absolute2}). Further, as $\gamma > 0$ on $(\sfrac{1}{2},\sfrac{1}{2}+\epsilon)$ and $\gamma <0$ on $(\sfrac{1}{2}-\epsilon,\sfrac{1}{2})$, we have that $s=\sfrac{1}{2}$ is a discontinuity point of $H$ if $G(\sfrac{1}{2})\neq -G(\sfrac{1}{2})$. By assumption we have that $\mathbb{P}[G(\sfrac{1}{2})\neq -G(\sfrac{1}{2})]=\mathbb{P}[G(\sfrac{1}{2})\neq 0]=1$. Thus, with probability one, $H$ is discontinuous whenever $\gamma$ changes sign.\\

        \vspace*{0.4cm}
        
{\centering
    \includegraphics[width=0.9\textwidth]{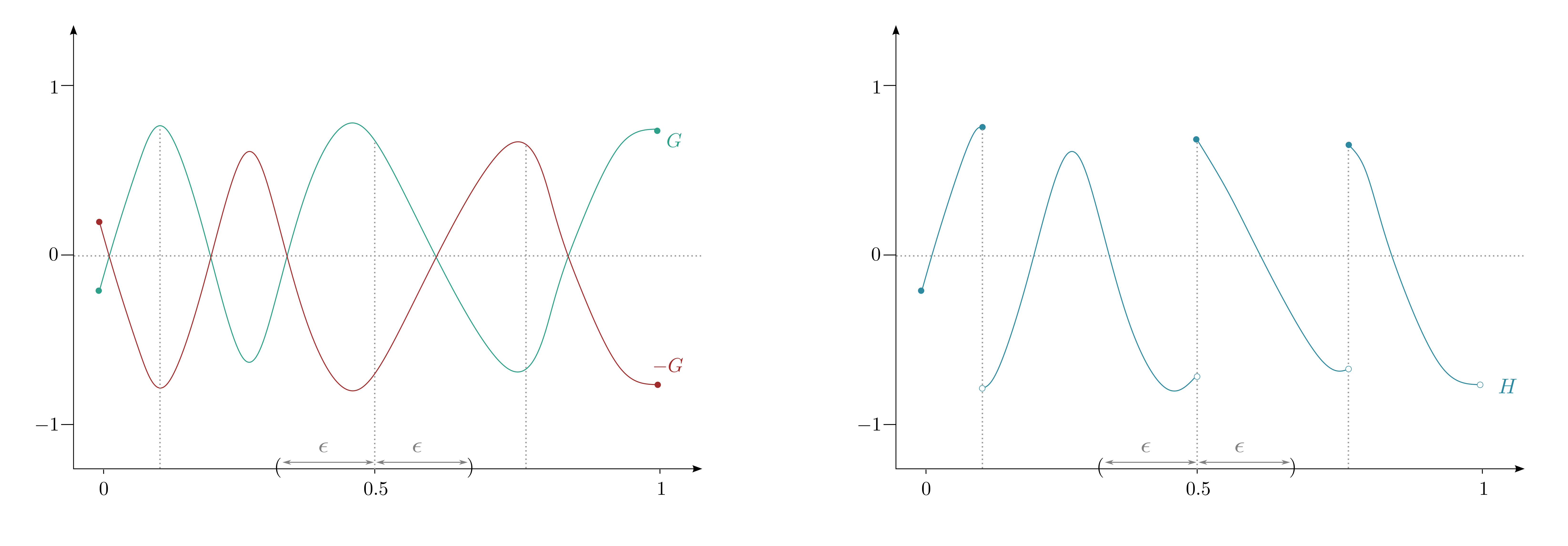}
    \captionof{figure}{Left: An instance of the limiting noise process $G$, alongside $-G$. Right: The corresponding instance of the limiting process $H$. Whenever $\gamma$ changes sign this process switches from $-G$ to $+G$, from which it can be seen that $H$ is discontinuous with probability one.}
    \label{fig:absolute2}}
    \vspace*{0.4cm}
        \noindent 
        However, instances of $\hat{H}_n$ are clearly continuous for all $n$ by definition. As a continuous function cannot converge uniformly to a discontinuous function we are left with the conclusion $\hat{H}_n$ does not converge to $H$ uniformly on $\mathcal{S}$.\\
        \\
        \noindent
        That is, $\tau_n^{-1}\big(|\hat{\gamma}_n(s)|-|\gamma(s)|\big)$ does not, and in general cannot, satisfy a functional CLT!
    \end{example}
    \vspace*{0.3cm}
\end{examplebox}

 Although the process $\tau_n^{-1}\big(|\hat{\gamma}_n(s)|-|\gamma(s)|\big)$ does not satisfy a fCLT, it is clear that it still is `well-behaved' in some sense. For instance, the final process $H$ does not have fractal-like properties nor does it ever possess vertical asymptotes.  Given these characteristics, one might expect such processes to have been well studied in the literature. However, this is not the case; in fact, as we detail below, the prevailing approaches for modelling piecewise processes cannot be used to model this example.

The conventional approach for modelling random processes resembling the above is to employ the Skorokhod topology. For signals sampled over a one-dimensional compact domain (typically the interval $[0,1]$), this involves treating instances of the process as points in a suitably topologized Skorokhod space. With this abstraction, convergence can be measured using a suitable choice of Skorokhod metric, results such as the Arzel\`{a}-Ascoli theorem become applicable, and central limit theorems for certain piecewise processes can be established \cite{skorokhod1956, van1996weak, whitt2002stochastic, kern2022skorokhod}. As a result, such methods have been employed to analyse many applied problems including image alignment, market equilibrium, wealth gaps, traffic flow and comparison of graphical models for dynamical systems (\cite{Glasbey2001,Liao2024,Goldie1977,Jabari2012,Deshmukh2017}). Extensions to processes defined over higher-dimensional spaces have also been proposed in the literature, but appear to be less widely employed (\cite{neuhaus1971weak,bickel1971convergence}).  

All Skorokhod-based approaches share a common assumption which limits their applicability to our setting. In one dimension, these approaches assume piecewise processes are right continuous with left limits, and in higher dimensions, this assumption is typically applied independently to each axis. The piecewise processes which motivate this work, however, are not typically equipped with a natural choice of left-right orientation. The example above illustrates this: the continuity of $H$ depends on the value of $G$, making it sometimes left-continuous and sometimes right-continuous.  

In fact, piecewise processes without a predetermined orientation arise naturally in a range of real-world settings. In disease mapping, socioeconomic conditions often exhibit discrete changes at provincial borders (\cite{GOEPP2024}). In geology, abrupt discontinuous jumps in rock formations frequently impact modelling considerations (\cite{Kim2005}). And in climatology settings, temperature variations are influenced by geographical features, such as large bodies of water, creating piecewise continuous noise processes (\cite{Yang2017,SSS,Munawar2024}). 

Furthermore, analytical considerations can also give rise to such `unorientable' piecewise processes. For example, noise terms in piecewise regression models for spatial data often exhibit discontinuous behaviour near boundary points (\cite{Qiu2004, lehner2024}), as do predictions from models incorporating spatial boundaries derived from classification tasks (\cite{Park2022}). Another common, but more subtle, source of such processes occurs when comparing images to one another, a practice commonly conducted manually in the field of functional Magnetic Resonance Imaging (fMRI) (\cite{Friston1999, Nichols2005}). We shall explore this setting in further detail in Section \ref{sec:applications}. 

In sum, we are unaware of any existing framework for analysing piecewise processes of the kind illustrated in Example \ref{example:bg_abs}. The existing approaches based on the Skorokhod topology require assumptions on the orientation of the processes which are not tenable for many settings of interest. To address this limitation, in this work we shall introduce a novel notion of convergence for functions, \textit{convergence with restraint} (see Section \ref{sec:convergence}), which can be used to describe the behaviour of random processes such as that of Example \ref{example:bg_abs} (see Section \ref{sec:random}).

\subsection{Confidence Regions}\label{sec:background_crs}

In this section, we provide a detailed background on confidence regions, situating our work within the existing literature. A wealth of literature on CRs has been published over the last three decades. Earlier works focused on the generation of CRs for specific applications including; Gaussian white noise models under isotropy or stationarity, wavelet processes that are uniform on Besov balls and certain classes of models defined over Sobelov spaces (\cite{Lindgren1995, Siegmund1995, Genovese2005, Cai2006}). In more recent literature, generalized constructions have been proposed allowing for CR generation across a range of domains. As noted in the introduction, although these approaches vary substantially, each makes restrictive assumptions which prevents their use in analysing the piecewise processes that we are interested in.

For instance, \cite{Reffaello2004} and \cite{Chen2017} propose general frameworks for CR generation in which the quantile $q$ in \eqref{eq:cr_def1} is evaluated via a Hausdorff loss criterion. To do so, both place strong assumptions on the continuity and differentiability of their processes to ensure the convergence in distribution of the Hausdorff distances between estimated and true sets (c.f. Theorem 2 of \cite{Reffaello2004} and Assumptions $K1$-$K2$ of \cite{Chen2017}). In our setting, however, Hausdorff convergences of this form cannot hold, as small perturbations in horizontal distance no longer correspond to continuous changes in $G$ when $G$ is allowed to be spatially discontinuous. For this reason, the methods of \cite{Reffaello2004} and \cite{Chen2017} cannot be applied to analyse piecewise processes. We emphasize here that, although many arguments in this work do involve the Hausdorff distance, our notion of restrained convergence does not require, nor does it necessitate, the Hausdorff convergence assumptions of \cite{Reffaello2004} and \cite{Chen2017}.

Another approach to evaluating the quantile $q$ is to first establish a relation of the form $\mathbb{P}[\hat{\mathcal{L}}_c \subseteq \mathcal{L}_c~ \wedge ~\hat{\mathcal{U}}_c \subseteq \mathcal{U}_c]\approx \mathbb{P}[\sup_{s\in \mu^{-1}(c)}|G(s)| \leq q]$ and then estimate the distribution $\sup_{s\in \mu^{-1}(c)}|G(s)|$ empirically. Such `spatial suprema'-based approaches have been pursued, for example, by \cite{Bugni2010}, \cite{Mammen2013} and \cite{Qiao2019}. However, the approaches adopted in these works all require that, for some appropriately defined empirical process $\hat{G}_n$, quantiles of the distribution of $\sup_{s\in \hat{\mu}_n^{-1}(c)}|\hat{G}_n(s)|$ converge to those of $\sup_{s\in \mu^{-1}(c)}|G(s)|$. To ensure such conditions are met, \cite{Bugni2010} places an assumption of stochastic equicontinuity on $\hat{G}_n$ (see \cite{Bugni2010} Assumption A4), \cite{Mammen2013} assumes the quantile convergence explicitly alongside specific convergence rates (see \cite{Mammen2013} Assumptions P1-3), and \cite{Qiao2019} place a range of differentiability, continuity and convergence assumptions on their processes (see \cite{Qiao2019} Assumptions F1-2, K, A, H1-2). In any case, such quantile convergence typically does not hold in our setting, as we shall often find ourselves in the situation where realizations of $\hat{G}_n$ are spatially continuous functions, but their limit is discontinuous. In such instances, as we shall see in Section \ref{sec:motivation}, it is easy to construct situations in which $\sup_{s\in \mu^{-1}(c)}|\hat{G}_n(s)|\not\xrightarrow{d}\sup_{s\in \mu^{-1}(c)}|G(s)|$. 

Many authors adopt the `spatial-suprema' approach but bypass the need for the above quantile convergence assumptions by instead assuming an fCLT of the form $\hat{G}_n\xrightarrow{d}G$, with additional regularity conditions placed on the signal $\mu$ and limiting field $G$. Examples of such an approach can be found in, for instance, \cite{French2013}, \cite{SSS}, \cite{JANKOWSKI2012}, \cite{Bowring:2019a} and \cite{MaullinSapey2022}. Alongside the aformentioned fCLT, these works all assume that $G$ has almost surely continuous sample paths (see, for example, Equation (3) of \cite{French2013}, Assumption 1b of \cite{SSS} and Assumption 2.2.1 of \cite{MaullinSapey2022}). In addition, \cite{French2013} assumes $\hat{G}_n$ has almost surely continuous sample paths, whilst \cite{SSS}, \cite{Bowring:2019a} and \cite{MaullinSapey2022} assume that $\hat{\mu}_n$ is almost surely continuous for $n$ large enough. Furthermore, \cite{SSS}, \cite{Bowring:2019a} and \cite{MaullinSapey2022} require topological conditions be placed on $\mu$ (see \cite{SSS} Assumption 1c), while \cite{French2013} uses a larger region for quantile estimates ($\mu^{-1}([c,\infty))$ instead of $\mu^{-1}(c)$) to avoid such concerns, resulting in conservative inference. \cite{JANKOWSKI2012} considers a class of Gaussian random fields with a pre-specified spatial covariance function, but, like \cite{French2013}, assumes $G$ and $\hat{G}_n$ have almost surely continuous sample paths, and bypasses the need for topological constraints of \cite{SSS} by using a larger region for quantile estimation. Unfortunately, as a result of the aforementioned assumptions, none of these works can be applied to our setting. This is with the exception of \cite{MaullinSapey2022} which, as we shall show in Section \ref{sec:applications}, can be viewed as a special case of our approach. 

Recently, a substantive generalisation and unification of the above approaches was proposed by \cite{telschow2023scope}. This work shows that CRs can be generated without the assumptions of almost sure continuity made by \cite{SSS} and \cite{Bowring:2019a}, nor the equicontinuity of \cite{Bugni2010} or fCLT assumptions of \cite{SSS} and \cite{Bowring:2019a}, and without adopting conservative approaches exhibited in \cite{French2013} and \cite{JANKOWSKI2012}. In fact, the central assumptions required by \cite{telschow2023scope} are only that of quantile convergence and tightness of $\hat{G}_n$ (see conditions (M1) and (M2) of \cite{telschow2023scope} Appendix C). Despite its being a substantial simplification of the existing assumptions employed in the literature, applying \cite{telschow2023scope} in our setting is difficult due to the abstract nature of Assumptions (M1) and (M2). It is not immediately clear under which assumptions the appropriate quantiles of piecewise processes converge, and this is a topic to which we shall devote much of the following text.

A wealth of altogether different approaches, which do not fall under the framework of Equation \eqref{eq:cr_def1}, have also been proposed. For instance, also requiring an fCLT alongside differentiability assumptions is \cite{Davenport2022}, which generates CRs specifically for peak locations under stationarity assumptions via MCMC methods. In the Bayesian setting, similar constructions have also been proposed to generate `credible' regions, with much literature focusing on computational approaches that simulate quantiles of an appropriately specified posterior distribution (\cite{Bolin2015, French2015}). All of the aforementioned constructions have found considerable application in a range of applied domains including, for instance, climatology, kernel density estimation, astrophysics and fMRI neuroscience (\cite{Bowring:2019a,Bowring:2020, Mejia2020, Hazra2021, Pelaez2022}). To the best of our knowledge, however, none can be applied in the piecewise setting.

%% file: Article/preliminaries.tex
Throughout this work, $\mathcal{S}$ shall denote a compact metric space, typically representing a known physical spatial domain or parameter space.
A random variable $G$ whose values are scalar functions from $\mathcal{S}$ to $\mathbb{R}$ shall be referred to as a \textit{`random process'} on $\mathcal{S}$. 
That is, given a probability space $(\Omega, \mathcal{F},\mathbb{P})$\nomenclature{$\Omega$}{Sample space}\nomenclature{$\mathcal{F}$}{Event space}\nomenclature{$\mathbb{P}$}{Probability Measure}, $G$ is a function which maps $\omega \in \Omega \mapsto f\in \mathbf{F}$, where $f$ is a function belonging to the function space $\mathbf{F}$\nomenclature{$\mathbf{F}$}{Arbitrary function space}. Following convention, the underlying sample space shall often be dropped from notation, so that $G(s)$ is used as shorthand for $G(\omega)(s)$; the value of the random process $G$ at spatial location $s\in \mathcal{S}$.  When $\mathbf{F}=C(\mathcal{S},\mathbb{R})$ (the space of continuous functions from $\mathcal{S}$ to $\mathbb{R}$), we say that $G$ is a \textit{continuous} random process.  

To develop a notion of convergence for piecewise processes, we shall first consider convergent sequences of fixed (non-random) functions. To make the distinction between fixed functions and random processes clear, we shall adopt the convention of using lower-case greek letters when discussing random processes converging to fixed limits (i.e. $\hat{\mu}_n\xrightarrow{d}\mu$), lower-case roman script letters when discussing sequences of fixed functions converging to fixed limits (i.e. $f_n\rightarrow f$), and upper-case roman script when denoting random processes converging to random processes (i.e. $G_n\xrightarrow{d} G$). Here, $\xrightarrow{d}$ represents weak convergence of the form described in \cite{van1996weak}, whilst $\rightarrow$ represents either (fixed) pointwise or uniform convergence, depending on the context.

In the following, we shall consider partitions of $\mathcal{S}$, denoted $\{V_i\}$, that are \textit{locally finite} - that is, for every $s \in \mathcal{S}$, there exists a neighbourhood of $s$ which intersects only finitely many of the sets $V_i$. Such partitions will be indexed by a countable set $\mathcal{I}$, often taken to be $\{1,...,m\}$ for some predefined $m\in\mathbb{N}$. In the non-random setting, we now define a \textit{piecewise continuous} function as follows.
\begin{mdframed}
    \begin{definition}[Piecewise Continuous Function] A function $f:\mathcal{S}\rightarrow \mathbb{R}$\nomenclature{$f$}{Non-random piecewise continuous function $f:\mathcal{S}\rightarrow \mathbb{R}$} is \textit{\textbf{piecewise continuous}} if there exists a locally finite collection of disjoint non-empty subsets $\{V_i\}_{i \in \mathcal{I}}$\nomenclature{$V_i$}{Spatial partition element. The partition $\{V_i\}_{i\in \mathcal{I}}$ of $\mathcal{S}$ is used to define piecewise continuous functions} such that $\mathcal{S}=\cup_i V_i$ and, for each $i$, $f|_{V_i}$ can be extended to a continuous function on $\overline{V_i}$. 
    \end{definition}
    \vspace*{0.25cm}
\end{mdframed}
We emphasize that the partition $\{V_i\}_{i \in \mathcal{I}}$ shall always be assumed to be fixed and locally finite. For brevity, we often write $\{V_i\}$ in place of $\{V_i\}_{i \in \mathcal{I}}$. Further, we shall denote the continuous extension of $f|_{V_i}$ to $\overline{V_i}$ as $f^i$ \nomenclature{$f^i$}{The continuous extension of $f\vert_{V_i}$ to $\overline{V_i}$} and,  when we wish to highlight the partition $\{V_i\}$, we say that $f$ is \textit{piecewise continuous subordinate to $\{V_i\}$}. Denoting the space of piecewise continuous functions subordinate to $\{V_i\}$ as $P(\{V_i\})$\nomenclature{$P(\{V_i\})$}{Space of piecewise continuous functions subordinate to $\{V_i\}$}, we can define a \textit{piecewise continuous} (random) process as a random process for which the function space $\mathbf{F}=P(\{V_i\})$.

The definition of random process we have used is based upon the work of \cite{van1996weak}. It is well known that within the framework of \cite{van1996weak}, many of the properties of random processes which we may be interested in studying are not measurable. To counter such issues, the notions of inner and outer probability are required. Given the probability space $(\Omega, \mathcal{F}, \mathbb{P})$, the outer and inner probabilities of arbitrary $B\subseteq\Omega$ are defined as:   \begin{equation}\nonumber
    \mathbb{P}^*[B]=\inf \{\mathbb{P}[A]: A \supseteq B, A \in \mathcal{F}\} \quad \text{and} \quad \mathbb{P}_{*}[B]=1-\mathbb{P}^{*}[\Omega\setminus B],
\end{equation}
respectively\nomenclature{$\mathbb{P}_{*},\mathbb{P}^{*}$}{Inner and outer probability, respectively}. Our use of inner and outer probability is motivated by the fact that the inner and outer probability is guaranteed to exist for all events we consider, even those which are non-measurable. In practical settings, the distinction between $\mathbb{P}^*/\mathbb{P}_*$ and $\mathbb{P}$ is often of little consequence, as the notions of inner, outer and standard probability are equivalent for measurable sets. In our derivations, the only consequence of employing inner/outer probability is that additivity cannot be automatically assumed, as it may fail on non-measurable sets.

Our ultimate aim is to quantify the uncertainty involved in estimating level sets of a fixed but unknown spatially varying function $\mu: \mathcal{S}\rightarrow \mathbb{R}$. Although this setting motivates the theory outlined in Sections \ref{sec:motivation} to \ref{sec:random}, we do not consider the function $\mu$ explicitly until Section \ref{sec:piecewise_copes}. To maximize the generality of our results, we assume in Section \ref{sec:piecewise_copes} that $\mu \in \ell^\infty(\mathcal{S})$\nomenclature{$\ell^\infty(X)$}{The set of bounded functions on a space $X$} where $\mu\in\ell^\infty(\mathcal{S})$ denotes the set of bounded, but not necessarily continuous, functions on $\mathcal{S}$. However, since our motivating applications typically involve continuous functions, we additionally assume in the applied examples of Section \ref{sec:applications} that $\mu$ is continuous, that is, $\mu \in C(\mathcal{S},\mathbb{R})$. Often, we shall treat $\mu$ as being constructed from a set of continuous functions. In such instances, which we shall denote this set as $\{\gamma^i\}$. For instance, in the application of Section \ref{sec:conjunc}, we consider $\mu$ of the form $\mu(s):=\min_{i}(\gamma^i(s))$. 

We denote by $\hat{\mu}_n$ an estimator of $\mu$ based on $n$ observations, modelled as a random process on $\mathcal{S}$. In a typical setting, we assume that $\hat{\mu}_n$ satisfies a pointwise (but not necessarily functional) Central Limit Theorem (CLT), i.e.,
\begin{equation}\nonumber
    \hat{H}_n(s):=\tau_n^{-1}(\hat{\mu}_n(s) -\mu(s))\xrightarrow{d}  H(s) \quad \text{pointwise},
\end{equation}
where $\hat{H}_n$ and $H$ are random processes taking values in $\ell^\infty(\mathcal{S})$ and $P(\{V_i\})$, respectively, and $\tau_n \to 0$ is a sequence of positive constants representing the convergence rate. A common choice is $\tau_n = \frac{1}{\sqrt{n}}$. For notational convenience, we omit an explicit variance term in all CLT expressions. The central idea of this work is to characterise the convergence of $\hat{H}_n$ by assuming uniform CLTs hold for estimators of $\{\gamma^i\}$, denoted $\{\hat{\gamma}^i_n\}$.

We shall employ of the following notation throughout the text. For arbitrary sets $X$ and $A$, we shall denote the distance between a point $x\in X$ and set $A\subseteq X$ as $d(x,A):=\inf_{a\in A}|x-a|$\nomenclature{$d(x,A)$}{Distance between a point $x$ and set $A$, defined as $d(x,A)=\inf_{a\in A}|x-a|$}. For an arbitrary spatial set $A\subseteq \mathcal{S}$, we denote its closure as $\overline{A}$\nomenclature{$\overline{A}$}{The topological closure of $A$ in $\mathcal{S}$} and complement as $A^c$ 
\nomenclature{$A^c$}{The complement of $A$ in $\mathcal{S}$}. In the case that a supremum is taken over the empty set, we define $\sup_{s \in\emptyset} f(s)=-\infty$ for arbitrary functions $f$.

%% file: Article/motivation.tex
As discussed in Sections \ref{sec:intro} and \ref{sec:background}, prior CR literature has largely considered cases where the limiting field $G$ is a continuous random process on $\mathcal{S}$. Here we instead focus on the setting where the limit, denoted $H$, is a piecewise continuous random process. To describe the convergence of such processes, we must first develop a theory for non-random sequences $f_n \in \ell^\infty(\mathcal{S})$ whose limits $f \in P(\{V_i\})$ are piecewise continuous subordinate to a fixed partition $\{V_i\}$ of $\mathcal{S}$.

Over the next two sections, we shall build upon the observations of Example \ref{example:bg_abs} to introduce a new notion of convergence. This convergence notion will be weaker than uniform convergence, yet sufficiently strong to ensure two informal criteria, which we shall explain and motivate below. These criteria will be essential for constructing CRs in later sections. We begin this section with some motivating examples that illustrate the behaviour we aim to capture.

\begin{examplebox}
\begin{example}\label{example:basic}
    Let $ \mathcal{S}=[0,1]$, $V_1=[0,\sfrac{1}{2}]$ and $V_2=(\sfrac{1}{2},1]$. Define;
    \begin{equation}\nonumber
    f_n(s)= \begin{cases} 1 & \text{if }s \in V_1, \\ (2-2s)^n & \text{if }s \in V_2, \\
    \end{cases} \quad \text{ and } \quad 
    f(s)= \begin{cases} 1 & \text{if }s \in V_1, \\ 0 & \text{if }s \in V_2. \\
    \end{cases} 
    \end{equation}
    Clearly, $f_n\rightarrow f$ pointwise,\nomenclature{$f_n$}{Functions converging to $f$. Mode of convergence is dependent upon context} but not uniformly (see Fig. \ref{fig:conv_fn1}). However, $\sup_A f_n(s)\rightarrow \sup_A f(s)$ for any closed $A\subseteq \mathcal{S}$. 
    
\centering
    \includegraphics[width=0.55\textwidth]{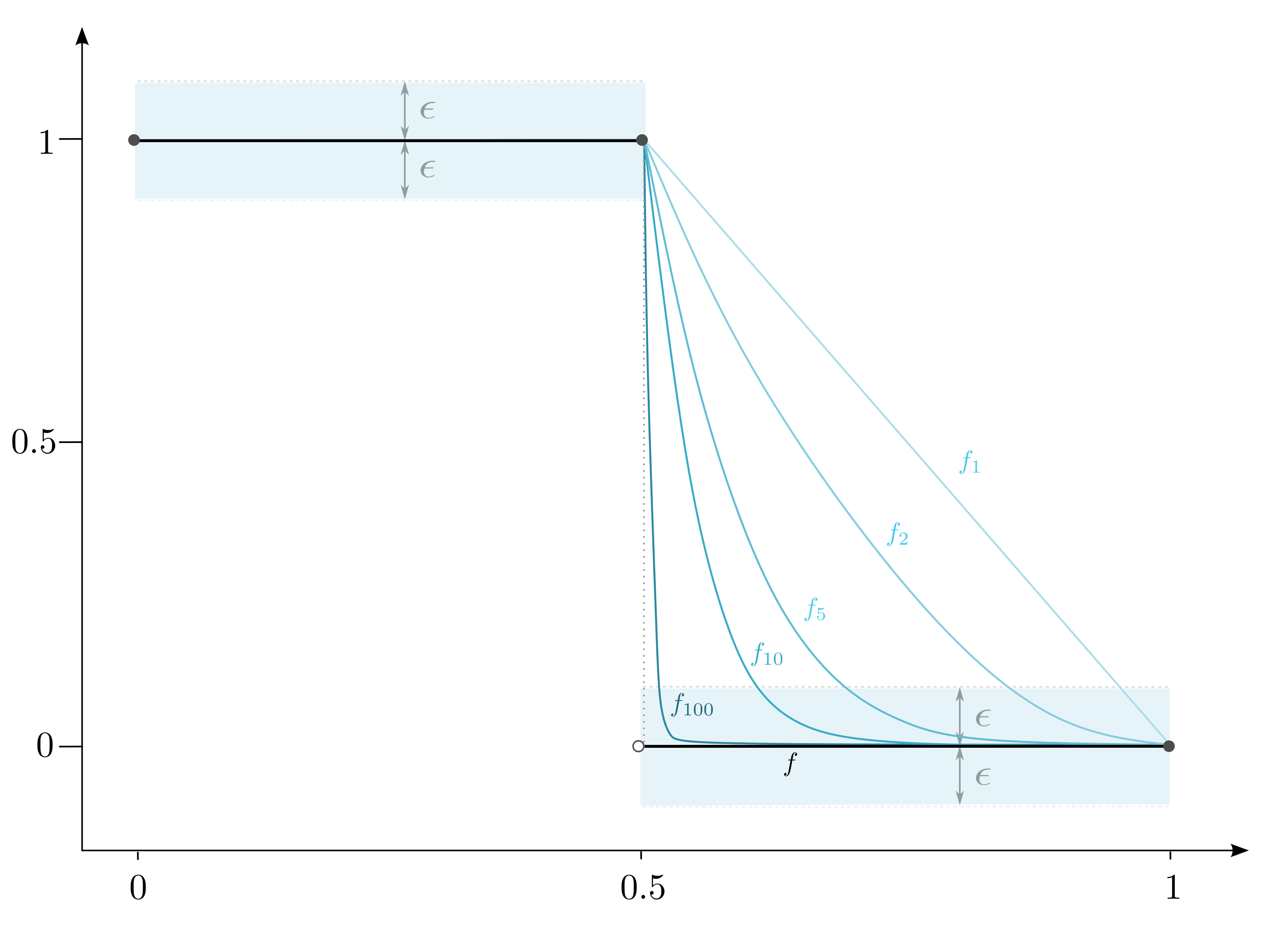}
    \captionof{figure}{An illustration of $f$ and $\{f_n\}$ for $n=1,2,5,10$ and $100$ in Example \ref{example:basic}. Also shown is a shaded blue region of uniform width $\epsilon$ about $f$. If it were the case that $f_n\rightarrow f$ uniformly, $f_n$ would eventually have to lie within the shaded blue region. However, in this example, for any $n$, $f_n$ will always leave the shaded blue region near $s=0.5$.}
    \label{fig:conv_fn1}
\end{example}
    \vspace*{0.25cm}
\end{examplebox}

Example~\ref{example:basic} illustrates the first of our two informal criteria; a property we will later build into our formal definition of convergence. Specifically, for any closed set $A \subseteq \mathcal{S}$, the limit superior of the suprema of $f_n$ over $A$ does not exceed the supremum of $f$ over $A$. Intuitively, this ensures that, for large $n$, the sequence $f_n$ cannot overshoot $f$ arbitrarily or diverge to infinity. Near discontinuities of $f$ (such as $s = 0.5$ in Fig.~\ref{fig:conv_fn1}), uniform convergence may fail; however, in the limit, $f_n$ remains bounded above by the larger of the left and right limits of $f$ at the point of discontinuity. For comparison, an example for which this criterion fails is given below.

\begin{examplebox}
\begin{example}\label{example:bad_converge}
    Define $\mathcal{S}$, $V_1$, $V_2$ and $f$ as in Example \ref{example:basic} and let;
    \begin{equation}\nonumber
    f_n(s)= \begin{cases} 1 & \text{if }s \in V_1, \\ (2-2s)^n + \Phi_n(s) & \text{if }s \in V_2, \\
    \end{cases} 
    \end{equation}
    where $\Phi_n(s)$ is the bump function of width $\sfrac{1}{5n}$ and height $1$ centered around $\sfrac{1}{2}+\sfrac{1}{10n}$.
    \begin{equation}\nonumber
        \Phi_n(s) =
\begin{cases}
\exp \left(1-(1-(10ns-5n-1)^2)^{-1}\right), & \text{ if }\frac{1}{2}\leq s \leq \frac{1}{2}+ \frac{1}{5n}, \\
0, & \text { otherwise. } \\
\end{cases}      
    \end{equation}
    As in Example \ref{example:basic}, we have that $f_n\rightarrow f$ pointwise, but not uniformly. However, for $A=[\sfrac{1}{2},1]$, the limit of $\sup_A f_n(s)$ is much larger than $\sup_A f(s)$ (see Fig. \ref{fig:bad_converge}).
    
    \centering
    \includegraphics[width=0.55\textwidth]{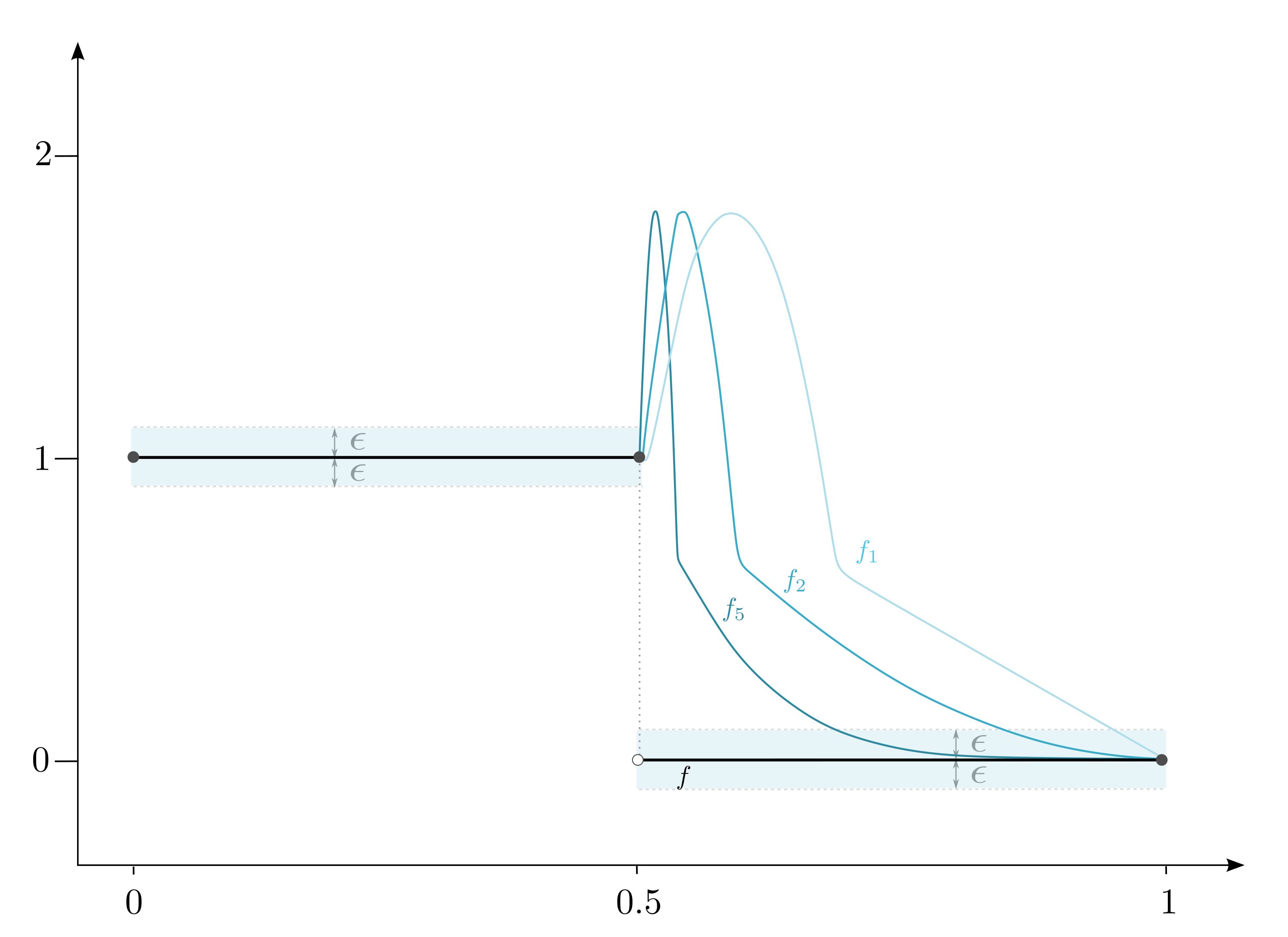}
    \captionof{figure}{An illustration of $f$ and $\{f_n\}$ for $n=1,2$ and $5$ in Example \ref{example:bad_converge}. Also shown is a shaded blue region of uniform width $\epsilon$ about the graph of $f$. For any $n$, $f_n$ will always leave the shaded blue region near $s=0.5$ and therefore $f_n$ does not converge to $f$ uniformly.}
    \label{fig:bad_converge}
\end{example}
    \vspace*{0.25cm}
\end{examplebox}

As can be seen in the above example, the assumption of pointwise convergence of $f_n$ to $f$ is not sufficient to guarantee that the suprema convergence criterion holds. In one dimension, a natural approach to ensuring this criterion is to first define the \textit{completed graph} of $f$ as
\begin{equation}\nonumber
    \Gamma_f := \left\{ \big(s,\alpha f^i(s) + (1-\alpha) f^j(s)\big) \in \mathcal{S} \times \mathbb{R} \;:\; i,j\in \mathcal{I},~ s \in V_i \cap \overline{V_j},\ \alpha \in [0,1]
    \right\}.
\end{equation}
That is, $\Gamma_f$ is the set of all points $(s',x)$ for which either $f(s') = x$ or $x$ lies between $f^i(s)$ and $f^j(s)$. Equivalently, one may view $\Gamma_f$ as the result of filling in all vertical segments in the graph of $f$ so as to connect the left and right limits at each point, analogous to drawing $f$ in a single continuous stroke without lifting the pen from the paper. Given this definition, the suprema criterion could be enforced by requiring that $f_n$ eventually lies within an $\epsilon$-neighbourhood of $\Gamma_f$ (as illustrated by the shaded blue region in Figure~\ref{fig:all3_res}).

\begin{examplebox}
    \begin{example}\label{example:all3_res}
        In Fig. \ref{fig:all3_res}, the functions from Examples \ref{example:basic} and \ref{example:bad_converge} are plotted alongside shaded $\epsilon$-neighbourhoods of the completed graphs of their pointwise limits. Here, the sequence of functions from Example \ref{example:basic} (left) exhibits the desired behaviour, with $f_n$ eventually lying in the shaded blue region, whilst the sequence of Example \ref{example:bad_converge} (right) does not.\\

\centering
    \includegraphics[width=\textwidth]{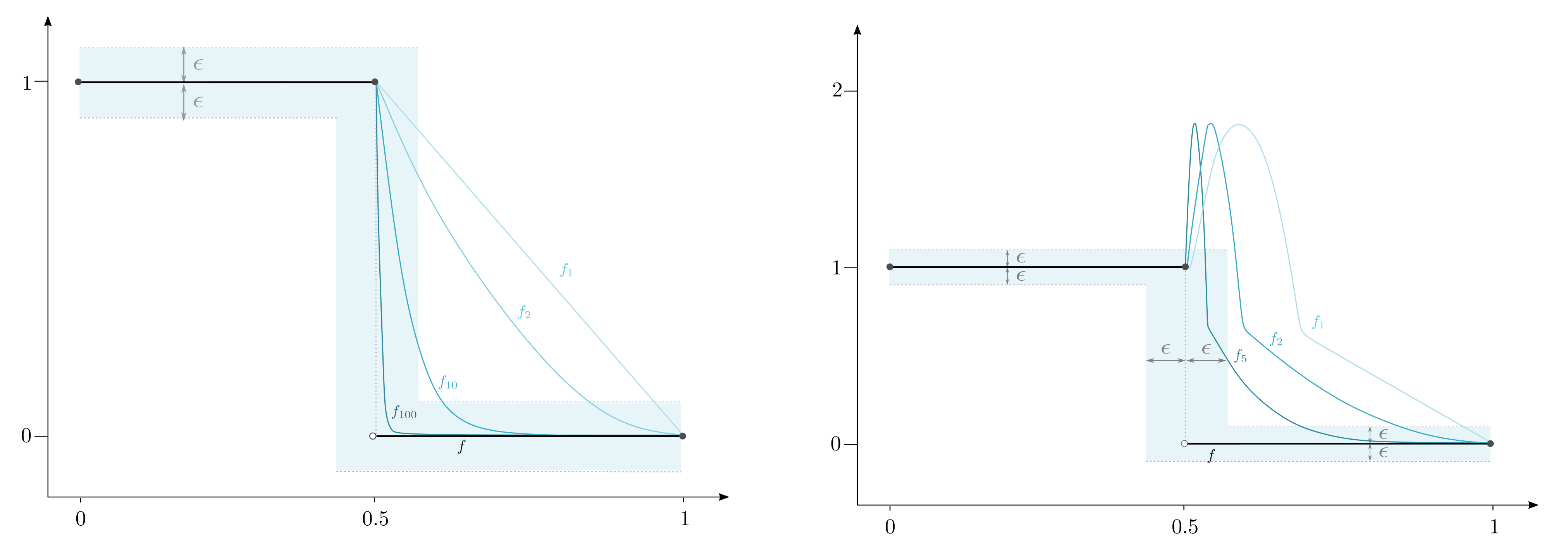}
    \captionof{figure}{Two sequences of functions which converge pointwise to $f=\mathbbm{1}_{[0,0.5]}$. Also shown is a shaded blue region representing the $\epsilon$-thickening of $\Gamma_f$. On the left, $f_n$ eventually lies within the blue region for all $n$. On the right, $f_n$ always exceeds the blue shaded region close to $0.5$. As noted in Examples \ref{example:basic} and \ref{example:bad_converge} an upper bound for $\sup f_n$ cannot be derived from $\sup f$ for the sequence on the right, but can be for those on the left.}
    \label{fig:all3_res}
        
    \end{example}
    \vspace*{0.25cm}
\end{examplebox}

The convergence to $\Gamma_f$ described above ensures our desired suprema criterion, but it may not be the suitable notion of convergence to employ in many practical scenarios. The reason for this is that the convergence of $f_n$ to $\Gamma_f$, as illustrated, does not exclude certain unusual oscillatory behaviours that are rarely encountered in practice. To clarify this point, consider the following example, in which $f_n$ oscillates near the discontinuity $s = 0.5$ with increasing frequency.

\begin{examplebox}
\begin{example}\label{example:bad_converge_weak}
    Let $ \mathcal{S}=[0,1]$, $V_1=[0,\sfrac{1}{2})$ and $V_2=[\sfrac{1}{2},1]$. Define $f$ as in Example \ref{example:mid_bdd} and let;
    \begin{equation}\nonumber
    f_n(s)= \begin{cases} 1 & \text{ if }0 \leq s \leq \frac{1}{2}-\frac{1}{4n}, \\ 
    \exp \left(1-(1-(4ns-2n+1)^2)^{-1}\right), & \text{ if } \frac{1}{2}- \frac{1}{4n} \leq s \leq \frac{1}{2}, \\
\exp \left(1-(1-(4ns-2n-1)^2)^{-1}\right), & \text{ if }\frac{1}{2}\leq s \leq \frac{1}{2}+ \frac{1}{2n}, \\
0, & \text { otherwise. } \\
    \end{cases} 
    \end{equation}
    In words, $f_n$ is equal to one up until $s=\sfrac{1}{2}-\sfrac{1}{4n}$. Following this, $f$ decreases in value until $s=\sfrac{1}{2}$, at which point $f(s)=0$. After $s=\sfrac{1}{2}$, $f_n$ follows the shape a bump function of width $\sfrac{1}{2n}$ and height $1$, centered around $\sfrac{1}{2}+\sfrac{1}{4n}$ (see fig \ref{fig:bad_converge_weak}).
    
    {\centering
    \includegraphics[width=0.55\textwidth]{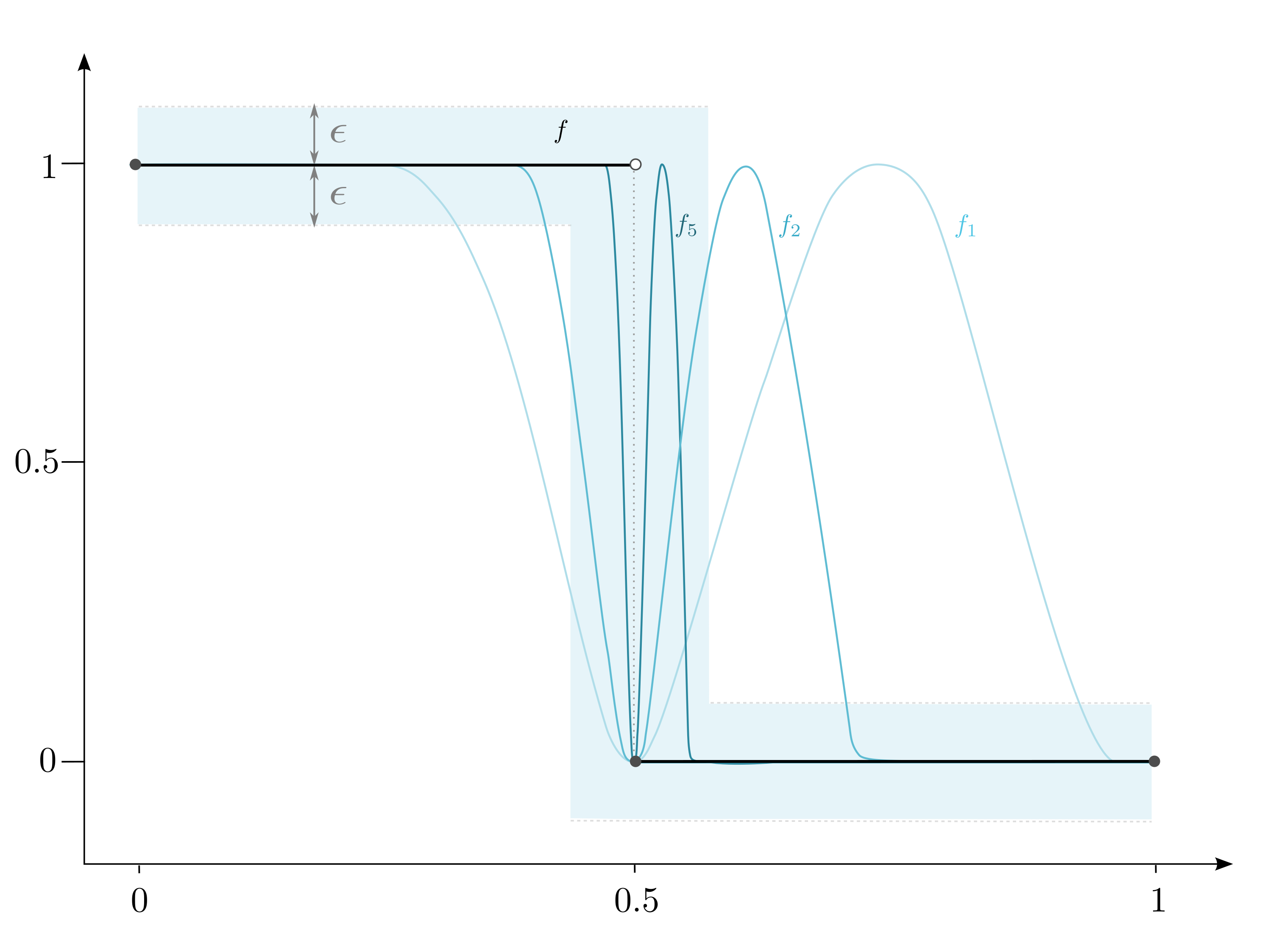}
    \captionof{figure}{An illustration of $f$ and $\{f_n\}$ for $n=1,2$ and $5$ in Example \ref{example:bad_converge_weak}. Also shown is a shaded blue region of uniform width $\epsilon$ about $f$. Although $f_n$ eventually is contained within the shaded blue region, it still exhibits strange, and potentially undesirable, oscillatory behaviour near $s=0.5$.}
    \label{fig:bad_converge_weak}}
    \vspace*{0.2cm}
    \noindent In this example, we still have that $f_n$ eventually lies within any $\epsilon-$thickening of $\Gamma_f$. However, as $n\rightarrow \infty$, $f_n$ becomes increasingly `pinched' near $s=\sfrac{1}{2}$. 
\end{example}
    \vspace*{0.25cm}
\end{examplebox}

The oscillatory behaviour observed in Example \ref{example:bad_converge_weak} is not only uncommon in practice, but allowing it can weaken the bounds that can be placed on CR coverage. Consequently, the second criterion motivating our new notion of convergence is the exclusion of such oscillations, thereby ensuring stricter control over local fluctuations $f_n$.

We emphasise that the above two criteria are intended only as motivating heuristics; our formal notion of convergence, introduced in the next section, will encode these ideas precisely in a form suitable for our setting. As noted in Section \ref{sec:background_pieceproc}, one possible approach to formalising the two criteria illustrated above in the one-dimensional setting is to employ the M1 Skorokhod topology (c.f. \cite{skorokhod1956,whitt2002stochastic,kern2022skorokhod}). However, the Skorokhod approach imposes restrictive orientation-based assumptions (see Section \ref{sec:background_pieceproc}) and thus extensions of the Skorokhod topology to the general $k$-dimensional case, suited to our purposes, have not been explored in the literature. The notion of convergence we shall propose is in fact slightly weaker than M1 Skorokhod, allowing greater scope for generalisation in higher dimensions. A full exploration of the relationship between our notion of convergence and Skorokhod convergence is beyond the scope of this paper but shall be pursued in future work.

The two criteria discussed above capture the type of \textit{`restrained'} behaviour we wish our convergence notion to describe. Having developed some intuition for this behaviour, in the next section we shall define the notion explicitly and examine its properties.

%% file: Article/convergence.tex
In this section, we introduce our new notion of convergence, which we call convergence with restraint. Before presenting its definition, we establish a few short results on pointwise convergence that will aid the discussion and serve as useful preliminaries for constructing CRs in Section \ref{sec:piecewise_copes}.

\subsection{Pointwise Convergence}

In the following sections, we shall develop theory allowing us to relate the suprema of $f$ to those of $f_n$ as $n \rightarrow \infty$. To do so, we first note that if $f_n \rightarrow f$ pointwise, then, in some instances, suprema of $f_n$ can be bounded below using suprema of $f$. This statement is made explicit in the following theorem.

\begin{mdframed}
    \begin{theorem}[Properties of Pointwise Convergence]\label{thm:pntwise_props} Suppose $\{A_k\}$\nomenclature{$A_k$}{Arbitrary subset of $\mathcal{S}$} are arbitrary subsets of $\mathcal{S}$. If $f_n \rightarrow f$ pointwise, then:
        \begin{equation}\nonumber 
        \sup_{s \in A^-_\infty}f(s) \leq \liminf_{n\rightarrow\infty}\bigg(\sup_{s\in A_n}f_n(s)\bigg), 
        \end{equation}
    where $A^-_\infty:=\liminf_{n\rightarrow\infty}A_n = \bigcup_{m \geq 1}\bigcap_{n\geq m} A_n$.\nomenclature{$A^-_\infty$}{The limit inferior of $A_n$ as $n\rightarrow \infty$} If, in addition, $\{A_k\}$ is decreasing (i.e. $A_1 \supseteq A_2 \supseteq ...$) then $A_\infty^{-}=\cap_{n\geq 1} A_n$. 
    \end{theorem}
    \vspace*{0.25cm}
\end{mdframed}
\begin{proof} First, suppose $A_n=\emptyset$ infinitely often. It follows that $A^-_\infty=\emptyset$ and thus $\sup_{A^-_{\infty}}f=-\infty$. In this case, the inequality holds. 
    
    For the case in which $A_n=\emptyset$ for only finitely many $n$, assume without loss of generality that $A_n\neq \emptyset$ for all $n$. Let $s^* \in A^-_\infty=\cup_{m \geq 1}\cap_{n\geq m} A_n$. It follows that there must exist an $m\in\mathbb{N}$ such that for all $n\geq m$, $s^* \in A_n$. It follows that for all $n\geq m$;
    \begin{equation}\nonumber
        f_n(s^*) \leq \sup_{s \in A_n}f_n(s).
    \end{equation}
    Applying $\liminf$ to both sides, and noting that $f_n(s^*)$ converges pointwise to $f(s^*)$, yields:
    \begin{equation}\nonumber
        f(s^*)=\liminf_{n\rightarrow\infty}f_n(s^*) \leq \liminf_{n\rightarrow\infty}\bigg(\sup_{s \in A_n}f_n(s)\bigg).       
    \end{equation}
    Noting that the above holds for all $s^\ast\in A^-_\infty$ gives:
    \begin{equation}\nonumber 
        \sup_{s \in A^-_\infty}f(s) \leq \liminf_{n\rightarrow\infty}\bigg(\sup_{s\in A_n}f_n(s)\bigg).
    \end{equation}
    The identity for decreasing $\{A_k\}$ follows directly from the definitions.
\end{proof}

Theorem \ref{thm:pntwise_props} shows that if $f_n\rightarrow f$ pointwise, the suprema of $f_n$ taken over arbitrary sets $A_n$ is bounded below by a function of $f$. However, as illustrated by the function in Example \ref{example:bad_converge}, the suprema of $f_n$ cannot be bounded above in a similar way. In fact, near a discontinuity of $f$, the supremum of $f_n$ can become arbitrarily larger than that of $f$. To handle this, we introduce the notion of \textit{`Convergence with Restraint'}.

\subsection{Convergence with Restraint}

In Section \ref{sec:motivation}, we outlined two informal criteria motivating our development of a new convergence notion: $(i)$ convergence of $f_n$ to $\Gamma_f$, and $(ii)$ exclusion of pathological oscillatory behaviours such as those in Example \ref{example:bad_converge_weak}. We now formalise these ideas in the definition of \textit{restrained bound}.
\begin{mdframed}
    \begin{definition}[Restrained Bound]\label{def:res_bdd1}
        Suppose $f$ is a piecewise continuous function subordinate to $\{V_i\}$. $f$ is said to be an \textit{upper restrained bound} for a sequence of functions, $\{f_n\}$, if the following condition is met:
    \begin{equation}\nonumber 
    \forall i,~\forall s \in \overline{V_i},~\lim\limits_{\delta \downarrow 0}\limsup\limits_{n\rightarrow\infty}\sup\limits_{s'\in B_{\delta}(s)\cap V_i}f_n(s')\leq \max\bigg(f^i(s),f(s)\bigg).
    \end{equation}
    If $-f$ is an upper restrained bound for the sequence $\{-f_n\}$, then we say that $f$ is a \textit{lower restrained bound} for $\{f_n\}$. If $f$ is both an upper and lower restrained bound for the sequence $\{f_n\}$, then we say that $f$ is simply a \textit{restrained bound} for $\{f_n\}$.
    \end{definition}
    \vspace{0.2cm}
\end{mdframed}

Several useful lemmas concerning elementary properties of restrained bounds may be found in Supplemental Material Section \ref{supp:resbdds}. As it shall be of particular use in Section \ref{sec:applications}, we include here the below theorem concerning sums of restrained bounds.

\begin{mdframed}
\begin{theorem}[Sums of Restrained Bounds]\label{thm:SumUpperRestraint}
	Let $\{f_n\},\{g_n\} \subset \ell^\infty(\mathcal{S})$ be arbitrary sequences with upper restrained bounds $f$ and $g$ subordinate to $\{U_i\}_{i\in\mathcal{I}}$ and $\{W_j\}_{j\in\mathcal{J}}$, respectively. Define $V_{ij} := U_i\cap W_j$ for $i \in \mathcal{I}, j \in \mathcal{J}$ and
\begin{equation*}
	\mathcal{N} :=
	\bigcup_{(i,j)\in \mathcal{I} \times \mathcal{J}} \overline{U_i}\backslash U_i \cap \overline{W_j}\backslash W_j \cap \overline{V_{ij}}.
\end{equation*}
If, for all $s \in \mathcal{N},i \in \mathcal{I}$ and $ j \in \mathcal{J}$  it holds that
\begin{equation}\label{eq:cond}
	\max\left(f^i(s) + g(s), f(s) + g^j(s) \right) \leq \max\left(f^i(s) + g^j(s), f(s) + g(s) \right),
\end{equation}
then it holds that $f + g$ is an upper restrained bound for $\{f_n + g_n\}$.
\end{theorem}
\vspace{0.2cm}
\end{mdframed}
\begin{proof}
    See Supplementary Material Section \ref{supp:resbdds}.
\end{proof}
Given the above definition of \textit{restrained bound},  we are now in a position to introduce the convergence notion which we wish to investigate: \textit{convergence with restraint}.

\begin{mdframed}
    \begin{definition}[Convergence with Restraint]\label{def:res_def} Suppose $f$ is a piecewise continuous function subordinate to $\{V_i\}$, and $\{f_n\}$ is a sequence of functions on $\mathcal{S}$.  We say that $\{f_n\}$ \textit{converges to $f$ with upper restraint}, denoted $f_n \xrightarrow{u.r.} f$, if $f_n \to f$ pointwise and $f$ is a restrained upper bound for $\{f_n\}$. Convergence with \textit{lower restraint} is defined analogously, and is denoted $f_n \xrightarrow{l.r.} f$. If $f_n$ converges to $f$ with both upper and lower restraint, we say that $\{f_n\}$ \textit{converges to $f$ with restraint}, denoted $f_n \xrightarrow{r} f$.
    \end{definition}
    \vspace*{0.3cm}
\end{mdframed}
An immediate consequence of the above definitions is that $f_n\xrightarrow{r} f$ if and only if $f_n\rightarrow f$ pointwise and $f$ is a restrained bound for $\{f_n\}$.
\begin{remark} It is straightforward to verify that the sequence of functions in Example~\ref{example:basic} converges with restraint, whereas those in Examples~\ref{example:bad_converge} and~\ref{example:bad_converge_weak} do not.  
For instance, to see that Example~\ref{example:bad_converge_weak} fails to converge with restraint, take $i = 2$ and $s = 0.5$.  
In this case, the limit on the left-hand side of Definition~\ref{def:res_def} is equal to $1$: for any $\delta > 0$, the set $B_\delta(s) \cap [\sfrac{1}{2}, 1]$ eventually contains points where $f_n$ equals $1$. However, the right-hand side simplifies to $f^2(s) = f(s) = 0$, showing that the bound does not hold.
\end{remark}

The notions of uniform convergence and convergence with restraint are closely related. While convergence with restraint does not, in general, imply uniform convergence, the converse does hold, as established in the following lemma.

\begin{mdframed}
    \begin{lemma}
        If $f$ is piecewise continuous subordinate to $\{V_i\}$ and $f_n\xrightarrow{}f$ uniformly, then $f_n \xrightarrow{r}f$.
    \end{lemma}
    \vspace{0.3cm}
\end{mdframed}
\begin{proof}
    First, note that if $f_n\xrightarrow{}f$ uniformly then $f_n\xrightarrow{}f$ pointwise. We must now show that $f$ is a restrained bound for $\{f_n\}$. To show the condition in Definition \ref{def:res_bdd1}, fix the values of $i\in\mathcal{I}$, $s \in \overline{V_i}$, and $\delta > 0$. By uniform convergence, we have that for all $s' \in B_\delta(s)\cap V_i$;
    \begin{equation}\nonumber
        f_n(s')-f^i(s') = f_n(s')-f(s') < \epsilon_n
    \end{equation}
    for some $\epsilon_n\rightarrow 0$ as $n\rightarrow \infty$. Thus, $f_n(s') < f^i(s') + \epsilon_n$. Noting that $f^i$ is continuous on $\overline{V_i}$, we obtain:
    \begin{equation}\nonumber
        f_n(s') - f^i(s) < f^i(s') - f^i(s) + \epsilon_n \leq \omega_\delta(f^i) + \epsilon.
    \end{equation}
    where $\omega_\delta(f^i)$ is the modulus of continuity of $f^i$, which tends to $0$ as $\delta\rightarrow 0$. Thus $f_n(s') < f^i(s) + \omega_\delta(f^i) + \epsilon$. Applying the appropriate limits and suprema to both sides, it now follows that:
    \begin{equation}\nonumber
        \lim\limits_{\delta \downarrow 0}\limsup\limits_{n\rightarrow\infty}\sup\limits_{s'\in B_{\delta}(s)\cap V_i}f_n(s')\leq f^i(s) \leq \max\bigg(f^i(s),f(s)\bigg).
    \end{equation}
    Thus, $f$ is an upper restrained bound for $\{f_n\}$. As identical logic shows that $f$ is a lower restrained bound, this completes the proof.
\end{proof}

Analogous to the lower bound obtained in Theorem~\ref{thm:pntwise_props}, the notion of restrained bound can be used to establish upper bounds on suprema of $f_n$. This idea is formalized in the theorem below.

\begin{mdframed}
    \begin{theorem}[Properties of Upper Restrained Bounds]\label{thm:res_prop}
    Suppose $\{A_k\}$ are arbitrary subsets of $\mathcal{S}$. If $f$ is a piecewise continuous function subordinate to $\{V_i\}$ and $f$ is an upper restrained bound for $f$, then:
         \begin{equation}\nonumber
         \limsup_{n\rightarrow\infty}\bigg(\sup_{s\in A_n}f_n(s)\bigg) \leq \sup_{i\in\mathcal{I}}\sup_{s \in A_\infty^{i,+}}f^i(s)
        \end{equation}
    where $A^{i,+}_\infty:=\cap_{k\geq 1}\overline{V_i \cap \overline{\cup_{n \geq k}A_n}}$. If, in addition, $\{A_k\}$ is decreasing (i.e. $A_1 \supseteq A_2 \supseteq ...$) then $A^{i,+}_\infty=\cap_{n\geq 1}\overline{V_i \cap \overline{A_n}}$. 
    \end{theorem}
    \vspace*{0.25cm}
\end{mdframed}
\begin{remark}
    The $+$ and $-$ notation in this theorem and Theorem \ref{thm:pntwise_props} has been chosen to highlight that the sets in question resemble limit suprema and inferior, respectively. To be specific, $A^-_{\infty}$ is the limit inferior of $\overline{A_n}$, whilst $A^{+}_\infty$ closely resembles, but is not equal to, the limit suprema of $V_i\cap\overline{A_n}$, defined as: 
    \begin{equation}\nonumber
        \limsup_{n\rightarrow\infty}\bigg(V_i\cap\overline{A_n}\bigg) := \bigcap_{k\geq 1}\bigg(V_i \cap \bigcup_{n \geq k}\overline{A_n}\bigg).
    \end{equation}
    The key difference between the above expression and the definition of $A^{i,+}_\infty$ is the placement of closures inside the expression. The set in the theorem contains two closures, neither of which is equivalent to the closure in the above. This nuance is a result of the fact that we have placed no assumptions on whether $\{V_i\}$ or $\{A_n\}$ need be open or closed. It is worth emphasizing that, without the assumption of $f$ being an upper restrained bound, no bound of the form given by Theorem \ref{thm:res_prop} can be placed on $f_n$ at all (c.f. Example \ref{example:bad_converge}). 
\end{remark} 
\begin{proof} To begin, suppose there were some $k^*$ such that $A_n=\emptyset$ for all $n\geq k^*$. In this case, we would have that, for each $i\in\mathcal{I}$:
\begin{equation}\nonumber
    A^{i,+}_\infty=\bigcap_{k\geq 1}\overline{ V_i\cap \overline{\bigcup_{n \geq k}A_n} }\subseteq \overline{ V_i\cap\overline{\bigcup_{n \geq k^{*}}A_n} } =\emptyset.
\end{equation}
It follows that in this case both sides of the inequality in Theorem \ref{thm:res_prop} equal $-\infty$ and the result trivially holds.

Now suppose $A_n\neq \emptyset$ for infinitely many $n$. For proof by contrapositive, we now assume the negation of the inequality in Theorem \ref{thm:res_prop}, i.e. we assume that:
\begin{equation}\nonumber
  \sup_{i\in\mathcal{I}}\sup_{s \in A_\infty^{i,+}}f^i(s)< \limsup_{n\rightarrow\infty}\bigg(\sup_{s\in A_n}f_n(s)\bigg) 
\end{equation}
For ease denote the left-hand side as $L$ and the right-hand side as $R$. Further, denote the left-hand side of the expression in Definition \ref{def:res_bdd1} as $K_{(i,s)}$ and the right-hand side as $K'_{(i,s)}$. To prove the result, we shall show the contrapositive. That is, we shall show;
\begin{equation}\nonumber
    L < R \implies K'_{(i,s^*)} < K_{(i,s^*)}
\end{equation}
for some $i \in \mathcal{I}$ and $s^* \in \overline{V_i}$. To begin, let $\{s_m\}\subset \mathcal{S}$ be such that $s_m \in A_{n_m}$ for each $m$ and:
\begin{equation}\nonumber
\begin{split}
     R& := \limsup_{n\rightarrow\infty}\bigg(\sup_{s\in A_n}f_n(s)\bigg) = \lim_{m\rightarrow \infty} f_{n_m}(s_m).
\end{split}
\end{equation}
This is possible by the definitions of $\limsup$ and $\sup$. Now as $\mathcal{S}$ is compact, $\{s_m\}$ must have a convergent subsequence tending to a limit $s^*$. Further, as the partition $\{V_i\}$ is locally finite, this subsequence must lie in $V_i$ infinitely often for some $i\in \mathcal{I}$. Without loss of generality assume that $\{s_m\}$ is such a subsequence (i.e. $s_m\rightarrow s^*$ and $s_m \in V_i$ for all $m$). It follows that $f^i(s_m)=f(s_m)$ for all $m$. Finally, let $j\in\mathcal{I}$ be such that $s^*\in V_j$ so that $f^j(s^*)=f(s^*)$ (note that $j$ may or may not equal $i$).

We now have that for any $k \geq 1$:
\begin{equation}\nonumber
    d\bigg(s^*, V_i \cap \bigcup_{n\geq k} A_n\bigg) \leq d(s^*,s_{m}) +d\bigg(s_{m},   V_i\cap\bigcup_{n\geq k} A_n\bigg)
    \xrightarrow{m\rightarrow \infty} 0
\end{equation}
where the convergence follows from the fact that $s_m\rightarrow s^*$ by construction, and, if $m$ is large enough that $n_m>k$, we have $s_m\in V_i \cap A_{n_m}\subseteq V_i \cap \cup_{n\geq k} A_n$. Thus $s^* \in \overline{ V_i\cap \cup_{n\geq k} A_n}$. As this holds for all $k\geq 1$, we have that:
\begin{equation}\nonumber
    s^* \in \bigcap_{k \geq 1}\overline{V_i \cap \bigcup_{n\geq k} A_n}\subseteq A^{i,+}_\infty.
\end{equation}
Further, as $s^* \in V_j$ and, by the above equation, $s^* \in \overline{\cup_{n\geq k}A_n}$ for all $k$, we have that:
\begin{equation}\nonumber
    s^* \in \bigcap_{k \geq 1}\bigg(V_j \cap \overline{\bigcup_{n\geq k} A_n}\bigg)\subseteq A^{j,+}_\infty.
\end{equation}
Now choose $\delta > 0$ and note that, as $s_m\in V_i$ and $s_m\in B_\delta(s^*)$ for $m$ large enough:
\begin{equation}\nonumber  
\begin{split}
    &\limsup_{n\rightarrow\infty}\sup_{s\in B_{\delta}(s^*)\cap V_i}f_{n}(s)- \max\bigg(f^i(s^*),f(s^*)\bigg)  \geq \lim_{m\rightarrow \infty}f_{n_m}(s_m)-\max\bigg(f^i(s^*),f(s^*)\bigg)
\end{split}
\end{equation}
Next, noting that $f^j(s^*)=f(s^*)$, $s^*\in A^{j,+}_\infty$ and $s^*\in A^{i,+}_\infty$, we see that the right hand side of the above is greater than or equal to:
\begin{equation}\nonumber
\begin{split}
& \lim_{m\rightarrow \infty}f_{n_m}(s_m)-\max\bigg(\sup_{s \in A^{i,+}_\infty}f^i(s),\sup_{s \in A^{j,+}_\infty}f^j(s)\bigg) \geq \lim_{m\rightarrow \infty}f_{n_m}(s_m)-\sup_{k\in\mathcal{I}} \sup_{s \in A^{k,+}_\infty}f^k(s) = R-L 
\end{split}
\end{equation}
Applying limits as $\delta$ tends to zero, we have that:
\begin{equation}\nonumber
    \begin{split}
    K_{(i,s^*)} - K'_{(i,s^*)} & = \lim_{\delta \downarrow 0}\limsup_{n\rightarrow\infty}\sup_{s\in B_{\delta}(s^*)\cap V_i}
    f_{n}(s)- \max\bigg(f^i(s^*),f(s^*)\bigg)  \geq R - L  > 0 \\
    \end{split}
\end{equation}
Thus, $K_{(i,s^*)}> K'_{(i,s^*)}$ as desired. The statements at the end of the theorem follow by identical logic to that in the proof of Theorem \ref{thm:pntwise_props}. This completes the proof.
\end{proof}

It is often the case that the suprema in question attain neither the upper nor the lower bound. This fact is illustrated by the following example.
\begin{examplebox}
\begin{example}\label{example:mid_bdd}
    Let $ \mathcal{S}=[0,1]$, $V_1=[0,\sfrac{1}{2})$ and $V_2=[\sfrac{1}{2},1]$. Define:
    \begin{equation}\nonumber
    f_n(s)= \begin{cases} 1-(2s)^n & \text{if }s \in V_1, \\ 0 & \text{if }s \in V_2, \\
    \end{cases} \quad \text{ and } \quad 
    f(s)= \begin{cases} 1 & \text{if }s \in V_1, \\ 0 & \text{if }s \in V_2. \\
    \end{cases} 
    \end{equation}
    Here we have that at $s=\sfrac{1}{2}$, $f(s)=0<1=f^1(s)$ and $f_n\rightarrow f$ pointwise but not uniformly. Further, in this example $f$ is a restrained bound for $\{f_n\}$, and thus $f_n\xrightarrow{r} f$. However, if we let $A_n=[\sfrac{1}{(2)^{1+\sfrac{1}{n}}},1]$ and $A_\infty=[\sfrac{1}{2},1]$, we have that $A_\infty^{1,+}=\{\sfrac{1}{2}\}$ and $A_\infty^{2,+}=A^{-}_\infty=A_\infty$ but;
    \begin{equation}\nonumber \lim_{n\rightarrow \infty}\sup_{s\in A_n} f_n(s) \geq \lim_{n\rightarrow \infty}f_n\bigg(\frac{1}{2^{1+\frac{1}{n}}}\bigg) = \frac{1}{2}, \quad\text{ and } 
    \end{equation}
    \begin{equation}\nonumber 
         \sup_{s \in A^-_\infty} f(s) = 0 < \frac{1}{2} < 1=\sup_{i\in\mathcal{I}}\sup_{s \in A_\infty^{+,i}}f^i(s). 
    \end{equation}
    
\centering
    \includegraphics[width=0.55\textwidth]{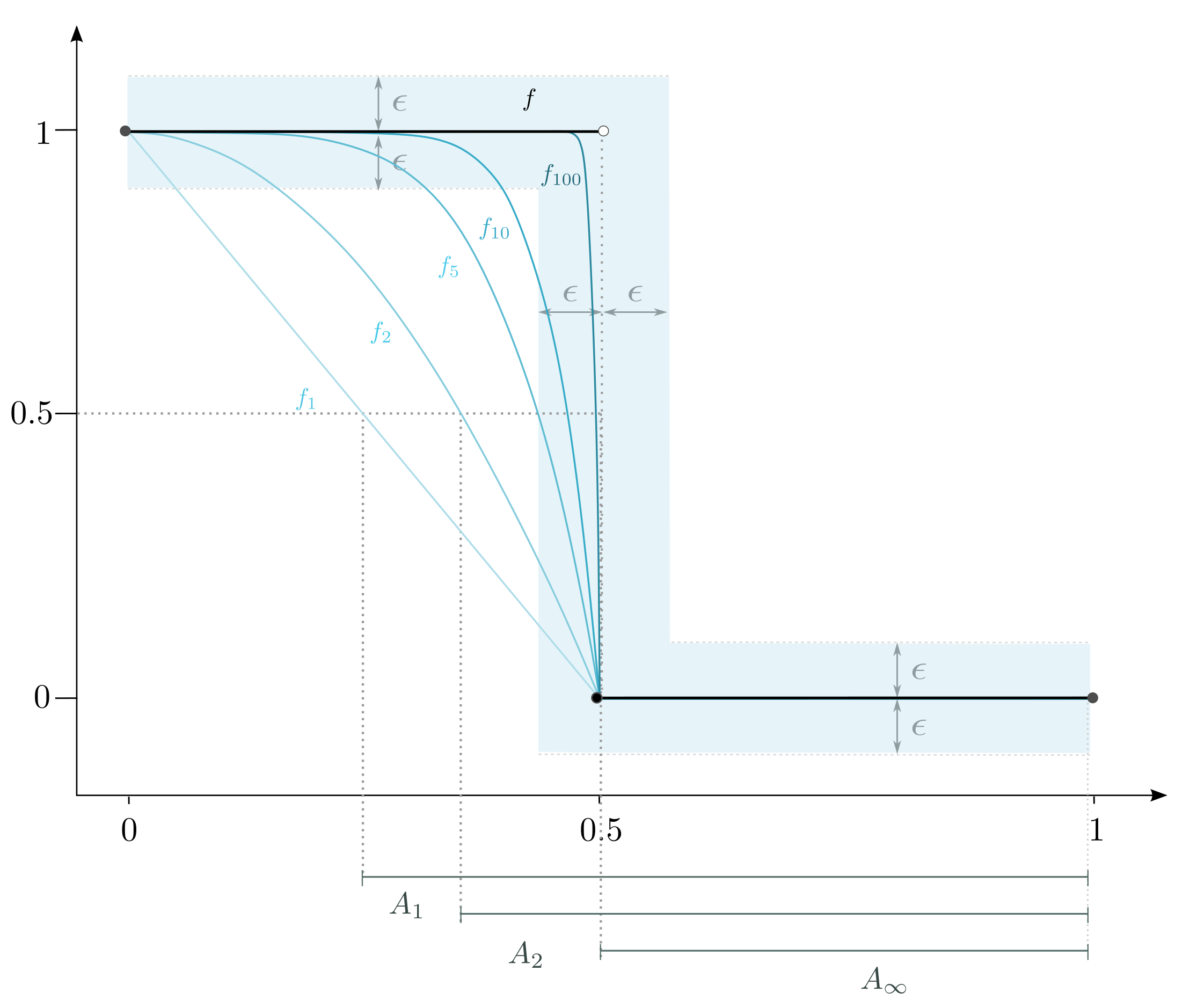}
    \captionof{figure}{An illustration of $f$ and $\{f_n\}$ for $n=1,2,5,10$ and $100$ in Example \ref{example:mid_bdd}. Also shown is a shaded blue region representing the  $\epsilon$ thickening of $\Gamma_f$ within which $f_n$ must eventually lie. The sets $A_1, A_2$ and $A_\infty$ are shown as intervals and it can be seen that, for finite $n$, $\sup_{A_n}f_n=\frac{1}{2}$.}
    \label{fig:conv_fn}
\end{example}
\end{examplebox}

%% file: Article/suprema.tex
We now have nearly all that we need to move to the study of random processes. The final step is to describe how we can construct a sequence of functions which converge with restraint, by using sequences of functions which converge uniformly.  As in Section \ref{sec:preliminaries}, we let $P(\{V_i\})$ denote the set of piecewise continuous functions on $\mathcal{S}$ that are subordinate to $\{V_i\}$.

Our goal now is to treat uniformly convergent sequences in each $\ell^\infty(\mathcal{S})$ as ``building blocks" for constructing functions which converge with restraint. As a convention, in the following text we employ the shorthand $l^\infty_k(\mathcal{S}):=\prod_{j=1}^k\ell^\infty(\mathcal{S})$\nomenclature{$l^\infty_k(\mathcal{S})$}{Shorthand for $\prod_{i=1}^k\ell^\infty(\mathcal{S})$} and denote elements of $l^\infty_k(\mathcal{S})$ with bold font; $\mathbf{f}=(f^1,...,f^k)$\nomenclature{$\mathbf{f}$}{Shorthand for $(f^1,...,f^k)\in l^\infty_k(\mathcal{S})$}. When discussing sequences $\mathbf{f}_n\rightarrow \mathbf{f}$ in $l^\infty_k(\mathcal{S})$, the convergence is assumed to be uniform on $\mathcal{S}$, unless stated otherwise. Whenever discussing mappings from $l^\infty_k(\mathcal{S})$ to $\ell^\infty(\mathcal{S})$, we shall use a tilde symbol (e.g. $\tilde{H}$)\nomenclature{$\tilde{H}$}{The $\sim$ symbol indicates that the function $\tilde{H}$ maps from $l^\infty_k(\mathcal{S})$ into $l^\infty(\mathcal{S})$} for the mapping itself and bold $\mathbf{F}\subseteq l^\infty_k(\mathcal{S})$ or $\mathbf{F}_n\subseteq l^\infty_k(\mathcal{S})$ to represent the domain of the mapping (e.g. $\tilde{H}:\mathbf{F}\rightarrow l^\infty(\mathcal{S})$). Such mappings will be used to to build sequences of functions on $\mathcal{S}$ with desirable convergence properties. This idea is formalised by the below definition.

\begin{mdframed} 
    \begin{definition}[Suprema Preserving]\label{def:SupremaPreserving}
    Let $\mathbf{F}\subseteq \ell^\infty_k(\mathcal{S})$ and $\mathbf{F}_n \subseteq \ell^\infty_k(\mathcal{S})$ for $n\in \mathbb{N}$.
    We say $\tilde{H}:\mathbf{F}\rightarrow P(\{V_i\})$ is \textit{suprema preserving} for $\tilde{H}_n:\mathbf{F}_n\rightarrow \ell^\infty(\mathcal{S})$, if for any convergent sequence $\mathbf{f}_n\rightarrow \mathbf{f}$ with $\mathbf{f}_n\in \mathbf{F}_n$ and $\mathbf{f}\in \mathbf{F}$, $\tilde{H}(\mathbf{f})$ is a restrained upper bound for $\tilde{H}_n\big(\mathbf{f}_n\big)$. Similarly, we say that it is \textit{infinima preserving} if  $\tilde{H}(\mathbf{f})$ is a restrained lower bound for $\tilde{H}_n\big(\mathbf{f}_n\big)$, and \textit{extrema preserving} if it is both suprema and infinima preserving.
    \end{definition}
    \vspace*{0.3cm}
\end{mdframed}

\begin{remark} 
    If $\tilde{H}$ is suprema preserving for $\tilde{H}_n$, we may also say that `$\tilde{H}$ preserves the suprema of $\tilde{H}_n$', or simply `$\tilde{H}, \tilde{H}_n$ are suprema preserving'. The reason for this naming convention is illustrated by the following corollary, which is a direct consequence of Theorems \ref{thm:pntwise_props} and \ref{thm:res_prop}, and the definition of suprema preserving.
\end{remark}
\begin{mdframed}
    \begin{corollary}[Properties of Suprema Preservation]\label{corr:sup_pres} Suppose that $\mathbf{f_n}\rightarrow \mathbf{f}$ with $\mathbf{f}_n \in \mathbf{F}_n$ and $\mathbf{f} \in \mathbf{F}$.  Let $\{A_k\}$ be an arbitrary sequence of sets and $A_\infty^-$ and $A_\infty^{i,+}$ be defined as in Theorems \ref{thm:pntwise_props} and \ref{thm:res_prop}, respectively. If $\tilde{H}_n(\mathbf{f_n})\rightarrow \tilde{H}(\mathbf{f})$ pointwise, then:
        \begin{equation}\nonumber
            \begin{split}
                \liminf_{n\rightarrow\infty}\bigg(\sup_{s\in A_n}\tilde{H}_n(\mathbf{f_n})(s)\bigg) \geq \sup_{s \in A_{\infty}^{-}}\tilde{H}(\mathbf{f})(s) \quad \text{and}
                            \end{split}
        \end{equation}
    if $\tilde{H}$ is suprema preserving for $\tilde{H}_n$ then:
        \begin{equation}\nonumber
            \begin{split}
                \limsup_{n\rightarrow\infty}\bigg(\sup_{s\in A_n}\tilde{H}_n(\mathbf{f_n})(s)\bigg) \leq \sup_{i\in\mathcal{I}}\sup_{s \in A^{i,+}_{\infty}}\tilde{H}(\mathbf{f})(s).
            \end{split}
        \end{equation}
    \end{corollary}
    \vspace*{0.25cm}
\end{mdframed}

At this point, it is natural to ask what type of functions are suprema preserving. The next theorem provides a partial answer to this question by giving a set of sufficient conditions for $\tilde{H}$ to be suprema preserving for $\tilde{H}_n$. Although apparently notationally dense, these conditions will greatly simplify our workings in Section \ref{sec:applications} when we turn to applications of our theory.

\begin{mdframed}
    \begin{theorem}[Sufficient Conditions for Suprema Preservation]\label{prop:suprema_charcter_sr} Let $\mathbf{F},\mathbf{F}_n\subseteq l^\infty_k(\mathcal{S})$, $\tilde{H}_n:\mathbf{F}_n\rightarrow \ell^\infty(\mathcal{S})$, $\tilde{H}:\mathbf{F}\rightarrow P(\{V_i\})$ and assume that $\tilde{H}$ has continuous extension to $l^\infty(\mathcal{S})$. Suppose the following conditions hold for all convergent $\mathbf{f}_n\rightarrow \mathbf{f}$ with $\mathbf{f}_n\in\mathbf{F}_n$ and $\mathbf{f}\in \mathbf{F}$:
        		\begin{equation}\nonumber
                \begin{split}
                (i)&~~\text{There exists $C>0$ and $\alpha\in\mathbb{R}$ such that for all sufficiently large $n$,} \\
                &~~\sup_{s \in \mathcal{S}}\tilde{H}(\mathbf{f}_n)(s) - \tilde{H}(\mathbf{f})(s) \leq C  \|  \mathbf{f}_n - \mathbf{f} \|_\infty^\alpha,\\
                (ii)&~~\text{For each $i \in \mathcal{I}$ and $s \in \overline{V_i}$.}\\
                &~~ 
            		\lim_{\delta \downarrow 0}\limsup_{n\rightarrow\infty}
            		\sup_{s'\in B_{\delta}(s)\cap V_i} \tilde{H}_n(\mathbf{f}_n)(s') - \tilde{H}(\mathbf{f}_n)(s')
            		\leq \max\bigg(0, \tilde{H}(\mathbf{f})(s) - \tilde{H}^{i}(\mathbf{f})(s)\bigg),\\
                    & ~~\text{where $\tilde{H}^i(\mathbf{f})$ denotes the continuous extension of $\tilde{H}(\mathbf{f})|_{V_i}$ to $\overline{V_i}$\nomenclature{$\tilde{H}^i(\mathbf{f})$}{Continuous extension of $\tilde{H}(\mathbf{f})|_{V_i}$ to $\overline{V_i}$}.}
            	\end{split}
        		\end{equation} 
        Then $\tilde{H}$ is suprema preserving for $\tilde{H}_n$.
    \end{theorem}
    \vspace*{0.25cm}
\end{mdframed}
\begin{remark} A drawback of the definition of suprema preserving is that, to verify it holds, one must relate the global behaviour of $\tilde{H}_n(\mathbf{f}_n)$ to that of $\tilde{H}(\mathbf{f})$ (see Definitions \ref{def:res_def} and \ref{def:SupremaPreserving}). The above theorem partially relieves this difficulty by providing conditions that are almost predominantly defined on the restrictions of $\tilde{H}_n(\mathbf{f}_n)$ and $\tilde{H}(\mathbf{f})$ to a given partition element $\overline{V_i}$. Such conditions are desirable as, in many practical settings, assumptions of local uniform convergence and continuity may be reasonable within each $\overline{V_i}$ and thus the above conditions might be easier to verify (see Section \ref{sec:applications} for examples). 

     It must be emphasized here that the claim \textit{``$\tilde{H}$ is suprema preserving for $\tilde{H}_n$''} in Theorem \ref{prop:suprema_charcter_sr} is a statement about $\tilde{H}$ as a mapping on $\mathbf{F}$ and not a statement about it's continuous extension to $l^\infty_k(\mathcal{S})$. It is not, in general, true that the conditions of Theorem \ref{prop:suprema_charcter_sr} imply that the \textit{continuous extension of $\tilde{H}$ to $l^\infty_k(\mathcal{S})$} is suprema preserving for $\tilde{H}_n$.
\end{remark}

\begin{proof} 
	Choose
	$i \in \mathcal{I}$ and $s \in \overline{V_i}$, let $j\in\mathcal{I}$ be such that $s \in V_j$, let $\delta>0$ and assume Condition $(i)$ holds. Now, by the triangle inequality for suprema and the fact that a supremum increases if the set over which it is taken gets larger, we have:
\begin{equation}\nonumber
      \sup_{s' \in B_{\delta}(s)\cap V_i}\tilde{H}_n(\mathbf{f}_n)(s') \leq \sup_{s'\in  B_{\delta}(s)\cap V_i}
      \tilde{H}_n(\mathbf{f}_n)(s') - \tilde{H}(\mathbf{f}_n)(s')
\end{equation}
\begin{equation}\nonumber
      + \sup_{s'\in \mathcal{S}} \tilde{H}(\mathbf{f}_n)(s') - \tilde{H}(\mathbf{f})(s') + \sup_{s'\in B_{\delta}(s)\cap V_i}\tilde{H}^{i}(\mathbf{f})(s').
\end{equation}
Taking the limit superior of both sides and noting that the second term on the right hand side is smaller than zero in the limit by condition $(i)$, we have that:
\begin{equation}\nonumber
      \limsup_{n\rightarrow \infty}\sup_{s'\in B_{\delta}(s)\cap V_i}\tilde{H}_n(\mathbf{f}_n)(s') \leq \limsup_{n\rightarrow \infty}\sup_{s'\in B_{\delta}(s)\cap V_i}
      \tilde{H}_n(\mathbf{f}_n)(s') - \tilde{H}(\mathbf{f}_n)(s') + \sup_{s'\in B_{\delta}(s)\cap V_i}\tilde{H}^i(\mathbf{f})(s').
\end{equation}
Next, applying the limit as $\delta\downarrow 0$ on both sides, we see that by the continuity of $\tilde{H}^i(\mathbf{f})$:
\begin{equation}\nonumber
      \lim_{\delta\downarrow 0}\limsup_{n\rightarrow \infty}\sup_{s'\in B_{\delta}(s)\cap V_i}\tilde{H}_n(\mathbf{f}_n)(s') \leq \lim_{\delta\downarrow 0}\limsup_{n\rightarrow \infty}\sup_{s'\in B_{\delta}(s)\cap V_i}
      \tilde{H}_n(\mathbf{f}_n)(s') - \tilde{H}(\mathbf{f}_n)(s') + \tilde{H}^i(\mathbf{f})(s).
\end{equation}
Subtracting $\max(\tilde{H}^i(\mathbf{f})(s),\tilde{H}(\mathbf{f})(s))$ from both sides we now see that:
	\begin{equation}\nonumber
	\begin{split}
      & \lim_{\delta\downarrow 0}\limsup_{n\rightarrow \infty}\sup_{s'\in B_{\delta}(s)\cap V_i}\tilde{H}_n(\mathbf{f}_n)(s') - \max\bigg(\tilde{H}^i(\mathbf{f})(s),\tilde{H}(\mathbf{f})(s)\bigg) \\
      & \leq \lim_{\delta \downarrow 0}\limsup_{n\rightarrow\infty}\sup_{s'\in B_{\delta}(s)\cap V_i}
      \tilde{H}_n(\mathbf{f}_n)(s') - \tilde{H}(\mathbf{f}_n)(s')- \max\bigg(0,\tilde{H}(\mathbf{f})(s) - \tilde{H}^{i}(\mathbf{f})(s)\bigg). \\
    \end{split}
    \end{equation}
    Thus, the condition stated in Definition \ref{def:res_bdd1} follows if:
    	\begin{equation}\nonumber
	\begin{split}
\lim_{\delta \downarrow 0}\limsup_{n\rightarrow\infty}\sup_{s'\in B_{\delta}(s)\cap V_i}
      \tilde{H}_n(\mathbf{f}_n)(s') - \tilde{H}(\mathbf{f}_n)(s')- \max\bigg(0,\tilde{H}(\mathbf{f})(s) - \tilde{H}^{i}(\mathbf{f})(s)\bigg) \leq 0\,,
    \end{split}
    \end{equation}
    which is precisely condition $(ii)$.
\end{proof}

\begin{remark} 
    If we want to verify that $\tilde{H}$ is both suprema and infimina preserving for $\tilde{H}_n$, then validating Condition $(i)$ for both cases is equivalent to showing that:
    \begin{equation*}
        ||\tilde{H}(\mathbf{f}_n) - \tilde{H}(\mathbf{f})||_\infty=\mathcal{O} \left( \|  \mathbf{f}_n - \mathbf{f} \|_\infty^\alpha
        \right).
    \end{equation*}
    In other words, it is equivalent to showing that the extension of $\tilde{H}$ to $l^\infty(\mathcal{S})$ is $\alpha-$H\"{o}lder continuous. It should also be noted here that Condition $(i)$ can in fact be weakened: the proof above only requires
    \begin{equation}\nonumber
        \limsup_{n\to\infty}\sup_{s'\in \mathcal{S}} 
        \big(\tilde{H}(\mathbf{f}_n)(s') - \tilde{H}(\mathbf{f})(s')\big) \leq 0.
    \end{equation}
    Nevertheless, in practice the stronger condition is often more convenient for analysis.
\end{remark}

It might also be of interest to ask whether the sum of two suprema preserving functions is again suprema preserving. The answer is, unfortunately, no - Example \ref{example:symdiff4} implicitly provides a counterexample, although some work is required to see this (see the proof of Theorem \ref{thm:symdiff_conv_with_res}in Supplementary Section \ref{supp:symdiff_conv_with_res} for further detail). However, the following theorem does provide a sufficient condition under which the sum of two suprema preserving functions is suprema preserving. This theorem will be invaluable to the proofs of Section \ref{sec:symdiff}. 

\begin{mdframed} 
    \begin{corollary}[Sum of Suprema Preserving Functions]\label{cor:SumSupPreserving} Let $\{U_i\}$ and $\{V_j\}$ be partitions of $\mathcal{S}$ and define $\mathcal{N}$ as in Theorem \ref{thm:SumUpperRestraint}.
    Suppose $\tilde{K}:\mathbf{F}\rightarrow P(\{U_i\})$ and $\tilde{L}:\mathbf{G}\rightarrow P(\{V_j\})$ are suprema
    preserving for $\tilde{K}_n:\mathbf{F}_n\rightarrow \ell^\infty(\mathcal{S})$ and  $\tilde{L}_n:\mathbf{G}_n\rightarrow \ell^\infty(\mathcal{S})$, respectively. If for all $s\in \mathcal{N}$, $f\in \mathbf{F}$ and $g\in \mathbf{G}$, it holds that
	\begin{equation}\label{eq:cond_sum_suprema_preserv}
	\begin{split}
		\max\Big(\tilde{K}^i(f)(s) + \tilde{L}(g)(s), &\tilde{K}(f)(s) + \tilde{L}^j(g)(s) \Big)\\
		&\leq \max\left(\tilde{K}^i(f)(s) + \tilde{L}^j(g)(s), \tilde{K}(f)(s) + \tilde{L}(g)(s) \right),
	\end{split}
	\end{equation}
	then $\tilde{K} + \tilde{L}: \mathbf{F} \times \mathbf{G} \rightarrow \ell^\infty(\mathcal{S})$ is suprema preserving for $\tilde{K}_n + \tilde{L}_n: \mathbf{F}_n \times \mathbf{G}_n \rightarrow \ell^\infty(\mathcal{S})$. 
	\end{corollary}
    \vspace*{0.3cm}
\end{mdframed}
\begin{proof}
    The result follows directly from the definition of suprema preserving and Theorem \ref{thm:SumUpperRestraint}.
\end{proof}

As we have seen, the notion of suprema preserving functions is useful as it allows us to describe convergence properties of piecewise continuous sequences of functions. This idea can be naturally extended to the setting of random processes, in which we have a probability space $(\Omega, \mathcal{F},\mathbb{P})$, and are concerned with sequences of random processes $\hat{H}_n:\Omega\rightarrow l^\infty(\mathcal{S})$ that `converge with restraint' to a piecewise continuous process $H: \Omega\rightarrow P(\{V_i\})$. To formalise the notion of convergence with restraint for random processes, we now provide the following definition.

\begin{mdframed}
\begin{definition}[Restrained Convergence of Random Variables]\label{def:sup_pres_rand}Let $\hat{H}_n:\Omega \rightarrow l^\infty(\mathcal{S})$, $H:\Omega \rightarrow P(\{V_i\})$ be such that $\hat{H}_n(\mathcal{S})$ converges weakly to $H(s)$ for all $s \in \mathcal{S}$. Further, suppose $\hat{H}_n$ and $H$ can be represented in the form $\hat{H}_n=\tilde{H}_n\circ \mathbf{\hat{G}_n}$ and $H=\tilde{H}\circ \mathbf{G}$, respectively, where, for some $\mathbf{F},\mathbf{F}_n\subseteq l^\infty_k(\mathcal{S})$:
\begin{itemize}
    \item $\tilde{H}:\mathbf{F}\rightarrow P(\{V_i\})$ is suprema preserving for $\tilde{H}_n:\mathbf{F}_n\rightarrow \ell^\infty(\mathcal{S})$,
    \item $\mathbf{\hat{G}_n}:\Omega\rightarrow \mathbf{F}_n$ and $\mathbf{G}:\Omega\rightarrow \mathbf{F}$ satisfy a (functional) weak convergence of the form $\mathbf{\hat{G}_n}\xrightarrow{d}\mathbf{G}$\nomenclature{$\mathbf{G}$}{Random variable taking its values in $\ell^\infty_k(\mathcal{S})$}\nomenclature{$\mathbf{G}_n$}{Random variable taking its values in $\ell^\infty_k(\mathcal{S})$ converging in distribution to $\mathbf{G}$}. 
\end{itemize}
Then, we say that $\hat{H}_n$ converges to $H$ with upper restraint and write $\hat{H}_n\xrightarrow{u.r.} \tilde{H}$\nomenclature{$\xrightarrow{u.r.}$}{Convergence with upper restraint of random variables}. If $\tilde{H}_n$ is instead infinima preserving for $\tilde{H}$, we say that $\hat{H}_n$ converges with lower restraint to $H$ and write $\hat{H}_n\xrightarrow{l.r.}H$. Similarly, if $\tilde{H}_n$ and $\tilde{H}$ are extrema preserving, we say that $\hat{H}_n$ converges with restraint to $H$ and write $\hat{H}_n\xrightarrow{r}H$.
\end{definition}
\vspace*{0.25cm}
\end{mdframed} 
\begin{remark} Another way of stating the above is to say that $\hat{H}_n$ converges with restraint to $H$ if and only if $\hat{H}_n(s)\xrightarrow{d}H(s)$ for all $s\in \mathcal{S}$ and the following diagrams commute,
\begin{equation}\nonumber
\begin{tikzcd}
\Omega \arrow[r, "\hat{H}_n"] \arrow[d, "\mathbf{\hat{G}}_n"'] 
  & \ell^\infty(\mathcal{S}) \\
\mathbf{F}_n \arrow[ur, "\tilde{H}_n"'] &
\end{tikzcd} \quad \quad \quad
\begin{tikzcd}
\Omega \arrow[r, "H"] \arrow[d, "\mathbf{G}"'] 
  & P(\{V_i\}) \\
\mathbf{F} \arrow[ur, "\tilde{H}"'] &
\end{tikzcd} 
\end{equation}
for function spaces $\{\mathbf{F}_n\}_{n\in\mathbb{N}}$ and $\mathbf{F}$, $\mathbf{\hat{G}_n}\xrightarrow{d} \mathbf{G}$ and $\tilde{H}$ being suprema preserving for $\tilde{H}_n$. The function spaces $\mathbf{F}$ and $\mathbf{F}_n$ have been introduced in order to state our theory in full generality. However, in the following sections, we shall conventionally take $\mathbf{F}_n=l^\infty_k(\mathcal{S})$ for all $n\in\mathbb{N}$ and $\mathbf{F}=C(\mathcal{S},\mathbb{R})^k$ for some appropriate choice of $k \in \mathbb{N}$. 

\end{remark}

Informally, Definition \ref{def:sup_pres_rand} may be thought of as stating that $H$ and $\hat{H}_n$ converge with restraint if they can be `built' from uniformly convergent variables $\mathbf{\hat{G}_n}=(\hat{G}^1_n,...,\hat{G}^k_n)\rightarrow\mathbf{G}=(G^1,...,G^k)$ in a `reasonable' way. Explicit examples of the above construction are given in Section \ref{sec:applications}. 

For some applications, it may be the case that $\hat{H}_n$ does not converge with restraint, but can instead be bounded above and below by random processes which do converge with restraint. This notion is captured by the following definition. We shall see an explicit example of its use in Section \ref{sec:symdiff}. 

\begin{mdframed}
    \begin{definition}[Confinement of a Process]\label{def:conf_pair} Let $\hat{L}_n,\hat{H}_n,\hat{U}_n:\Omega\rightarrow l^\infty(\mathcal{S})$ and $L,U:\Omega\rightarrow P(\{V_i\})$ be random processes. Suppose that $\forall s \in \mathcal{S},~\hat{L}_n(s)\leq\hat{H}_n(s)\leq\hat{U}_n(s)$ almost surely, $\hat{L}_n\xrightarrow{l.r.}L$ and $\hat{U}_n\xrightarrow{u.r.}U$. Then we say that $(\hat{L}_n,\hat{U}_n)$ is a \textit{confinement} of $\hat{H}_n$ with limit $(L,U)$.
    \end{definition}
    \vspace{0.2cm}
\end{mdframed}
\begin{remark}
    When using the above definition, we shall often simply say that $(\hat{L}_n,\hat{U}_n)$ is a confinement of $\hat{H}_n$, with the implicit understanding that in such instances the restrained limit of $\hat{L}_n$ is denoted by $L$ and the restrained limit of $\hat{U}_n$ is $U$. In such instances, we also assume that, if it exists, the pointwise limit of $\hat{H}_n$ is $H$; that is, for each $s \in \mathcal{S}$, $\hat{H}_n(s)\xrightarrow{d} H(s)$ in the usual univariate sense.
\end{remark}

We are now ready to state our first theorem concerning restrained convergence of random variables, Theorem \ref{thm:random_converge}. 

\begin{mdframed}
    \begin{theorem}\label{thm:random_converge} Suppose $\hat{H}_n:\Omega\rightarrow l^\infty(\mathcal{S})$ converges pointwise with confinement $(\hat{L}_n,\hat{U}_n)$. Let $\{A_n\}$ and $\{B_n\}$ be sequences of subsets of $\mathcal{S}$, with $A_\infty^-,B_\infty^-,A_\infty^{i,+}$ and $B^{i,+}_\infty$ defined as in Theorems \ref{thm:pntwise_props} and \ref{thm:res_prop}. Then, for arbitrary $q\in \mathbb{R}$, we have:
   \begin{equation}\nonumber
            \begin{split}
                \limsup_{n\rightarrow\infty}\mathbb{P}_*\bigg[\max\Bigg(& \sup_{s\in A_n} -\hat{H}_n(s), \sup_{s\in B_n}\hat{H}_n(s) \Bigg) \leq q \bigg]\\
                & \leq \mathbb{P}\bigg[\max\Bigg(\sup_{s \in A^{-}_\infty}-H(s), \sup_{s \in B^{-}_\infty}H(s) \Bigg)\leq q\bigg]\,, \quad \text{and}
            \end{split}
        \end{equation}
        \begin{equation}\nonumber
            \begin{split}
                \liminf_{n\rightarrow\infty}\mathbb{P}^*\bigg[\max\Bigg(& \sup_{s\in A_n} -\hat{H}_n(s), \sup_{s\in B_n}\hat{H}_n(s) \Bigg) < q\bigg]\\
                & \geq  \mathbb{P}\bigg[\sup_{i\in \mathcal{I}}\max\Bigg(\sup_{s \in A^{i,+}_\infty}-L^i(s),
                \sup_{s \in B^{i,+}_\infty}U^i(s) \Bigg)< q\bigg]\,.
            \end{split}
        \end{equation}
    \end{theorem}
    \vspace*{0.25cm}
\end{mdframed}
\begin{proof}  The first claim follows directly from Theorem \ref{thm:pntwise_props} combined with the continuous mapping theorem given by Supplemental Lemma \ref{lem:cts_mapping}. For the second claim, by the definition of convergence with restraint, we can write $\hat{L}_n=\tilde{L}_n(\mathbf{\hat{G}_n})$ and $\hat{U}_n=\tilde{U}_n(\mathbf{\hat{G}'_n})$ for some  $\tilde{L}_n:\mathbf{F}_n\rightarrow l^\infty(\mathcal{S})$ and $\tilde{U}_n:\mathbf{F}'_n\rightarrow l^\infty(\mathcal{S})$ with suprema preserving limits. Let $x_n\rightarrow x$ and $y_n\rightarrow y$ be convergent sequences in $\{\mathbf{F}_n\}_{n\in\mathbb{N}}$ and $\{\mathbf{F}_n'\}_{n\in\mathbb{N}}$ with limits $x \in \mathbf{F}$ and $y \in \mathbf{F}'$, respectively, and define the following functions:
    \begin{equation*}
    \begin{split}
        g_n(x_n,y_n)&= \max\Big( \sup_{s\in A_n} -\tilde{L}_n(x_n)(s), \sup_{s\in B_n} \tilde{U}_n(y_n)(s) \Big),\\
        g^+(x,y)
& = \sup_{i\in\mathcal{I}}\max\Big(\sup_{s\in A^{i,+}_\infty}-\tilde{L}(x)(s), \sup_{s\in B^{i,+}_\infty}\tilde{U}(y)(s)\Big)
    \end{split}
    \end{equation*} 
    Corollary \ref{corr:sup_pres} combined with Supplemental Lemma \ref{lem:max_min_limsup} now shows that $\limsup_{n \rightarrow \infty} g_n(x_n,y_n)\leq g^+(x,y)$. Applying Supplemental Lemma \ref{lem:cts_mapping} it follows that:
        \begin{equation}\nonumber
            \begin{split}
                \liminf_{n\rightarrow\infty}\mathbb{P}^*\bigg[\max\Bigg(& \sup_{s\in A_n} -\hat{L}_n(s), \sup_{s\in B_n}\hat{U}_n(s) \Bigg) < q\bigg]\\
                & \geq  \mathbb{P}\bigg[\sup_{i\in\mathcal{I}}\max\Bigg(\sup_{s \in A^{i,+}_\infty}-L^i(s),
                \sup_{s \in B^{i,+}_\infty}U^i(s) \Bigg)< q\bigg]\,.
            \end{split}
        \end{equation} 
    The result now follows from the definition of confinement.
\end{proof}

If $\hat{H}_n\xrightarrow{r}H$ then it is trivially true that $(\hat{H}_n,\hat{H}_n)$ is a confinement of $\hat{H}_n$ with limit $(H,H)$. Using this fact, we obtain the below corollary.
\begin{mdframed}
    \begin{corollary}
       If $\hat{H}_n\xrightarrow{r}H$ then the second bound of Theorem \ref{thm:random_converge} may be replaced with:
        \begin{equation}\nonumber
            \begin{split}
                \liminf_{n\rightarrow\infty}\mathbb{P}^*\bigg[\max\Bigg(& \sup_{s\in A_n} -\hat{H}_n(s), \sup_{s\in B_n}\hat{H}_n(s) \Bigg) < q\bigg]\\
                & \geq  \mathbb{P}\bigg[\sup_{i\in\mathcal{I}}\max\Bigg(\sup_{s \in A^{i,+}_\infty}-H^i(s),
                \sup_{s \in B^{i,+}_\infty}H^i(s) \Bigg)< q\bigg]\,.
            \end{split}
        \end{equation}    
    \end{corollary}
    \vspace{0.2cm}
\end{mdframed}
\begin{remark}
    For convenience, in Section \ref{sec:applications}, we shall assume that $\mathcal{I}$ is finite. As a result, for the applied examples discussed therein we shall write $\max_{i \in \mathcal{I}}$ instead of $\sup_{i \in \mathcal{I}}$ when employing the above theorem and corollary.
\end{remark}
To conclude this section, we introduce one final definition: the \textit{restrained Central Limit Theorem (rCLT)}. In the next section, we return to our motivating context, constructing CRs for random processes that satisfy rCLTs. 
\begin{mdframed}
    \begin{definition}[Restrained CLT] Suppose $\hat{\mu}_n:\Omega\rightarrow l^\infty(\mathcal{S})$ is an estimator of function $\mu\in l^\infty(\mathcal{S})$ satisfying a pointwise CLT of the form:
    \begin{equation}\nonumber
        \hat{H}_n(s):=\tau_n^{-1}(\hat{\mu}_n(s)-\mu(s))\xrightarrow{d} H(s),
    \end{equation}
    for some $\hat{H}_n:\Omega\rightarrow l^\infty(\mathcal{S})$, $H:\Omega\rightarrow P(\{V_i\})$ and some sequence of positive constants $\tau_n\rightarrow 0$.
    We say that $\hat{\mu}_n$ satisfies a \textit{restrained Central Limit Theorem (rCLT)} if, in addition, $\hat{H}_n\xrightarrow{r}H$. Alternatively, if $\hat{H}_n$ has a given confinement, we say that $\hat{\mu}_n$ satisfies a \textit{pointwise Central Limit Theorem (pCLT) with confinement}.
    \end{definition}
    \vspace{0.2cm}
\end{mdframed}

%% file: Article/crs.tex
In the previous sections, we developed a theory of restrained convergence for piecewise continuous random processes. In this section, we use this theory to construct confidence regions for supra- and super-threshold sets of random processes that satisfy either a restrained CLT or a pointwise CLT with confinement.
     
That is, throughout this section, we let
$\hat \mu_n:\Omega\rightarrow l^\infty(\mathcal{S})$ be an estimator of some spatially varying unknown function $\mu:\mathcal{S}\rightarrow\mathbb{R}$, let $\tau_n$ be a positive sequence with $\tau_n \rightarrow 0$, and define $\hat{H}_n := \tau_n^{-1}(\hat{\mu}_n-\mu)$. The pointwise limit of $\hat{H}_n$ is denoted $H$, and assumed to be a piecewise continuous process subordinate to $\{V_i\}$. We emphasize here that our theory differs from the existing CR literature as we do not assume an fCLT for $\hat{H}_n$ converging to $H$, only an rCLT or pCLT with confinement. This difference, although subtle, allows the theory to be applied in a wider range of practical scenarios, as will be illustrated in Section \ref{sec:applications}. It should also be noted that the only assumption of continuity employed in this section is that of $H$ being a piecewise continuous process. 

As discussed in Section \ref{sec:intro}, our aim is to generate lower and upper CRs, $\hat{\mathcal{L}}$ and $\hat{\mathcal{U}}$, for the sets:
\begin{equation}\nonumber
   \mathcal{L}=\{s \in \mathcal{S}: \mu(s)<0\} \quad \text{ and } \quad  \mathcal{U}=\{s \in \mathcal{S}: \mu(s)>0\},
\end{equation}
such that $\mathbb{P}[\,\hat{\mathcal{L}} \subseteq \mathcal{L}~\land~\hat{\mathcal{U}} \subseteq \mathcal{U}\,]\approx 1-\alpha$ as $n\rightarrow\infty$ for some pre-specified confidence level, $\alpha\in[0,1]$. This is analogous to confidence intervals in one-dimensional settings, for which a common convention is to use $\alpha=0.05$. 

However, defining CRs in this way presents difficulties: 
the probability above may fail to converge as $n \rightarrow \infty$, or worse, the inclusion event, $\{\hat{\mathcal{L}}\subseteq\mathcal{L}~\wedge~\hat{\mathcal{U}}\subseteq\mathcal{U}\}$, may not even be measurable. 
For these reasons, we instead aim to construct CRs that satisfy
\begin{equation}\label{eq:cr_def}
    \liminf_{n\rightarrow\infty}\mathbb{P}_*\big[\,\hat{\mathcal{L}} \subseteq \mathcal{L}~ \wedge ~\hat{\mathcal{U}} \subseteq \mathcal{U}\,\big] \geq 1-\alpha_1 \quad \text{ and } \quad \limsup_{n\rightarrow\infty}\mathbb{P}^*\big[\,\hat{\mathcal{L}} \subseteq \mathcal{L}~ \wedge ~\hat{\mathcal{U}} \subseteq \mathcal{U}\,\big] \leq 1-\alpha_2
\end{equation}
for some known $\alpha_1, \alpha_2 \in [0,1]$. 
This formulation provides inner and outer probability bounds on the 
inclusion event without requiring the event to be measurable or its limiting probability to exist. 
In the special case where measurability holds and the limit exists, 
\eqref{eq:cr_def} implies the probability of the inclusion event lies in 
$[1-\alpha_1,1-\alpha_2]$. 

To ensure \eqref{eq:cr_def} holds, we shall construct CRs according to the rule
\begin{equation}\nonumber
   \hat{\mathcal{L}}:=\{s \in \mathcal{S}: \tau_n^{-1}\hat{\mu}_n(s)<-q\} \quad \text{ and } \quad  \hat{\mathcal{U}}:=\{s \in \mathcal{S}: \tau_n^{-1}\hat{\mu}_n(s)>+q\},
\end{equation}
where $q \in \mathbb{R}$ will be chosen so as to ensure that at least one, if not both, of the bounds in \eqref{eq:cr_def} hold. To relate the construction above to \eqref{eq:cr_def}, and thus develop a method of choosing the value of $q$, we intend to employ the proof strategy of the Simultaneous Confidence Region Excursion (SCoRE) set metatheorem of \cite{telschow2023scope}. Specifically, adapting this framework to the context of piecewise continuous processes yields the following theorem:

\begin{mdframed} 
    \begin{theorem}[Confidence Regions for Piecewise Continuous Functions]\label{thm:scope_for_piecewise}
       Suppose $\hat{\mu}_n:\Omega\rightarrow l^\infty(\mathcal{S})$ satisfies a pCLT of the form $\hat{H}_n(s):=\tau_n^{-1}(\hat{\mu}_n(s)-\mu(s))\xrightarrow{d}H(s)$ with confinement $(\hat{L}_n, \hat{U}_n)$. Define: \begin{equation}\nonumber
            \mathfrak{u}^+_i(0) = \bigcap_{\eta>0}\overline{ \overline{\mu^{-1}\big([0, \eta]\big)} \cap V_i}\quad \text{ and } \quad \mathfrak{u}^-_i(0) = \bigcap_{\eta>0}\overline{ \overline{\mu^{-1}\big([-\eta, 0]\big)} \cap V_i}.
        \end{equation}
        The below statements now hold:
			        \begin{equation}\nonumber
			        \begin{split}
			            (i)\quad \liminf_{n\rightarrow\infty}
			            	\mathbb{P}_*\big[\,\hat{\mathcal{L}} \subseteq \mathcal{L}~ \wedge ~\hat{\mathcal{U}} \subseteq \mathcal{U}\,\big]
			            \geq \mathbb{P}\bigg[\sup_{i\in\mathcal{I}}\max\bigg(\sup_{s \in  \mathfrak{u}^+_i(0)} \!\!\!-L^i(s),\sup_{s \in \mathfrak{u}^-_i(0)} U^i(s)\bigg)<q\bigg],
			                    \end{split}
			        \end{equation}
            	    \begin{equation}\nonumber  
	        		\begin{split}
           				(ii)\quad \limsup_{n\rightarrow\infty}\mathbb{P}^*\big[\,\hat{\mathcal{L}} \subseteq \mathcal{L}~ \wedge ~\hat{\mathcal{U}} \subseteq \mathcal{U}\,\big]
           				\leq \mathbb{P}\bigg[\sup_{s \in  \mu^{-1}(0)} |H(s)| \leq q\bigg].
           			\end{split}
        			\end{equation}
    \end{theorem}
    \vspace*{0.3cm}
\end{mdframed}
\begin{remark} If $\hat{\mu}_n$ satisfies an rCLT, rather than a pCLT with confinement, the right hand side of (i) simplifies to
    \begin{equation}\nonumber
        \mathbb{P}\bigg[\sup_{i\in\mathcal{I}}\max\bigg(\sup_{s \in  \mathfrak{u}^+_i(0)} \!\!\!-H^i(s),\sup_{s \in \mathfrak{u}^-_i(0)} H^i(s)\bigg)<q\bigg].
    \end{equation}
    This follows directly from the definitions of confinement and convergence with restraint for random variables (see Definitions \ref{def:sup_pres_rand} and \ref{def:conf_pair}).
    If in addition, the sets $\mathfrak{u}^+_i(0)$ and $\mathfrak{u}^-_i(0)$ happen to be equal, further simplifications can be made by replacing the $H^i(s)$ and $-H^i(s)$ terms above with $|H^i(s)|$, thus highlighting the similarity between the right hand sides of $(i)$ and $(ii)$.  We shall explore conditions under which such simplifications arise in specific applications 
    in Section \ref{sec:applications} (see Corollaries \ref{corr:conds_abs}, \ref{cor:conj_conditions}, 
    and \ref{corr:conds_sd}).  
\end{remark}
\begin{proof}
    First, we prove $(i)$. To do so, we first show that if the inclusions $\hat{\mathcal{L}}\subseteq \mathcal{L}$ and $\hat{\mathcal{U}}\subseteq \mathcal{U}$ hold,  then $\sup_{s\in\mu^{-1}(0)}\vert \hat H_n(s) \vert \leq q$. Hence, suppose $s\in\mu^{-1}(0)$ satisfies $\vert \hat H_{n}(s) \vert > q$.  It follows that $\tau_n^{-1}(\hat{\mu}_n(s)-\mu(s))=\tau_n^{-1}\hat{\mu}_n(s)>q$ and thus $\hat\mu_n(s)>\tau_n q$. Therefore, $s\in\hat{\mathcal{U}}$. As $\mathcal{U}\cap\mu^{-1}(0)$ is empty, we have $\hat{\mathcal{U}}\not\subseteq \mathcal{U}$. A similar argument shows that, if there is an $s\in\mu^{-1}(0)$ such that $\hat H_n(s) < -q$ then $\hat{\mathcal{L}}\not\subseteq \mathcal{L}$.
    Consequentially, we have that 
    \begin{equation}\nonumber  
		\begin{split}
           	\limsup_{n\rightarrow\infty}\mathbb{P}^*\big[\,\hat{\mathcal{L}} \subseteq \mathcal{L}~ \wedge ~\hat{\mathcal{U}} \subseteq \mathcal{U}\,\big]
           				\leq            	\limsup_{n\rightarrow\infty}\mathbb{P}^*\bigg[\,\sup_{s \in  \mu^{-1}(0)} |\hat H_{n}(s)| \leq q\,\bigg]
           				\leq
                        \mathbb{P}\bigg[\sup_{s \in  \mu^{-1}(0)} |H(s)| \leq q\bigg]\,,
    	\end{split}
    \end{equation}
    where the last inequality follows from the first statement in Theorem \ref{thm:random_converge}. 

    We turn to the proof of $(ii)$. Assume that $\eta_n$ is a sequence such that $\tau_n^{-1}\eta_n \rightarrow \infty$ for $n$ tending to $\infty$. First, note that
    \begin{equation}\label{eq:CopeInclusion}
        \begin{split}
            \sup_{s\in \mu^{-1}([-\eta_n,0])\cap V_i} \hat{H}_n(s) < q~~\wedge~~ \sup_{s\in \mathcal{S}} \vert \hat{H}_n(s) \vert < q + \tau_n^{-1}\eta_n ~~\Longrightarrow~ ~\hat{\mathcal{U}} \subseteq \mathcal{U}\,.
        \end{split}
    \end{equation}
    To see this, let $s\in \big( \mathcal{S}\setminus\mathcal{U}\big) \cap \mu^{-1}([-\eta_n,0])\cap V_i$. This means that $q > \hat H_n(s) > \tau_n^{-1}\hat\mu_n(s)$ as $-\eta_n \leq \mu(s) \leq 0$. Thus, $s\in \mathcal{S} \setminus \hat{\mathcal{U}}$. If instead $s\in \big( \mathcal{S}\setminus\mathcal{U}\big) \cap \mathcal{S} \setminus  \big(\mu^{-1}([-\eta_n,0])\cap V_i \big)$, then $q + \tau_n^{-1}\eta_n > \hat{H}_n(s) \geq \tau_n^{-1}( \hat\mu_n(s) + \eta_n)$ as $\mu(s)\leq -\eta_n$.  From this, we have that $\hat\mu_n(s) < q\tau_n$ and thus $s\in \mathcal{S} \setminus \hat{\mathcal{U}}$. Putting the arguments together, we see that the left hand side of \eqref{eq:CopeInclusion} implies $\mathcal{U}^c\subseteq \hat{\mathcal{U}}^c$ and therefore $\hat{\mathcal{U}} \subseteq \mathcal{U}$. A similar argument shows that
    \begin{equation}\label{eq:CopeInclusion2}
        \begin{split}
            \sup_{s\in \mu^{-1}([0,\eta_n])\cap V_i} -\hat{H}_n(s) < q~~\wedge~~ \sup_{s\in \mathcal{S}} \vert \hat{H}_n(s) \vert < q + \tau_n^{-1}\eta_n ~~\Longrightarrow~~\hat{\mathcal{L}} \subseteq \mathcal{L}.
        \end{split}
    \end{equation}
    From \eqref{eq:CopeInclusion} and \eqref{eq:CopeInclusion2} we obtain
    \begin{equation}\nonumber
    \begin{split}
            \mathbb{P}_*\big[\,\hat{\mathcal{L}} \subseteq \mathcal{L}~ \wedge ~\hat{\mathcal{U}} \subseteq \mathcal{U}\,\big]
        \geq
            \mathbb{P}_*\bigg[
\max\Big( \sup_{s\in \mu^{-1}([0,\eta])\cap V_i} -\hat{H}_n(s)&,
 \sup_{s\in \mu^{-1}([-\eta,0])\cap V_i} \hat{H}_n(s)\Big) < q         
            \,\bigg]\\ &+ 
            1 -  \mathbb{P}_*\bigg[
            \sup_{s\in \mathcal{S}} \vert \hat{H}_n(s) \vert < q + \tau_n^{-1}\eta_n          
            \,\bigg]\,.
    \end{split}
    \end{equation}
    Applying the limit inferior to both sides yields that the first summand on the right hand side is bounded from below by 
    \begin{equation}\nonumber
        \begin{split}  
        \mathbb{P}\bigg[\sup_{i\in\mathcal{I}}\max\bigg(\sup_{s \in  \mathfrak{u}^+_i(0)} \!\!\!-L^i(s),\sup_{s \in \mathfrak{u}^-_i(0)} U^i(s)\bigg)<q\bigg]
                \end{split}
    \end{equation}
    by the second statement of Theorem \ref{thm:random_converge}. The remaining terms on the right hand side converge to zero. To see this, we show that $\sup_{s \in \mathcal{S}}|\hat{H}_n(s)|$ is asymptotic tight in the sense of \cite[p.21]{van1996weak}. Fix $\varepsilon > 0$ be arbitrary, and let $\kappa > 0$ be such that:
	\begin{equation}\label{eq:tightness_bound}
	\begin{split}
		\mathbb{P}\bigg[\sup_{i\in\mathcal{I}}\sup_{s \in \bar V_i} \big\vert H^i(s) \big\vert < \kappa\bigg]
				\geq 1-\varepsilon\,.
	\end{split}
	\end{equation}		  The above is possible as any random variable on $\mathbb{R}$ is tight. Thus, as $\tau_n^{-1}\eta_n\rightarrow \infty$, for $n$ large enough we have that
	\begin{equation}\nonumber
	   \liminf_{n\rightarrow\infty} \mathbb{P}_*\bigg[
			\sup_{s \in \mathcal{S}}  \big\vert \hat{H}_n(s) \vert < q + \tau_n^{-1}\eta_n \bigg]\\
			\geq
		\liminf_{n\rightarrow\infty} \mathbb{P}_*\bigg[
			\sup_{s \in \mathcal{S}}  \big\vert \hat{H}_n(s) \big\vert < \kappa \bigg]
	\end{equation}	
	\begin{equation}\nonumber
			\geq
			\mathbb{P}\bigg[\sup_{i\in\mathcal{I}}\sup_{s \in  \bar V_i} \big\vert H^i(s) \big\vert < \kappa\bigg]
			\geq 1-\varepsilon\,.
	\end{equation}	
	Here the second inequality follows from Theorem \ref{thm:random_converge}
	applied to	$A_n = B_n = \mathcal{S}$ and the last inequality follows from \eqref{eq:tightness_bound}.
\end{proof}

\begin{remark} 
    The strength of Theorem \ref{thm:scope_for_piecewise} is that it can be used to generate CRs. For instance, part $(ii)$ shows that if we set $q$ to be the $(1-\alpha_2)\%$ quantile of $\sup_{\mu^{-1}(0)}|H|$, then the (outer) probability of inclusion is asymptotically bounded above by $(1-\alpha_2)\%$. An analogous statement holds for the lower bound in $(i)$.  Consequently, to construct CRs with target coverage in $[1-\alpha_1,1-\alpha_2]$, it suffices to estimate the distribution of the suprema statistics appearing inside the probabilities on the right hand sides of $(i)$ and $(ii)$. In practice, such estimation is typically achieved via bootstrap methods. As these procedures are well established in the literature, we omit further discussion here and refer the reader to Section \ref{sec:background_crs} for additional details and references.
\end{remark}

%% file: Article/applications.tex
In this section, we illustrate the utility of our theorems for constructing CRs in a range of practical examples. Although the examples differ substantively from one another, each can be modelled within our framework using the same general approach. Broadly, this approach proceeds as follows.

First, we formulate the problem in terms of an estimator $\hat{\mu}_n:\Omega \rightarrow l^\infty(\mathcal{S})$ of some target function $\mu:\mathcal{S} \rightarrow \mathbb{R}$. Next, we define $\hat{H}_n:=\tau_n^{-1}(\hat{\mu}_n-\mu)$ and identify its pointwise (distributional) limit $H(s)$. We then express $\hat{H}_n$ and $H$ in the forms $\hat{H}_n:=\tilde{H}_n(G_n^i)$ and $\hat{H}:=\tilde{H}(G^i)$, where $\mathbf{G_n}\xrightarrow{d}\mathbf{G}$ and $\tilde{H}_n$ and $\tilde{H}$ are known fixed functions. The domains of $\tilde{H}_n$ and $\tilde{H}$ will be $\mathbf{F}_n:=l_k^\infty(\mathcal{S})$, for all $n$, and $\mathbf{F}=C(\mathcal{S},\mathbb{R})^k$ throughout, where $k$ is a known constant that may vary between examples. Using these representations, we establish either an rCLT for $\hat{\mu}_n$, by showing that $\hat{H}_n\xrightarrow{r} H$, or a pCLT with confinement, by deriving a confinement for $\hat{H}_n$. We shall show these convergence results by verifying that Conditions $(i)$ and $(ii)$ of Theorem \ref{prop:suprema_charcter_sr} are satisfied. Finally, we apply Theorem \ref{thm:scope_for_piecewise} to construct the desired CRs.

During the final step, computational methods are employed to approximate the distribution of the suprema appearing on the right hand side of Theorem~\ref{thm:scope_for_piecewise} $(i)$ and $(ii)$. Based on these distributions, we select the quantile $q$ to achieve the desired coverage probability.  Since our aim here is to illustrate the theoretical framework, we do not discuss the computational details of quantile estimation but refer the interested reader to \cite{SSS, Bowring:2019a} for further discussion. To facilitate computation, it is often useful to establish conditions under which the sets $\mathfrak{u}^{\pm}_i(0)$ can be simplified and interpreted in a meaningful way. We also address this question within the context of each example. 

Before presenting the examples, we offer a final remark on the continuity assumptions underlying this work. In Theorem~\ref{thm:scope_for_piecewise} of Section~\ref{sec:piecewise_copes}, no continuity is assumed for $\mu$, $\hat{\mu}_n$, or $\hat{H}_n$. The only requirement is that $H$ is a piecewise continuous process on $\mathcal{S}$, typically constructed from a collection of continuous random processes $\{G^i\}$. In the examples that follow, we shall employ one additional assumption: that $\mu$ is continuous. Although this assumption is not strictly necessary for applying Theorem~\ref{thm:scope_for_piecewise}, continuity assumptions of this form are commonplace in the applied literature and, in the following examples, enable the derivation of tighter confidence bounds. 

\subsection{Absolute Value}\label{sec:absolute}

We begin with the canonical example of a piecewise process, discussed at the beginning of the document in Section \ref{sec:background_pieceproc}; the absolute value. A simple motivating example for this process is given below:
\begin{examplebox}
\begin{example}[Uncertainty in Decision Boundary Estimation] 
Suppose a researcher is interested in a classification problem where two predictors $X_1$ and $X_2$ are used to predict a binary outcome $Y$. A standard approach in this setting is to fit a logistic regression model, yielding an estimate of $\mathbb{P}[\,Y=1\,|\,X_1=x_1,X_2=x_2\,]$ given by
\begin{equation}\nonumber
    \hat{p}_n(x_1,x_2) := \frac{1}{1+e^{-\hat{f}_n(x_1,x_2)}},
\end{equation}
where $\hat{f}_n(x_1,x_2)$ is a linear combination of features whose coefficients, $\{\hat{\beta}_i\}$, are estimated from a dataset containing $n$ observations. In the simplest case $\hat{f}_n(x_1,x_2) = \hat{\beta}_0 + \hat{\beta}_1 x_1 + \hat{\beta}_2 x_2$, but more flexible analyses may incorporate polynomial expansions, interaction terms, or basis spline functions of $x_1$ and $x_2$ to capture nonlinear effects. Given a new pair of predictors $(x_1,x_2)$, the researcher predicts $Y=1$ if $\hat{p}_n(x_1,x_2)>0.5$ and $Y=0$ otherwise. This rule corresponds to the decision boundary shown in Fig.~\ref{fig:class} (a).

Now suppose classification errors carry a high cost, as in medical diagnosis when $Y$ indicates the presence of a disease. Both misclassifying a healthy patient as needing treatment and failing to detect the disease in an affected patient incur significant penalties. In such cases, it may be preferable to account for uncertainty in the estimated decision boundary $\hat{p}_n(x_1,x_2)=0.5$ and abstain from predicting when confidence is low, perhaps instead recommending further tests to be conducted.

        \vspace{0.2cm}
    \begin{centering}\includegraphics[width=0.85\textwidth]{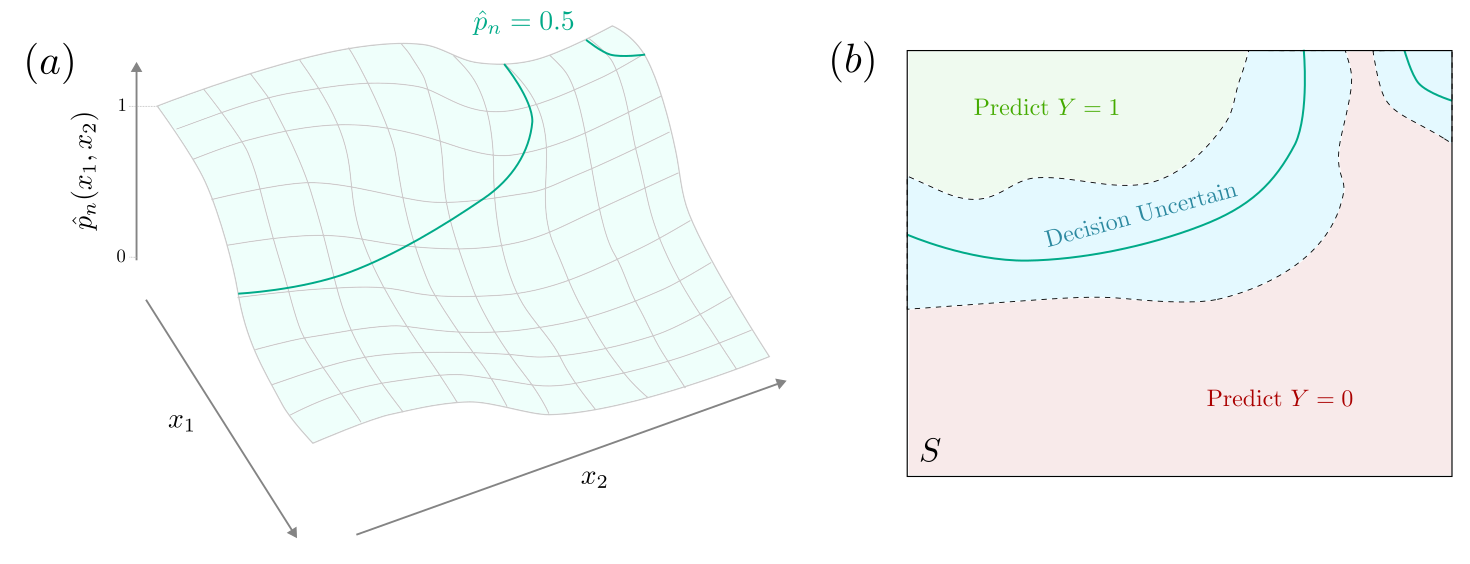}
    \captionof{figure}{(a) Estimated probability $\hat{p}_n$ as a function of $x_1$ and $x_2$, with the decision boundary $\hat{p}_n(x_1,x_2)=0.5$ shown in green. (b) Top-down view of the decision boundary with a confidence region (CR) shaded in blue. Predictions are: $Y=1$ in the green region, $Y=0$ in the red region, and uncertain within the blue CR (where more data may be needed).}\label{fig:class}
    \end{centering}
    
\end{example}
\end{examplebox}

In this example, we aim to quantify the spatial uncertainty of the estimated contour $\hat{p}_n(s) = 0.5$, where $s := (x_1, x_2)$ lies in a compact region $\mathcal{S}$ and $\hat{p}_n$ is an estimator of some unknown true $p\in l^\infty(\mathcal{S})$. To do so, we first define $\gamma(s) := p(s) - 0.5$ and $\mu:=|\gamma|$, so that the decision boundary corresponds to the set $\mathcal{U}^c$, where $\mathcal{U} := \{ s \in \mathcal{S} : \mu(s) > 0 \}$. Our goal is to construct a CR for $\mathcal{U}^c$ by first constructing a CR, $\hat{\mathcal{U}}$, for $\mathcal{U}$ and noting that $\hat{\mathcal{U}} \subseteq \mathcal{U}$ if and only if $\mathcal{U}^c \subseteq \hat{\mathcal{U}}^c$ (that is, we wish to generate the blue CR shown in Fig.~\ref{fig:class}(b)).

To construct a CR for the absolute value of $\gamma : \mathcal{S} \to \mathbb{R}$, we shall assume $\gamma$ is continuous and consider the estimator $\hat{\gamma}_n := \hat{p}_n - 0.5$. Note that the setting described is now identical to that discussed in Section~\ref{sec:background_pieceproc}. 
Following the approach outlined there, we partition $\mathcal{S}$ according to the sign of $\gamma$, defining
\begin{equation}\nonumber
    V_i := \{ s \in \mathcal{S} : \mathrm{sgn}(\gamma(s)) = i \}, \quad \text{for}\quad i \in \mathcal{I}:=\{-1, 0, 1\}.
\end{equation}
To generate CRs, we shall consider the process $\hat{H}_n := \tau_n^{-1} \big( |\hat{\gamma}_n| - |\gamma| \big)$
where $\tau_n=1/\sqrt{n}$. Assuming an fCLT of the form $\hat{G}_n := \tau_n^{-1} (\hat{\gamma}_n - \gamma) \xrightarrow{d} G$, where $G$ is a continuous random process, we see that
the pointwise limit of $\hat{H}_n$ is given by the piecewise process
\begin{equation}\nonumber
    H(s) :=
    \begin{cases}
        G(s) & \text{if } s \in V_1, \\[6pt]
        -G(s) & \text{if } s \in V_{-1}, \\[6pt]
        |G(s)| & \text{if } s \in V_0.
    \end{cases}
\end{equation}
With these definitions in place, we can now state the following result.

\begin{mdframed}
\begin{theorem}[RCLT for the Absolute Value]\label{thm:abs_val_res}
    Given the definitions above, $\hat{\mu}_n =\vert \hat\gamma_n \vert$ satisfies a restrained CLT. That is, $\hat{H}_n\xrightarrow{r}H$.
\end{theorem}
\vspace{0.2cm}
\end{mdframed}

\begin{proof}
    See Supplementary Material Section \ref{supp:abs_val_res}.
\end{proof}

By applying Theorem \ref{thm:scope_for_piecewise}, we may now generate CRs for the absolute value as follows.
\begin{mdframed}
\begin{theorem}[Confidence Regions for the Absolute Value]\label{thm:abs_val_crs} Suppose $\gamma:\mathcal{S}\rightarrow \mathbb{R}$ and $G,\hat{\gamma}_n:\Omega \rightarrow l^\infty(\mathcal{S})$ are defined as above and $\mathcal{U}$ and $ \hat{\mathcal{U}}$ are given by:
\begin{equation}\nonumber
    \mathcal{U}:=\{s \in \mathcal{S}: |\gamma(s)|>0\} \quad \text{ and }\quad \hat{\mathcal{U}}:=\{s \in \mathcal{S}:\tau_n^{-1}|\hat{\gamma}_n(s)|>+q\}.
\end{equation}
Then:
\begin{equation}\nonumber
    (i)\quad\liminf_{n\rightarrow\infty}
        \mathbb{P}_*\big[~\hat{\mathcal{U}} \subseteq \mathcal{U}\,\big]
    \geq\mathbb{P}\bigg[\sup_{s \in  \gamma^{-1}(0)} \vert G(s) \vert <q\bigg],\quad 
    \end{equation}
    \begin{equation}\nonumber (ii)\quad\limsup_{n\rightarrow\infty}\mathbb{P}^*\big[~\hat{\mathcal{U}} \subseteq \mathcal{U}\,\big]
    \leq \mathbb{P}\bigg[\sup_{s \in  \gamma^{-1}(0)} \vert G(s) \vert \leq q\bigg].
\end{equation}
\end{theorem}
\vspace{0.2cm}
\end{mdframed}

\begin{proof}
    See Supplementary Material Section \ref{supp:abs_val_crs}.
\end{proof}

By employing standard properties of the inner and outer measure, the following corollary is immediate.
\begin{mdframed}
    \begin{corollary}\label{corr:conds_abs} Suppose $\gamma:\mathcal{S}\rightarrow \mathbb{R},~G,\hat{\gamma}_n:\mathcal{S}\times \Omega,~\mathcal{U}$ and $ \hat{\mathcal{U}}$ are defined as above. If $\mathbb{P}[\sup_{s\in\gamma^{-1}(0)}|G(s)|=q]=0$, then: 
        \begin{equation}\nonumber \lim_{n\rightarrow\infty}\mathbb{P}\big[~\hat{\mathcal{U}} \subseteq \mathcal{U}\,\big]
    = \mathbb{P}\bigg[\sup_{s \in  \gamma^{-1}(0)} \vert G(s) \vert \leq q\bigg].
\end{equation}
\vspace{0.2cm}
    \end{corollary}
\end{mdframed}

It is important that we acknowledge that many variants of Theorem \ref{thm:abs_val_crs} are already known in the CR literature. For example, Corollary~1 of \cite{SSS} provides a similar, though not identical, result under additional assumptions of continuity on $\hat{\gamma}_n$ and the topology of $\partial \mathcal{U}$. The result may also be derived using the results of  \cite{telschow2023scope}, which removes the topological assumptions of \cite{SSS} and substitutes the continuity assumption for the more general condition of quantile convergence of $G$. Beyond these references, CRs for contour estimation have also been studied extensively in \cite{Lindgren1995,Polfeldt1999,French2014,Bolin2017} and \cite{Aquino2024} (of particular relevance is Equation~(3) of \cite{French2014}).

However, the results of this section differ from the existing literature in two key ways. First, Theorem~\ref{thm:abs_val_crs} requires only an fCLT for $G$ and continuity of $\gamma$ (we conjecture that this latter assumption may also be removed using confinements. However, we shall not consider this extension in detail here) making it applicable in a wider range of settings. Second, and more importantly, we establish the novel convergence result $\hat{H} \xrightarrow{r} H$ for $H$ defined using absolute value functions. This is not just a theoretical curiosity: it is crucial for the developments that follow, particularly in proofs employed in Section~\ref{sec:symdiff}, where more complex functions involving sums of absolute values are analysed.

\subsection{Intersections and Unions}\label{sec:conjunc}

A notable application of our theory is to the setting in which excursion sets are derived from multiple images and researchers wish to compare how such excursion sets overlap or differ from one another. This concept is made explicit by the following example.

\begin{examplebox}
    \begin{example}[Climate Conjunction Analysis] 
    Suppose a climate researcher is interested in identifying geographical regions that exhibit abnormally large estimated temperature increases, consistent across multiple studies. Specifically, they aim to compare images of temperature increase drawn from $m$ different studies (illustrated for $m=2$ in Fig.~\ref{fig:example8}) in order to identify regions where the temperature increase exceeds a predefined threshold $c$ in \textbf{all studies}.

    Let $\hat{\gamma}^i_n(s)$ denote the estimated temperature increase at location $s$ in study $i$. Formally, the researcher is interested in the overlap of the excursion sets $\{s \in \mathcal{S} : \hat{\gamma}^i_n(s) > c\}$ across all studies, for a given $c \in \mathbb{R}$. This overlapping region, where all studies agree the threshold is exceeded, is referred to as a \textbf{conjunction}.
        \vspace{0.2cm}

    \begin{centering}
    \includegraphics[width=\textwidth]{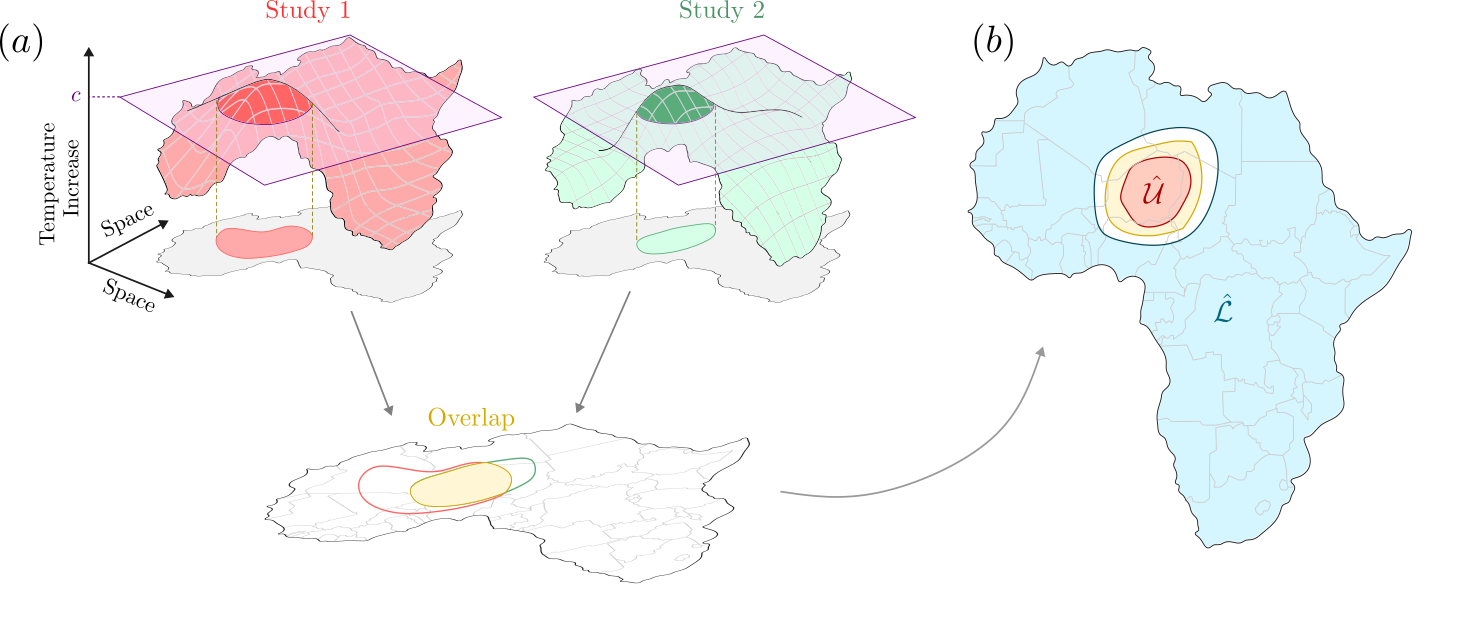}
    \captionof{figure}{$(a)$ Top: estimated temperature increase maps from two studies, represented as red and green surfaces, thresholded at level $c$. Bottom: boundaries of the excursion sets $\{s \in \mathcal{S} : \hat{\gamma}^i_n(s) > c\}$ for $i \in \{1,2\}$, with their conjunction (overlap) highlighted in yellow. $(b)$ Top-down view of the conjunction set in $(a)$, with corresponding CRs, $\hat{\mathcal{L}}$ and $\hat{\mathcal{U}}$, shaded in blue and red, respectively.}
    \label{fig:example8}
    \end{centering}
    \end{example}
\end{examplebox}

In the above example, the researcher possesses $m$ spatially-varying estimates, $\hat{\gamma}_n^i:\Omega \rightarrow l^\infty(\mathcal{S})$ for $i\in\mathcal{M}:=\{1,...,m\}$, of $m$ unknown, true, signals $\gamma^i\in C(\mathcal{S},\mathbb{R})$. As in previous section, we assume a fCLT of the following form:
\begin{equation}\nonumber
    \hat{J}_n^i:=\tau_n^{-1}(\hat{\gamma}_n^i-\gamma^i)\xrightarrow{d}J^i
\end{equation}
where the fCLT is understood to hold jointly over $i$ and $\mathbf{\hat{J}_n}$ and $\mathbf{J}$ take their values in $l^\infty_m(\mathcal{S})$ and $C(\mathcal{S}, \mathbb{R})^m$, respectively. Empirical interest now lies in identifying where all $m$ true signals exceed some predetermined threshold $c\in\mathbb{R}$. Alternatively, this can be written in terms of the minimum function, $\min_i(\gamma^i)$. Asking whether $\min_i(\gamma^i)>c$ is equivalent to asking where \text{all} functions exceed $c$. Without loss of generality, in the following we assume that $c=0$. Our aim is to generate CRs for the sets $\mathcal{U}:=\{s \in \mathcal{S}: \min_i\gamma^i(s)>0\}$ and $\mathcal{L}:=\{s \in \mathcal{S}: \min_i\gamma^i(s)<0\}$ (c.f. Fig.~\ref{fig:example8} (b)). Following our previous conventions, in this example, we define $\mu:\mathcal{S}\rightarrow \mathbb{R}$ as $\mu(s):=\min_{i \in \mathcal{M}} \gamma^i(s)$.

The most direct approach to estimating $\mu=\min_i(\gamma^i)$ is to use $\min_i(\hat{\gamma}_n^i)$. However, as in the previous examples, 
\begin{equation}\label{eq:Hn_def_min}
    \hat{H}_n(s):=\tau_n^{-1}\bigg(\min_{i\in\mathcal{M}}\hat{\gamma}_n^i(s)-\min_{i\in\mathcal{M}}\gamma^i(s)\bigg)
\end{equation}
does not satsify an fCLT. Instead, by direct computation $\hat{H}_n$ can be seen to converge pointwise to a piecewise process. To express this limiting process, we define the shorthand for $\sigma\subseteq \mathcal{M}$: 
\begin{equation}\nonumber
    \gamma^\sigma:=\min_{j \in \sigma}\gamma^j\quad \text{ and }\quad V_\sigma:=\{s \in \mathcal{S}: \gamma^\mathcal{M}(s)=\gamma^j(s) ~\text{if and only if}~ j \in \sigma\}. 
\end{equation}
In words, $V_\sigma$ is spatial region over which the smallest $\{\gamma^j\}$ are exactly those for which $j\in\sigma$. Another way of saying this is that $s\in V_\sigma$ if and only if ``$\min_i(\gamma^i(s))=\gamma^j(s)$ for all $j\in\sigma$". By direct computation, it can be seen that the pointwise limiting process is given by:
\begin{equation}\label{eq:H_def_min}
    H(s) := \min_{i \in \sigma} J^i(s)  \quad \text{ if }\quad s \in V_\sigma.
\end{equation}
Note that, in this example, the partition $\{V_\sigma\}$ is indexed by $\mathcal{I}:=\mathcal{P}(\mathcal{M})$, where $\mathcal{P}$ represents the powerset. We emphasize here that there is no requirement for the above process to be continuous across the boundaries $\{\partial V_\sigma\}$, as the minimum may jump discretely when indices are added to or removed from $\sigma$. We are now able to state the following convergence result.

\begin{mdframed}
\begin{theorem}[RCLT for the Minimum Estimator]\label{thm:min_example} 
    Given the definitions above, $\hat{\mu}_n=\min_{i \in \mathcal{M}} \hat\gamma^i_n$ satisfies a restrained CLT. That is, $\hat{H}_n\xrightarrow{r}H$.
\end{theorem}
\vspace*{0.3cm}
\end{mdframed}

\begin{proof}
    See Supplementary Material Section \ref{supp:min_example}.
\end{proof}

We are now in a position to obtain CRs for this setting as follows.
\begin{mdframed}
\begin{theorem}[Confidence Regions for the Conjunction Example]\label{thm:min_crs} 
Suppose $\mu,\gamma^i:\mathcal{S}\rightarrow \mathbb{R}$, $\hat{\gamma}_n^i,\hat{J}^i_n:\Omega\rightarrow l^\infty(\mathcal{S})$ and $J^i:\Omega\rightarrow C(\mathcal{S},\mathbb{R})$ are defined as above and $\mathcal{L}, \hat{\mathcal{L}}, \mathcal{U}$ and $ \hat{\mathcal{U}}$ are given by:
\begin{equation}\nonumber
    \mathcal{L}:=\bigg\{s \in \mathcal{S}:\min_{i\in\mathcal{M}}\big(\gamma^i(s)\big)<0\bigg\},\quad\mathcal{U}:=\bigg\{s \in \mathcal{S}:\min_{i\in\mathcal{M}}\big(\gamma^i(s)\big)>0\bigg\},
\end{equation}
\begin{equation}\nonumber
    \hat{\mathcal{L}}:=\bigg\{s \in \mathcal{S}:\tau_n^{-1}\min_{i\in\mathcal{M}}\big(\hat{\gamma}_n^i(s)\big)<-q\bigg\}\quad \text{ and }\quad\hat{\mathcal{U}}:=\bigg\{s \in \mathcal{S}:\tau_n^{-1}\min_{i\in\mathcal{M}}\big(\hat{\gamma}_n^i(s)\big)>+q\bigg\}.
\end{equation}
If $\mathfrak{u}^\pm_\sigma(0)$ are defined as in the statement of Theorem \ref{thm:scope_for_piecewise}, with the index $i$ substituted for $\sigma$ throughout, and $G^\sigma:=\min_{i\in\sigma}J^i$ and $\hat{G}_n^\sigma:=\min_{i\in\sigma}\hat{J}_n^i$, then:
\begin{equation}\nonumber
\begin{split}
    (i)\quad \liminf_{n\rightarrow\infty}
        \mathbb{P}_*\big[\,\hat{\mathcal{L}} \subseteq \mathcal{L}~ \wedge ~\hat{\mathcal{U}} \subseteq \mathcal{U}\,\big]
    \geq \mathbb{P}\bigg[\max_{\sigma\in \mathcal{I}}\max\bigg(\sup_{s \in  \mathfrak{u}^+_\sigma(0)} -G^\sigma(s),\sup_{s \in \mathfrak{u}^-_\sigma(0)} G^\sigma(s)\bigg)<q\bigg],
            \end{split}
\end{equation}
\begin{equation}\nonumber  
\begin{split}
    (ii)\quad \limsup_{n\rightarrow\infty}\mathbb{P}^*\big[\,\hat{\mathcal{L}} \subseteq \mathcal{L}~ \wedge ~\hat{\mathcal{U}} \subseteq \mathcal{U}\,\big]
    \leq \mathbb{P}\bigg[\max_{\sigma\in \mathcal{I}}\sup_{s \in  \mu^{-1}(0) \cap V_\sigma} \vert G^\sigma(s) \vert \leq q\bigg].
\end{split}
\end{equation}
\end{theorem}
\vspace{0.1cm}
\end{mdframed}
\begin{proof}
    The result follows by direct application of Lemma \ref{thm:min_example} and Theorem \ref{thm:scope_for_piecewise}. 
\end{proof}

The problem of constructing CRs for `single-study' climate data has been discussed in \cite{French2013, French2015, SSS} and \cite{Hazra2021}. Similarly, `multiple-condition' conjunction sets have been a topic of investigation in the fMRI literature for several decades (c.f. \cite{Friston1999, Nichols2005}). However, few authors have addressed the spatial uncertainty inherent in conjunction sets, and the prevailing practice in fMRI remains rudimentary visual inspection, without formal statistical assessment (see, e.g., \cite{Zhang2021, Dijkstra2021, Ferreira2021}). We believe this issue has been largely overlooked in the literature. As a result, \cite{MaullinSapey2022} introduced CRs for conjunction analysis. However, the framework in \cite{MaullinSapey2022} relied upon much stronger assumptions than those employed here, including differentiability of $\hat{\gamma}^i_n$ and $\gamma$ with nonzero gradients in pre-specified regions, as well as several restrictive topological assumptions on the boundary of $\mathcal{U}$ (see Section 2.2 of \cite{MaullinSapey2022}). The present proofs thus offer a substantial simplification and extension of the earlier approach.

Under additional assumptions, the expressions provided by Theorem \ref{thm:min_crs} can be further simplified. The following corollary offers such a simplification.

\begin{mdframed}
    \begin{corollary}\label{cor:conj_conditions}
        Suppose the conditions of Theorem \ref{thm:min_crs} hold. If either of the following conditions hold:
        \begin{itemize}
            \item[\textbf{C1:}] For every $\sigma \in \mathcal{I}$, we have that $\overline{\mu^{-1}(0)\cap V_\sigma}=\mu^{-1}(0)\cap \overline{V_\sigma}$.
            \item[$\widetilde{\textbf{C1}}$\textbf{:}] For every $s\in \mu^{-1}(0)$, if $s \in \overline{V_\sigma}$, then for any sufficiently small open ball around $s$, $B_\delta(s)$ with $\delta>0$, we have that the sets $B_\delta(s)\cap \mu^{-1}([0,\infty))\cap V_\sigma$ and $B_\delta(s)\cap \mu^{-1}((-\infty,0])\cap V_\sigma$ are non-empty, and the set $B_\delta(s)\cap V_\sigma$ is path-connected.
        \end{itemize}
        Then part $(i)$ of Theorem \ref{thm:min_crs} simplifies to:
        \begin{equation}\nonumber
            \liminf_{n\rightarrow\infty}
        \mathbb{P}_*\big[\,\hat{\mathcal{L}} \subseteq \mathcal{L}~ \wedge ~\hat{\mathcal{U}} \subseteq \mathcal{U}\,\big]
    \geq \mathbb{P}\bigg[\max_{\sigma\in \mathcal{I}}\sup_{s \in  \mu^{-1}(0)\cap V_\sigma} |G^\sigma(s)|<q\bigg].
        \end{equation}
    If, in addition, the following condition holds:\begin{itemize}
  \item[\textbf{C2:}] 
  $\mathbb{P}\bigg[\displaystyle\max_{\sigma\in \mathcal{I}} \displaystyle\sup_{s \in \mu^{-1}(0) \cap V_\sigma} \big| G^\sigma(s) \big| = q \bigg] = 0$.
\end{itemize}
    then parts $(i)$ and $(ii)$ of Theorem \ref{thm:min_crs} simplify to:
    \begin{equation}\nonumber
        \lim_{n\rightarrow\infty}
        \mathbb{P}\big[\,\hat{\mathcal{L}} \subseteq \mathcal{L}~ \wedge ~\hat{\mathcal{U}} \subseteq \mathcal{U}\,\big]
    =\mathbb{P}\bigg[\max_{\sigma\in \mathcal{I}}\sup_{s \in  \mu^{-1}(0)\cap V_\sigma} |G^\sigma(s)|\leq q\bigg].
    \end{equation}
    \end{corollary}
    \vspace{0.1cm}
\end{mdframed}

\begin{proof}
    See Supplementary Material Section \ref{supp:conj_conditions}.
\end{proof}

\begin{remark}
    As shown in Supplementary Material Section \ref{supp:conj_conditions}, we have that $\widetilde{\textbf{C1}}$ implies \textbf{C1}. Although \textbf{C1} is the weaker assumption, we have opted to include both conditions here, as assumptions similar to $\widetilde{\textbf{C1}}$ are common in the CR literature. For examples, see Assumption (f1) of \cite{Cuevas2006}, Assumption 1a of \cite{SSS}, and Assumption 2.2.3a of \cite{MaullinSapey2022}. In fact, we conjecture that Assumptions 2.2.3a-c of \cite{MaullinSapey2022} imply $\widetilde{\textbf{C1}}$, and thus $\widetilde{\textbf{C1}}$ is the weaker assumption. We do not prove this claim here, but note that the non-empty intersections of $\widetilde{\textbf{C1}}$ are supersets of the sets $\mathcal{F}_0^\circ$ and $\mathcal{J}^\sigma_0$ in the notation of \cite{MaullinSapey2022}, and thus Assumption 2.2.3a of \cite{MaullinSapey2022} directly implies the statements about non-empty sets in $\widetilde{\textbf{C1}}$.
\end{remark}

We conclude this section by observing that, as the section title suggests, the above approach can be extended to construct CRs for unions (or \textit{`disjunctions'}) of excursion sets. This is done by replacing each occurrence of \texttt{max} in the preceding arguments with \texttt{min}, and carefully reversing signs and inequalities where required in the proof of Theorem~\ref{thm:min_example} to account for the asymmetry. The resulting statement parallels Theorem~\ref{thm:min_crs}, but applies to the disjunction sets
\begin{equation}\nonumber
\mathcal{L} := \{s \in \mathcal{S} : \max_{i\in\mathcal{M}} \gamma^i(s) < 0\} 
\quad \text{and} \quad 
\mathcal{U} := \{s \in \mathcal{S} : \max_{i\in\mathcal{M}} \gamma^i(s) > 0\},
\end{equation}
with $\mu := \max_{i\in\mathcal{M}} \gamma^i$. A similar idea can also be used to derive CRs for set differences of the form
\begin{equation}\nonumber
\{s \in \mathcal{S} : \gamma^1(s) > 0\} \setminus \{s \in \mathcal{S} : \gamma^2(s) > 0\},
\end{equation}
though we omit the details here. Based on these observations, it might be expected that the above proofs can be generalized to handle arbitrary set-theoretic combinations of excursion sets. However, this is unfortunately not the case for set-theoretic combinations such as unions of intersections or intersections of unions. In such cases, convergence with restraint results analogous to Theorem~\ref{thm:min_example} cannot be established, and we must instead use the notion of confinement introduced in Definition~\ref{def:conf_pair}. This issue is explored further in the following section.

\subsection{Symmetric Difference}\label{sec:symdiff}

The final application we shall consider is that of the \textit{symmetric difference}, motivated by the below example. 

\begin{examplebox}
\begin{example}[Image Comparison in fMRI]
A common aim in the analysis of fMRI data is to compare brain activity between two groups of subjects, such as healthy controls and patients with a neurological condition, to identify brain regions that are active in one group but not the other. For example, a researcher might present both groups with the same memory recall task, record their brain activity, and average the results over each group to obtain two images, each showing regions of significant activation for one group. In this setting, it is of interest to identify the active brain regions observed during the experiment which differ between the groups.    

Formally, we denote the true (unknown) spatially-varying brain activation for the two groups as $\gamma^1$ and $\gamma^2 : \mathcal{S} \to \mathbb{R}$, where $\mathcal{S}$ represents the human brain. For a predefined activation threshold $c$, we define the excursion sets $\{s \in \mathcal{S} : \gamma^i(s) > c\}$ for $i \in \{1, 2\}$, which correspond to the active regions of brain activity observed for each group.  The aim is to identify the symmetric difference of these sets, that is, the regions that are active in exactly one group but not both (see Fig.~\ref{fig:symdiffX}).

    \begin{centering}
    \includegraphics[width=\textwidth]{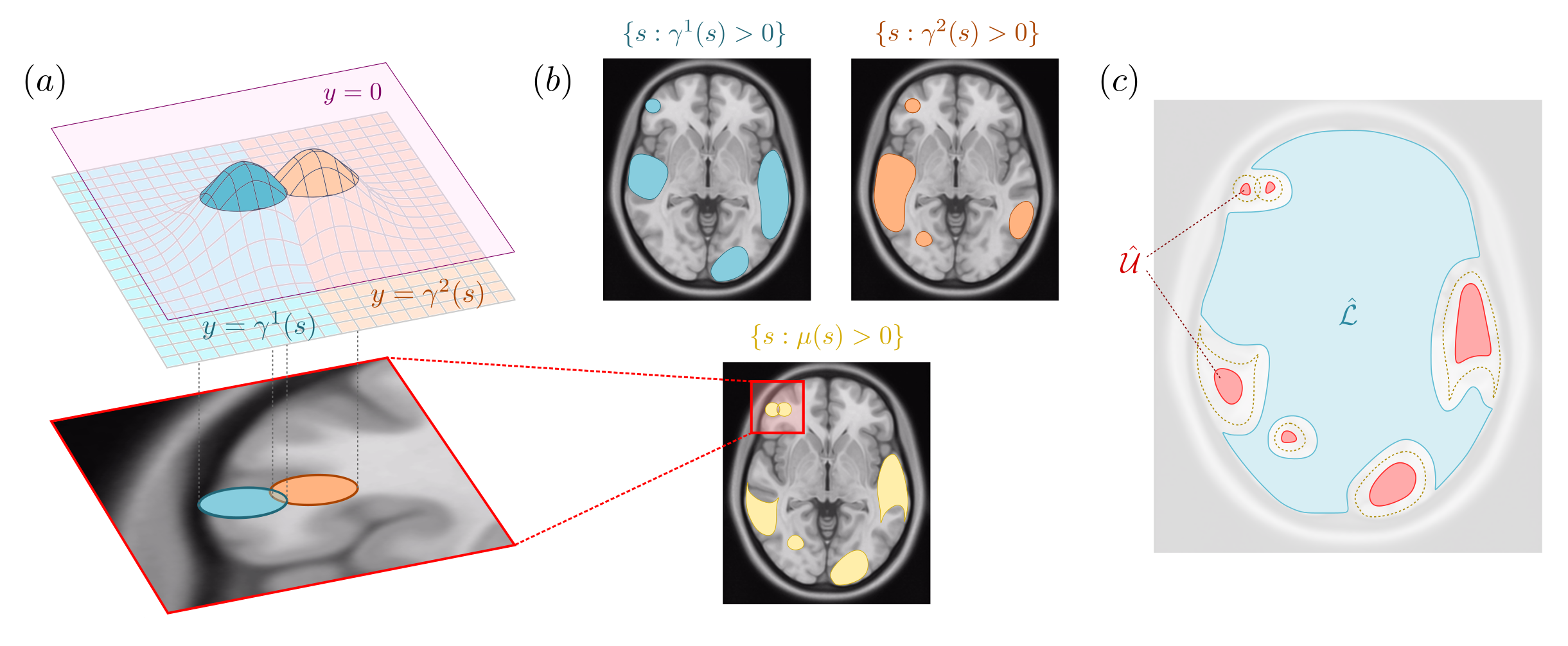}
    \captionof{figure}{$(a)$ Top: Two spatially varying signals, $\gamma^1$ (blue) and $\gamma^2 : \mathcal{S} \rightarrow \mathbb{R}$ (orange), representing the true unknown brain functions of healthy and control subjects, respectively, both thresholded at the level $c=0$. Bottom: The corresponding excursion sets, defined as $\{s \in \mathcal{S} : \gamma^i(s) > c\}$ for $i \in {1,2}$. $(b)$ Top: A broader view of the excursion sets $\{s \in \mathcal{S} : \gamma^i(s) > c\}$ for $i \in {1,2}$, overlaid on a structural image of the human brain. Bottom: The symmetric difference set, representing regions that belong to one excursion set but not both. $(c)$ CRs for the symmetric difference set, with the lower and upper bounds denoted by $\hat{\mathcal{L}}$ and $\hat{\mathcal{U}}$, respectively, and shown in blue and red shading.}
    \label{fig:symdiffX}
    \end{centering}
\end{example}
\end{examplebox} 

In this example, the researcher possesses two estimators, $\hat{\gamma}_n^1, \hat{\gamma}_n^2 : \Omega \rightarrow l^\infty(\mathcal{S})$, corresponding to the unknown true signals $\gamma^1, \gamma^2 : \mathcal{S} \rightarrow \mathbb{R}$, respectively. For these signals, we again assume a functional central limit theorem (fCLT) of the form:
\begin{equation}\nonumber
    \hat{J}_n^i:=\tau_n^{-1}(\hat{\gamma}_n^i-\gamma^i)\xrightarrow{d}J^i
\end{equation}
which holds jointly over $i \in \{1,2\}$, where the processes $\mathbf{\hat{J}}_n$ and $\mathbf{J}$ take values in $l^\infty_2(\mathcal{S})$ and $C(\mathcal{S}, \mathbb{R})^2$, respectively.

Our focus is on the symmetric difference between the excursion sets $\{s \in \mathcal{S} : \gamma^i(s) > c\}$ for $i \in \{1,2\}$. We begin by defining the \textit{symmetric difference of two functions} $\gamma^1$ and $\gamma^2$ as:
\begin{equation}\nonumber
    (\gamma^1 \Delta \gamma^2)(s) := \max\bigg(\min(-\gamma^1(s), \gamma^2(s)),\min(\gamma^1(s), -\gamma^2(s))\bigg).
\end{equation}
It is straightforward to verify that the symmetric difference between the sets $\{s \in \mathcal{S} : \gamma^i(s) > c\}$ for $i \in \{1,2\}$ is equivalent to the set $\{s \in \mathcal{S} : (\gamma^1 \Delta \gamma^2)(s) > c\}$. This is the set for which we aim to construct CRs. To this end, we treat $\hat\mu_n :=\hat{\gamma}_n^1 \Delta \hat{\gamma}_n^2$ as an estimator for $\gamma^1 \Delta \gamma^2$ and, as in previous sections, assume without loss of generality that $c = 0$.

Next, we define a partition of $\mathcal{S}$, $\{V_{\rho}\}$, indexed by $\rho\in\mathcal{I}:=\{-1,0,1\}^2$ as:
\begin{equation}\nonumber
    V_{\rho}:=\big\{s \in \mathcal{S}: \text{sgn}(\gamma^1(s)-\gamma^2(s)) = \rho_1~\text{ and }~\text{sgn}(\gamma^1(s)+\gamma^2(s)) = \rho_2\big\},
\end{equation}
where $\rho_k$ is the $k^{th}$ component of $\rho$ and sgn represents the sign function; i.e. sgn$(x)=x/|x|$ for $x\neq 0$ and sgn$(0)=0$ otherwise.
In this instance, our interest lies in the random process $\hat{H}_n=\tau_n^{-1}(\hat{\gamma}_n^1 \Delta \hat{\gamma}_n^2 - \gamma^1 \Delta \gamma^2)$, which, similarly to previous examples, does not satisfy a fCLT. In fact, the pointwise limit of $\hat{H}_n$ can be derived via direct computation as follows:\newpage
\begin{floatingfigure}[r]{0.56\textwidth}\centering
    \includegraphics[width=0.52\textwidth]{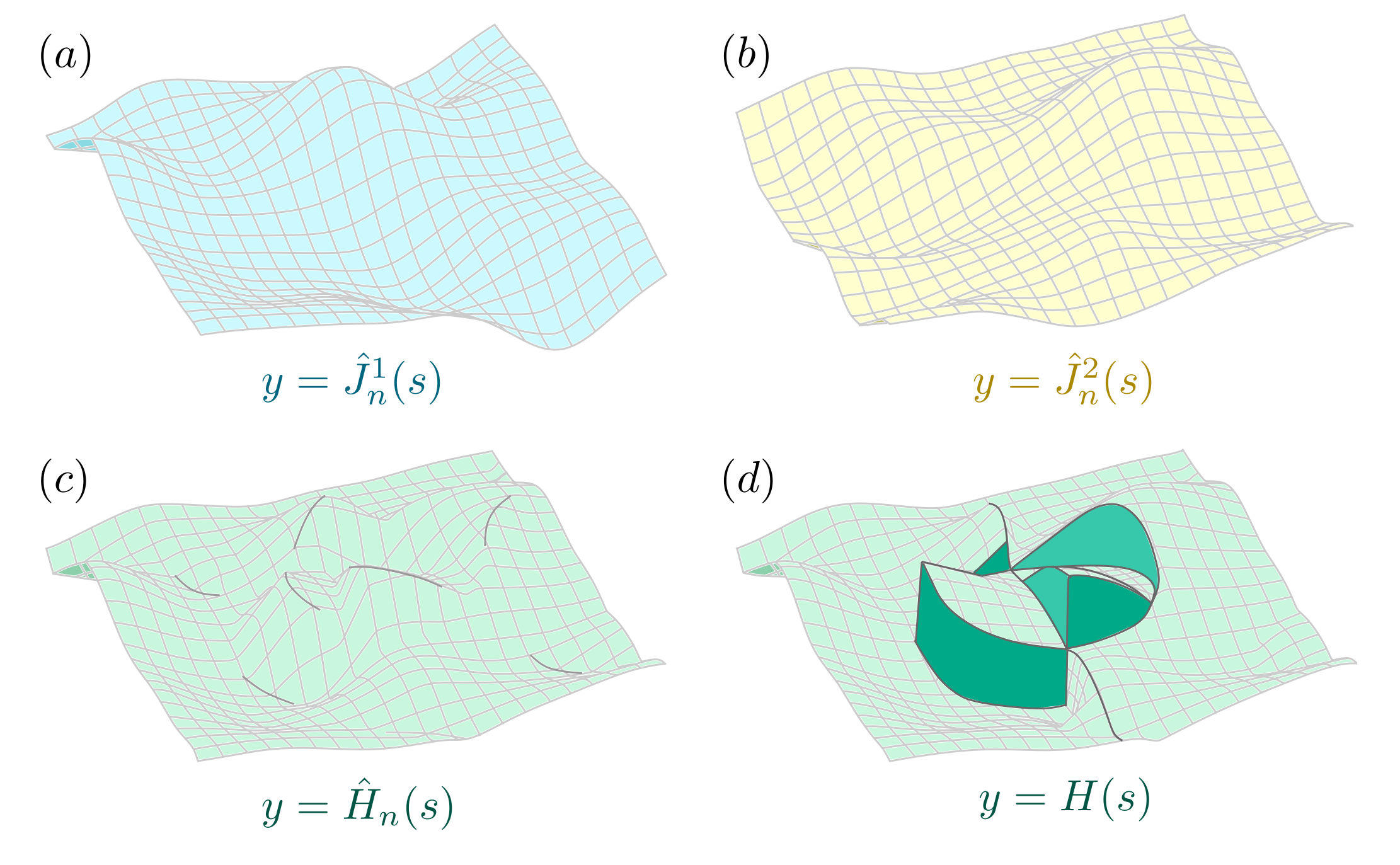}
    \captionof{figure}{Realizations of the random processes $(a)\ \hat{J}_n^1$ and $(b)\ \hat{J}_n^2$, alongside the corresponding process $(c)\ \hat{H}_n$. Shown in $(d)$ is the process $H$, the pointwise limit of $\hat{H}_n$ as $n \rightarrow \infty$. As $n$ increases, the `creases' in the realisations of $\hat{H}_n$ become steeper, until in the limiting case, $H$, they become directly vertical discontinuities. Note that, although depicted as such, we do not require that $\hat{J}_n^1$ and $\hat{J}_n^2$ be continuous processes.}
    \label{fig:symdiff3}
\end{floatingfigure}

\begin{equation}\nonumber
    H(s)=\scalebox{0.84}{$\displaystyle\begin{cases}
       -J_1(s) & \text{ if } s \in V_{(-1,+1)},\\
       -J_2(s) & \text{ if } s \in V_{(+1,+1)},\\
       J_1(s) & \text{ if } s \in V_{(+1,-1)},\\
       J_2(s) & \text{ if } s \in V_{(-1,-1)},\\
       \max(J_1,J_2)(s) & \text{ if } s \in V_{(0,-1)},\\
       \min(-J_1,J_2)(s) & \text{ if } s \in V_{(-1,0)},\\
       \min(J_1,-J_2)(s) & \text{ if } s \in V_{(1,0)},\\
       \max(-J_1,-J_2)(s) & \text{ if } s \in V_{(0,1)},\\
       (J_1\Delta J_2)(s) & \text{ if } s \in V_{(0,0)}.\\
    \end{cases}$}
\end{equation}

 Unfortunately, however, unlike in the previous examples, $\hat{H}_n$ does not converge with either upper or lower restraint to $H$. This fact is demonstrated by the following one dimensional example, Example \ref{example:symdiff4}. As such, for this example, instead of establishing an rCLT, we shall construct a confinement $(\hat{L}_n,\hat{U}_n)$ for $\hat{H}_n$ (see Definition \ref{def:conf_pair}). 

\begin{examplebox}
\begin{example}\label{example:symdiff4}
    Let $\mathcal{S}=[-2,2]$ and, for $i \in \{1,2\}$, define $\gamma^i$ and $\hat{\gamma}^i_n$ as:
    \begin{equation}\nonumber
        \gamma^1(s) := 2|s|, \quad \gamma^2(s):=s, \quad \hat{\gamma}^1_n(s):=\gamma^1(s)-\frac{2}{\sqrt{n}}, \quad \text{and}\quad \hat{\gamma}^2_n(s):=\gamma^2(s).
    \end{equation}
    It is clear that $\hat{\gamma}_n^i\xrightarrow{}\gamma^i$ uniformly. Defining $\hat{H}_n(s):=\sqrt{n}((\hat{\gamma}_n^1 \Delta \hat{\gamma}_n^2)-(\gamma^1 \Delta \gamma^2))$ and $H(s):=0$, we again have that $\hat{H}_n\rightarrow H$ pointwise (see Fig. \ref{fig:symdiff4}). However, for all $n$, we have that $\sup_{s \in [0,2/\sqrt{n}]}|\hat{H}_n(s)|=5/3>1$ and thus $\hat{H}_n$ cannot converge with upper or lower restraint to $H$ as $|H(s)|=0\leq 1$ for all $s$.
    
\centering
    \includegraphics[width=0.55\textwidth]{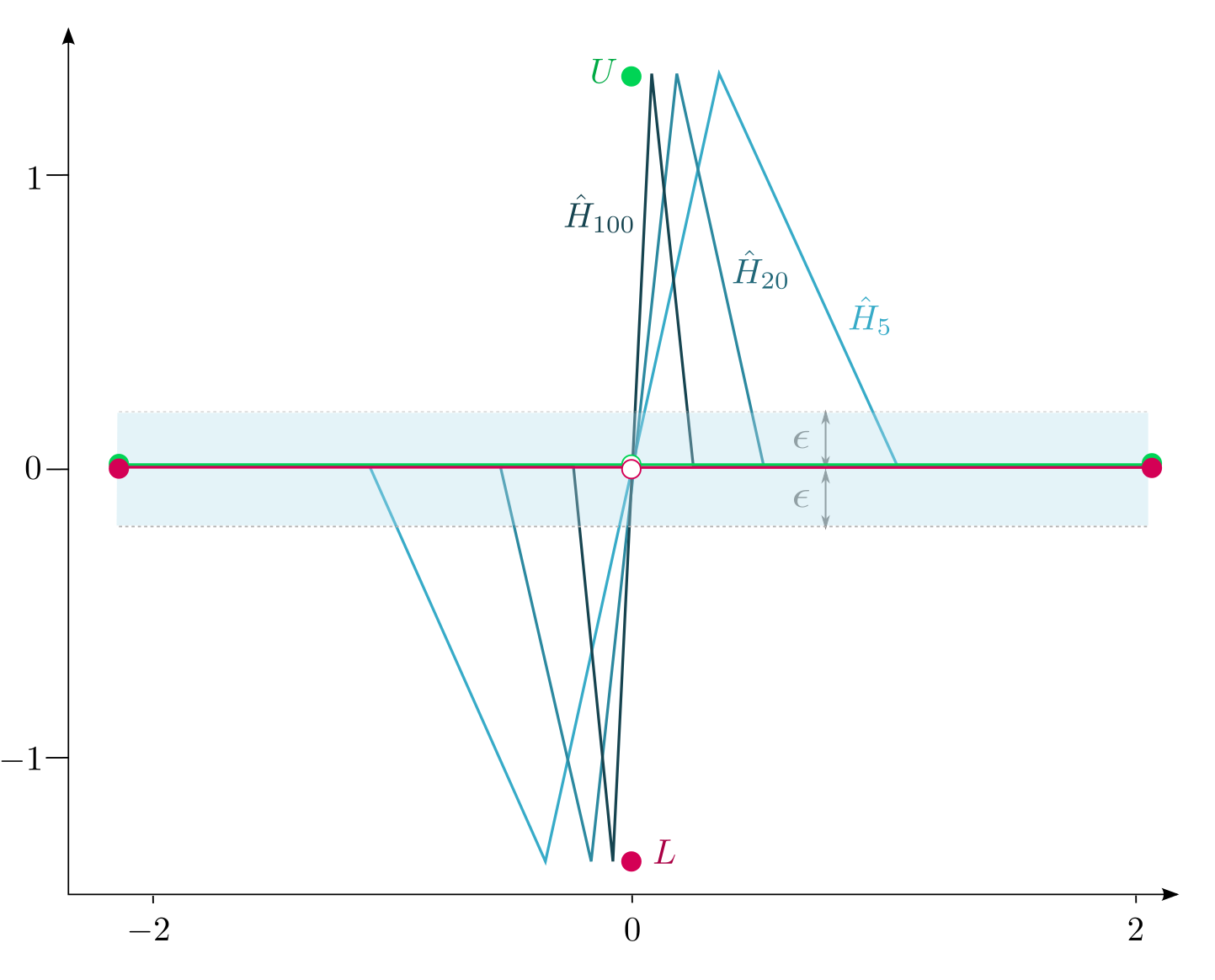}
    \captionof{figure}{An illustration of $H$ and $\{\hat{H}_n\}$ for $n=5,20$ and $100$ in Example \ref{example:symdiff4}. If it were the case that $\hat{H}_n\rightarrow H$ upper (lower) restraint, $\hat{H}_n$ would eventually have to lie below the upper (lower) boundary of the shaded blue region. However, in this example, for any $n$, $\hat{H}_n$ will always dip both below and above the shaded blue region near $s=0$.}
    \label{fig:symdiff4}
\end{example}
    \vspace*{0.25cm}
\end{examplebox}

 To generate CRs for the symmetric difference, we must account for and isolate the type of \textit{``unrestrained''} behaviour exhibited near the origin in Fig. \ref{fig:symdiff4}. 
To do so, it is first worth noting that the problem becomes substantially more tractable if we employ the following change of variables:
\begin{equation}\nonumber
    (\gamma^1,\gamma^2) \mapsto (d, m):= \bigg(\frac{1}{2}(\gamma^1-\gamma^2),\frac{1}{2}(\gamma^1+\gamma^2)\bigg).
\end{equation}
Under this transformation, Supplementary Lemma \ref{lem:symdiff} shows that
 \begin{equation}\nonumber
    \hat H_n(s) = \tau_n^{-1} \big( \vert \hat d_n \vert  - \vert d \vert \big) - \tau_n^{-1} \big(\vert m \vert  - \vert \hat m_n \vert   \big),
\end{equation}
which is readily seen to be the sum of two random processes, each satisfying an rCLT by Theorem~\ref{thm:abs_val_res}. Although, as previously noted, the fact that $\hat H_n$ is a sum of two processes satisfying rCLTs does not necessarily imply that it satisfies an rCLT itself (cf. Example~\ref{example:symdiff4}), this representation is nonetheless useful as it enables the application of Corollary~\ref{cor:SumSupPreserving}. In particular, we can use Corollary~\ref{cor:SumSupPreserving} to identify the problematic subset $\mathcal{N}$ over which a restrained bound is not guaranteed. As shown in Supplementary Material Section \ref{supp:symdiff_conv_with_res}, in the context of the symmetric difference, this subset is given by
\begin{equation}\nonumber
    \mathcal{N}:= \bigcup_{\kappa \in \{-1,+1\}^2}V_{(0,0)} \cap \overline{V_\kappa}.
\end{equation}
Having identified the region over which convergence is not guaranteed to hold, we now construct a confinement $(\hat{L}_n,\hat{U}_n)$ for $\hat{H}_n$ (c.f. Definition \ref{def:conf_pair}). Specifically, we choose $\hat{L}_n$ and $\hat{U}_n$ such that they coincide with $\hat{H}_n$ on $\mathcal{S} \setminus \mathcal{N}$, but act as lower and upper bounds for $\hat{H}_n$ on $\mathcal{N}$. 

Revisiting the definitions of $H$ and $\hat{H}_n$, it follows that, within $V_{(0,0)}$, we have $|H|\leq \max(|J^1|,|J^2|)$ and $|\hat{H}_n|\leq \max(|\hat{J}_n^1|,|\hat{J}_n^2|)$. We shall use these bounds to ``patch over" the ``difficult region'' $\mathcal{N}$ by defining $\hat{U}_n(s), U(s),\hat{L}_n(s)$ and $L(s)$ as follows:
\begin{equation}\nonumber
    \hat{U}_n(s):=\begin{cases}
        \hat{H}_n(s) & s\not\in \mathcal{N} \\
        \max\big(|\hat{J}_n^1(s)|,|\hat{J}_n^2(s)|\big) & s\in \mathcal{N} \\
    \end{cases} \quad\text{ and }\quad  U(s):=\begin{cases}
        H(s) & s\not\in \mathcal{N} \\
        \max\big(|J^1(s)|,|J^2(s)|\big) & s\in \mathcal{N},
    \end{cases}
\end{equation}
\begin{equation}\nonumber
    \hat{L}_n(s):=\begin{cases}
        \hat{H}_n(s) & s\not\in \mathcal{N} \\
        -\max\big(|\hat{J}_n^1(s)|,|\hat{J}_n^2(s)|\big) & s\in \mathcal{N} \\
    \end{cases} \quad\text{ and }\quad  L(s):=\begin{cases}
        H(s) & s\not\in \mathcal{N} \\
        -\max\big(|J^1(s)|,|J^2(s)|\big) & s\in \mathcal{N}.
    \end{cases}
\end{equation}
An immediate consequence of the above definitions is that $\hat{U}_n$, $U$, $\hat{L}_n$ and $L$ are not piecewise continuous over $\{V_\kappa\}_{\kappa \in \{-1,+1\}^2}$ (as we have partitioned $V_{(0,0)}$ into two). For this reason we must modify our partition to obtain $\{V^\ast_\kappa\}_{\kappa \in \mathcal{I^*}}$ where $\mathcal{I}^\ast:=\mathcal{I}\cup \{\ast\}$ and:
\begin{equation}\nonumber
    V^\ast_\kappa := \begin{cases}
        V_\kappa, & \text{if } \kappa \neq (0,0)\text{ and }\kappa\neq \ast,\\
        V_{(0,0)}\setminus \mathcal{N}, &\text{if } \kappa=(0,0),\\
        \mathcal{N}, & \text{if }\kappa=\ast.
        
    \end{cases}    
\end{equation}
We now are in a position to state the following convergence result.

\begin{mdframed}
\begin{theorem}[PCLT with Confinement for the Symmetric Difference]\label{thm:symdiff_conv_with_res} 
    Given the definitions above, $\hat{\mu}_n= \hat{\gamma}_n^1 \Delta \hat{\gamma}_n^2$ satisfies a pointwise CLT with confinement $(\hat{L}_n,\hat{U}_n)$.
\end{theorem}
\vspace*{0.3cm}
\end{mdframed}
\begin{proof}
    See Supplementary Material Section \ref{supp:symdiff_conv_with_res}.
\end{proof}
We can now generate CRs as follows:
\begin{mdframed}
    \begin{theorem}[Confidence Regions for the Symmetric Difference Example]\label{thm:crs_symdiff}
        Suppose $\mu,\gamma^i:\mathcal{S}\rightarrow \mathbb{R}$, $\hat{\gamma}_n^i,\hat{J}^i_n:\Omega\rightarrow l^\infty(\mathcal{S})$, and $J^i:\Omega\rightarrow C(\mathcal{S},\mathbb{R})$ are defined as above and $\mathcal{L}, \hat{\mathcal{L}}, \mathcal{U}$ and $ \hat{\mathcal{U}}$ are given by:
\begin{equation}\nonumber
    \mathcal{L}:=\bigg\{s \in \mathcal{S}:(\gamma^1 \Delta \gamma^2)(s)<0\bigg\},\quad\mathcal{U}:=\bigg\{s \in \mathcal{S}:(\gamma^1 \Delta \gamma^2)(s)>0\bigg\},
\end{equation}
\begin{equation}\nonumber
    \hat{\mathcal{L}}:=\bigg\{s \in \mathcal{S}:\tau_n^{-1}(\hat{\gamma}_n^1 \Delta \hat{\gamma}_n^2)(s)<-q\bigg\}\quad \text{ and }\quad\hat{\mathcal{U}}:=\bigg\{s \in \mathcal{S}:\tau_n^{-1}(\hat{\gamma}_n^1 \Delta \hat{\gamma}_n^2)(s)>+q\bigg\}.
\end{equation}
Further assume $\mathfrak{u}^\pm_\sigma(0)$ is defined for $\sigma\in\mathcal{I}$ as in the statement of Theorem \ref{thm:scope_for_piecewise}, and $\mathcal{N}$ is defined as above. Defining the set $T_i^k$ for $i\in\{1,2\}$ and $k\in\{-1,1\}$ as:
\begin{equation}\nonumber
    T_i^k:= \mathfrak{u}^-_{((3-2i)k,-k)}(0) \cup \mathfrak{u}^+_{((2i-3)k,k)} (0)\cup \mathfrak{u}^-_{(0,-k)}(0)\cup \mathfrak{u}^+_{(-k,0)}(0)\cup \mathcal{N},  \quad \text{and}
\end{equation}
\begin{equation}\nonumber
    R_{1} := \mu^{-1}(0)\cap \big(V_{(-1,1)}~\cup~V_{(1,-1)}\big) \quad \text{and} \quad R_{2} := \mu^{-1}(0)\cap \big(V_{(1,1)}~\cup~V_{(-1,-1)}\big),
\end{equation}
we now obtain that:
    \begin{equation}\nonumber
			        \begin{split}
			            (i)\quad \liminf_{n\rightarrow\infty}
			            	\mathbb{P}_*\big[\,& \hat{\mathcal{L}} \subseteq \mathcal{L}~ \wedge ~\hat{\mathcal{U}} \subseteq \mathcal{U}\,\big]
			         \geq \\
                        & \mathbb{P}\Bigg[\max\Bigg(\max_{\substack{\text{\scriptsize{$i \in \{1,2\}$}}\\
                        \text{\scriptsize{$k \in\{-1,1\}$}}}}\sup_{s \in T_i^k}kJ^i(s), \sup_{s \in V_{(0,0)}\setminus \mathcal{N}}|(J^1\Delta J^2)(s)|\Bigg)<q\Bigg],
			                    \end{split}
			        \end{equation}                    
            	    \begin{equation}\nonumber  
	        		\begin{split}
           				(ii)\quad  \limsup_{n\rightarrow\infty}~\mathbb{P}^*\big[\,&\hat{\mathcal{L}} \subseteq \mathcal{L}~ \wedge ~\hat{\mathcal{U}} \subseteq \mathcal{U}\,\big]
           				\leq \\ & \mathbb{P}\bigg[\max\bigg(\max_{i\in \{1,2\}}\sup_{s\in R_i}|J^i(s)|,\sup_{s\in V_{(0,0)}}|(J^1\Delta J^2)(s)|\bigg) \leq q\bigg].
           			\end{split}
        			\end{equation}
                    \vspace{0.2cm}
    \end{theorem}
\end{mdframed}
\begin{proof}
    See Supplementary Material Section \ref{supp:crs_symdiff}
\end{proof}

As with the previous examples, we now consider additional assumptions which may simplify the statement of Theorem \ref{thm:crs_symdiff}, allowing for an exact limit, rather than limit inferiors or superiors. On first viewing, one might consider the assumption $\mathcal{N} = \emptyset$ a natural candidate.  However, despite greatly simplifying the expressions in Theorem \ref{thm:crs_symdiff}, this assumption is in practice unrealistic. For example, even in the simple `venn-diagram' arrangement shown on the bottom of panel (a) in Figure \ref{fig:symdiffX}, the two points at the top and bottom of the central oval can be shown to lie in $\mathcal{N}$. Given this, we instead propose the following alternative conditions.

\begin{mdframed}
    \begin{corollary}\label{corr:conds_sd}
        Suppose the conditions of Theorem \ref{thm:abs_val_crs} and one of the following hold:
        \begin{itemize}
            \item[\textbf{S1:}] For every $\kappa \in \mathcal{I}$, we have that $\overline{\mu^{-1}(0)\cap V_\kappa}=\mu^{-1}(0)\cap \overline{V_\kappa}$.
            \item[$\widetilde{\textbf{S1}}$\textbf{:}] For every $s\in \mu^{-1}(0)$, if $s \in \overline{V_\kappa}$, then for any sufficiently small open ball around $s$, $B_\delta(s)$ with $\delta>0$, we have that the sets $B_\delta(s)\cap \mu^{-1}([0,\infty))\cap V_\kappa$ and $B_\delta(s)\cap \mu^{-1}((-\infty,0])\cap V_\kappa$ are non-empty, and the set $B_\delta(s)\cap V_\kappa$ is path-connected.
        \end{itemize}
        If, in addition, each of the below conditions also hold:
        \begin{itemize}
            \item[\textbf{S2:}] $\mathbb{P}\bigg[\displaystyle\sup_{s \in \mu^{-1}(0)} -L(s) = \displaystyle\sup_{s \in \mu^{-1}(0)} -H(s) \bigg] = 1$ and  $\mathbb{P}\bigg[\displaystyle\sup_{s \in \mu^{-1}(0)} U(s) = \displaystyle\sup_{s \in \mu^{-1}(0)} H(s) \bigg] = 1$,
            \item[\textbf{S3:}] The below probability is zero:
            \begin{equation}\nonumber
                \mathbb{P}\bigg[\max\bigg(\max_{i\in \{1,2\}}\sup_{s\in R_i}|J^i(s)|,\sup_{s\in V_{(0,0)}}|(J^1\Delta J^2)(s)|\bigg) = q\bigg]=0,
            \end{equation}
        \end{itemize}
        then we have that:
        \begin{equation}\nonumber
            \lim_{n\rightarrow\infty}~\mathbb{P}\big[\,\hat{\mathcal{L}} \subseteq \mathcal{L}~ \wedge ~\hat{\mathcal{U}} \subseteq \mathcal{U}\,\big]
           				=   \mathbb{P}\bigg[\max\bigg(\max_{i\in \{1,2\}}\sup_{s\in R_i}|J^i(s)|,\sup_{s\in V_{(0,0)}}|(J^1\Delta J^2)(s)|\bigg) \leq q\bigg].
        \end{equation}
    \end{corollary}
\end{mdframed}

\begin{proof}
    See Supplementary Material Section \ref{supp:conds_sd}.
\end{proof}

\begin{remark}
    
	In words, the above conditions say the following. Condition \textbf{S1} handles degenerate situations in which the boundaries of $\{V_\kappa\}$ happen to be tangential to $\mu^{-1}(0)$, or exhibit unusual fractal like behaviour. Condition \textbf{S2} acknowledges the possibility of the spiking behaviour shown in Fig. \ref{fig:symdiff4}, but assumes that such behaviour affects the computation of the supremum with probability zero. This assumption is partly motivated by the fact that, while it is possible to construct symmetric difference examples such as Example \ref{example:symdiff4} that exhibit erratic, non-convergent spiking, doing so appears to be difficult. We conjecture that in many practical settings, such behaviour is confined to small regions of the spatial domain and thus has minimal impact on the overall computation. Further work is needed to rigorously verify this claim. Finally, condition \textbf{S3} ensures that the distribution of the maximum statistic does not have a discrete component. For further discussion of the ball assumption, assumption $\widetilde{\textbf{S1}}$, see the remark following Corollary \ref{cor:conj_conditions}. 
\end{remark}

%% file: Article/discussion.tex
The development of valid inferential tools is fundamental in statistical science.
Although CRs for excursion sets have recently gained popularity, the current available literature often relies on restrictive assumptions such as that of a continuous error processes. We resolved this gap by introducing a novel convergence concept for piecewise-continuous processes, specifically designed to establish asymptotic results on quantile convergence for functionals involving the supremum and infimum of a random process. Among others, this new framework allowed us to eliminate the stringent condition of differentiability and the ball condition imposed in the previous results of \cite{MaullinSapey2022}. Furthermore, it facilitates the construction of CRs in a range of new, previously unattainable, contexts, such as the symmetric difference of excursion sets. 

Our developments allow for a wider range of applications of CRs of excursion sets; this is important as the question of spatial inference on logical combinations of excursion sets arises naturally in applications such as neuroimaging,  climatology, and cosmology and other branches involving spatial patterns. As outlined in Section \ref{sec:applications}, we believe the method can be utilised to answer a wide range of practical questions concerning the spatial similarity of excursion sets in imaging data. We intend to pursue these applications in greater depth in future work. 

Additional future work will also explore the theory and applications of our novel convergence notion in greater detail. 
Preliminary investigations reveal that convergence with restraint appears to be closely related to convergence in Skorokhod metrics, with the key distinction that the former does not require a designated direction from which limits must exist, at each point $s \in \mathcal{S}$. Another direction for future development is to extend our results to allow for locally countable partitions (i.e.~to the setting where, for all $s\in \mathcal{S}$, there is a neighbourhood of $s$ which intersects potentially countably infinitely many $V_i$). Such an extension requires additional technical assumptions that go beyond the scope of this paper. We note here that these technical challenges also pose an issue for extensions of Skorokhod metrics to higher-dimensional domains, which similarly require piecewise continuous functions be defined on locally countable partitions. Another challenging yet important avenue for future research is to investigate whether convergence with restraint can be adapted to define a metric and to explore whether in such a framework an Arzelà–Ascoli type theorem can be proven in order to establish a natural weak convergence theory for piecewise continuous functions over arbitrary dimensional domains.

%% file: Article/acknowledgement.tex
F.T. is funded by the Deutsche Forschungsgemeinschaft (DFG) under Excellence Strategy The
Berlin Mathematics Research Center MATH+ (EXC-2046/1, project ID:390685689)